# A Fast Differential Grouping Algorithm for Large Scale Black-Box Optimization

Zhigang Ren, *Member, IEEE*, An Chen, Yaochu Jin, *Fellow, IEEE*, Wenhua Guo, Yongsheng Liang, *Student Member, IEEE*, and Zuren Feng, *Member, IEEE*

*Abstract*—Decomposition plays a significant role in cooperative co-evolution which shows great potential in large scale black-box optimization. However, current popular decomposition algorithms generally require to sample and evaluate a large number of solutions for interdependency detection, which is very time-consuming. To address this issue, this study proposes a new decomposition algorithm named fast differential grouping (FDG). FDG first identifies the type of an instance by detecting the interdependencies of a few pairs of variable subsets selected according to certain rules, and thus can rapidly complete the decomposition of a fully separable or nonseparable instance. For an identified partially separable instance, FDG converts the key decomposition process into a search process in a binary tree by taking corresponding variable subsets as tree nodes. This enables it to directly deduce the interdependency related to a child node by reutilizing the solutions sampled for corresponding parent and brother nodes. To support the above operations, this study designs a normalized variable-subset-oriented interdependency indicator, which can adaptively generate decomposition thresholds according to its distribution and thus enhances decomposition accuracy. Computational complexity analysis and experimental results verify that FDG outperforms popular decomposition algorithms. Further tests indicate that FDG embedded in a cooperative co-evolution framework can achieve highly competitive optimization results as compared with some state-of-the-art algorithms for large scale black-box optimization.

## I. INTRODUCTION

WITH the enhancement of people's ability to acquire, process, and analyze data, more and more large scale

Manuscript received March XX, 2019; revised XXXX XX, 20XX; accepted XXXX XX, 20XX. This work was supported in part by the National Nature Science Foundation of China under Grant 61873199 and in part by China Postdoctoral Science Foundation under Grants 2014M560784 and 2016T90922. (*Corresponding Authors: Zhigang Ren and Yaochu Jin*)

Z. Ren, A. Chen, and Y. Liang are with the Autocontrol Institute, School of Electronic and Information Engineering, Xi'an Jiaotong University, Xi'an 710049, China (e-mail: renzg@mail.xjtu.edu.cn; chenan123@stu.xjtu.edu.cn; liangyongsheng@stu.xjtu.edu.cn).

Y. Jin is with the Department of Computer Science, University of Surrey, Guildford GU2 7XH, U.K. (e-mail: yaochu.jin@surrey.ac.uk).

W. Guo is with the State Key Laboratory for Manufacturing Systems Engineering, School of Mechanical Engineering, Xi'an Jiaotong University, Xi'an 710049, China (e-mail: markguo@mail.xjtu.edu.cn).

Z. Feng is with the State Key Laboratory for Manufacturing Systems Engineering, School of Electronic and Information Engineering, Xi'an Jiaotong University, Xi'an 710049, China (e-mail: fzr9910@mail.xjtu.edu.cn).

black-box optimization (LSBO) problems have been emerging in scientific research and engineering applications since the last decade [1]-[5]. For example, some aerodynamic shape optimization problems may consist of thousands of decision variables [6], and the reconstruction of a gene regulatory network in bioinformatics needs to determine the status of thousands of interactions [7]. However, mathematical programming methods, which strictly depend on problem models, are unsuited for such problems due to their black-box nature; evolutionary algorithms (EAs), which do not require analytical objective functions, still work, but tend to lose their efficiency as the problem dimension increases [1]-[3]. This can be attributed to the fact that the solution space of a problem exponentially grows with the increase of its dimension and the traditional EAs cannot adequately explore the solution space of a LSBO problem within acceptable computation time [8], [9].

To address this issue, some researchers developed cooperative co-evolution (CC) methods [5], [8]-[13]. Taking *"divide-and-conquer"* as the basic idea, CC solves a LSBO problem by first dividing it into a set of smaller and simpler subproblems and then cooperatively optimizing them with a traditional EA. Since its inception, CC has attracted much research attention and achieved great progress in algorithmic components, including the decomposition method, the optimizer for subproblems, and the computation resource allocation strategy [1], [3], [5], [14]. It has also been shown that CC is really superior to the traditional EAs in tackling LSBO problems [14]-[17]. Nevertheless, some key challenges, such as the development of effective and efficient decomposition methods, remain open.

It is understandable that decomposition is fundamental to CC. A proper decomposition can greatly reduce the optimization difficulty of a LSBO problem without changing its optimum. However, an improper decomposition may cause CC to converge to a Nash equilibrium rather than a real optimum [18]. Static decomposition methods [8]-[11], which are mainly employed by early CC algorithms, directly divide decision variables into several fixed subcomponents without considering their interdependencies at all. Consequently, they perform poorly on nonseparable and partially separable problems. Random decomposition methods remedy this defect to a certain extent by stochastically grouping variables and dynamically changing the grouping result in every cycle of CC [12], but their performance still deteriorates on LSBO problems that have many interacting variables [19].

Different from the two classes of decomposition methods listed above, learning-based methods conduct decomposition by explicitly investigating the interdependencies among variables, and thus significantly improve the decomposition accuracy at the cost of sampling and evaluating a certain number of solutions [13], [20], [21]. As a representative of this type of method, differential grouping (DG) provides a clear and simple criterion for decomposition and can outperform some other learning-based decomposition algorithms such as variable interaction learning (VIL) [13], [20]. Due to its advantage and potential, the basic DG has attracted a great deal of research effort, and several variants have been suggested, including global DG (GDG) [22], extended DG (XDG) [23], graph-based DG (gDG) [24], fast interdependency identification (FII) [25], DG2 [26], and recursive DG (RDG) [17]. Among these variants, RDG, the most recently developed one, provides a reference. With a recursive decomposition strategy, it reduces the fitness evaluation (FE) requirements of GDG and DG2 from $O(n^2)$ to $O(n \log_2 n)$ for an $n$-dimensional problem. This result seems exciting, but is still unacceptable for many real-world LSBO problems whose objective functions cannot be analytically described and can be evaluated only by numerical simulations, which are very expensive in terms of both time and economy [16], [27]. Therefore, how to further reduce the FE requirement becomes a challenging and realistic task. On the other side, most learning-based decomposition algorithms demand users to specify a threshold for the decomposition indicator. However, the threshold is generally instance-dependent or even variable-dependent, so is hard to set [13], [22], [26]. DG2 alleviates this issue by adaptively determining a threshold for each pair of variables [26]. Nevertheless, it requires the full interdependency information in advance so as to determine the weights of the lower and the upper bounds of the computational roundoff error in the corresponding threshold.

This study aims to address the above issues by proposing a fast DG (FDG). The main contributions of this work are:

1) It designs a normalized variable-subset-oriented interdependency indicator accompanied with a simple but efficient indicator distribution analysis method. This enables FDG to perform decomposition from the perspective of variable subsets and to detect interdependencies according to the indicator distribution without a user specified threshold. As a consequence, FDG needs to calculate much fewer indicators and can be adapted to different LSBO instances well.

2) Different from existing decomposition algorithms which handle different types of LSBO instances in the same manner, FDG first identifies the type of separability of an instance by detecting the interdependencies of a few pairs of variable subsets selected according to a sophisticated strategy, and then further performs decomposition only if the instance is partially separable. This implies that the fully separable or nonseparable instances can be tackled with a very small number of FEs.

3) FDG converts the key decomposition process on a partially separable instance into a search process in a binary tree by taking corresponding variable subsets as nodes, which facilitates the reutilization of the evaluated solutions. As a result, FDG can directly deduce the interdependency related to a child node by reusing the solutions sampled for the corresponding parent and brother nodes.

This study also theoretically analyzes the computational complexity of FDG on different types of LSBO instances and evaluates its performance on two benchmark suites. The results indicate that FDG is superior to existing DG variants in terms of both efficiency and accuracy.

The remainder of this paper is organized as follows. Section II reviews the related work on CC. Section III presents FDG in detail. Section IV reports experimental settings and results. Finally, Section V concludes the paper and discusses some promising directions for the future research.

## II. RELATED WORK

CC has been successfully adopted to tackle a variety of LSBO problems. Algorithm 1 presents the general framework of CC, where three main algorithmic components are involved, including the decomposition method, the optimizer, and the computation resource allocation strategy.

The main challenge in applying CC consists in problem decomposition. Early CC algorithms mainly aim at enhancing the performance of conventional EAs on ordinary problems instead of scaling them up to large scale problems, and generally employ simple static decomposition methods. This class of methods divide an $n$-dimensional problem into $k$ $s$-dimensional subproblems with $s < n$ and $ks = n$ and keep the decomposition result fixed during the whole optimization process. Examples of such methods include the "$n$ to $n \times 1$" method [10], the splitting-into-half [11] stategy, and the general "$n$ to $k \times s$" method [9]. These methods really improve the performance of some EAs, such as evolutionary programming [8], particle swarm optimization (PSO) [9], genetic algorithm [10], and differential evolution [28], on certain kinds of problems. However, they lose their effectiveness on non-separable and partially separable problems due to the neglect of variable interdependence.

Directing against the rigidness of static decomposition methods, Yang *et al.* [12] proposed a random grouping (RG) scheme. RG stochastically selects decision variables for each subproblem in every cycle of CC with the aim of grouping any two interdependent variables into the same subproblem at least once. Omidvar *et al.* [19] showed that the probability of putting all the interdependent variables together sharply decreases with the increase of the number of interdependent variables, and suggested increasing the grouping frequency by reducing the iteration times in a cycle. Moreover, it was revealed that RG also faces the difficulty of getting a proper subcomponent size.

---

**Algorithm 1**: General CC

1. Decompose original decision vector $x : x \rightarrow \{x_1, \dots, x_k\}$ ; // Decomposition
2. Initialize the population for each subproblem $i$: $P_i$, $i = 1, \dots, k$ ;
3. Initialize the best overall solution $x^*$ ;
4. **while** the termination condition is not met **do**
5.    Determine the subproblem $i$ to be optimized;      // Allocation of CR
6.    $(P_i, x^*) \leftarrow$ Optimizer$(P_i, x^*, feNum)$ ;      // Optimizer
7.    Redivide $x$ and adjust populations of new subproblems if necessary;
8. **return** $x^*$ .

Aiming at this issue, Yang *et al.* [29] and Omidvar *et al.* [30] developed their respective adaptive methods, both of which probabilistically select a subcomponent size for every new cycle from a set of candidates based on the idea of reinforcement learning. Instead of using a sophisticated adaptive mechanism, Omidvar *et al.* [19] suggested randomly selecting a subcomponent size from its candidates if no fitness improvement is yielded during the previous cyle, and otherwise leaving it unchanged. This simpler strategy was also employed by Li and Yao [31] when scaling up PSO with CC and RG. Different from above subcomponent size adaption approaches which adjust the subcomponent size in a serial way, the Cooperative Coevolution with Adaptive Subcomponents proposed by Trunfio [32] concurrently applies a pool of subcomponent sizes at a specified learning phase such that a more reliable evaluation of the candidate sizes can be achieved, and finally just employies the best performing subcomponent size during the subsequent optimization phase. Trunfio *et al.* [33] further expanded this idea to choose a proper population size for each subproblem besides the subcomponent size.

Despite their simplicity, the static and the random decomposition methods can hardly get a desirable decomposition since they do not detect the interdependencies among decision variables. To remedy this defect, learning-based decomposition methods were proposed. Some researchers suggested taking all the variables as random ones and conducting decomposition according to the correlation coefficience between each pair of them [34], [35]. Xu *et al.* [36] indicated that the commonly-used Pearson correlation coefficient cannot properly depict the nonlinear interdependencies among decision variables, and replaced it with mutual information when carrying out grouping. Sun *et al.* [21] proved that, for a differentiable function $f(\cdot)$, if its two variables $x_i$ and $x_j$ directly interact, then $\partial f / \partial x_i$ has a functional relationship with $x_j$, and they further employed the maximal information coefficient between $\partial f / \partial x_i$ and $x_j$ to quantify the interdependency between the two variables. K. Weicker and N. Weicker [37] deemed that there exists an interdependency between two variables if a solution of higher quality can be generated by concurrently changing both variables than by just changing one of them. Based on this idea, they proposed the VIL method. Chen *et al.* [20], Sun *et al.* [38], and Ge *et al.* [39] further improved this method and made it more suitable for large scale problems. To reduce the FE requirement of this kind of methods, Ge *et al.* [40] recently developed a two stage variable interaction reconstruction algorithm, which first explores part of variable interactions using the fast VIL algorithm presented in [39], and then trains a marginalized denoising model to construct the overall variable interactions. Instead of investigating the interdependency between two variables, some researchers defined heuristic information for each variable, such as the variation of each variable in two consecutive CC cycles [41] and the correlation coefficient between each variable and the objective function [42], and grouped the variables of similar heuristic values into the same subcomponent.

Omidvar *et al.* [13] theorized the original interdependency discriminate criterion proposed by Tezuka *et al.* [43] and developed the effective and well-known DG method. Its main idea is that two variables are interdependent if the fitness variation caused by the perturbation on one variable relies on the value of the other. It was verified that DG is superior to the VIL algorithm developed by Chen *et al.* [20], but shows quite low accuracy on some benchmark functions due to the neglect of indirect interdependencies. To tackle this issue, Mei *et al.* [22], Sun *et al.* [23], Ling *et al.* [24], and Omidvar *et al.* [26] developed GDG, XDG, gDG, and DG2, respectively. These four DG variants tend to detect all pairwise interdependencies and thus requires $O(n^2)$ FEs for an $n$-dimensional problem, where gDG reduces its FE requirement to some extent by avoiding sampling solutions for any two variables having been verified to be indirectly interacting and DG2 achieves this by systematically reutilizing a few part of evaluated solutions. To further save FEs, Hu *et al.* [25] proposed the FII algorithm. FII first excludes separable variables from nonseparable ones by detecting the interdependency between each variable and all the other ones, and then subdivides nonseparable variables also without learning the full interdependency information among nonseparable variables. Much more recently, Sun *et al.* [17] proposed the RDG algorithm and achieved more excellent performance. For a variable subset $X_j$ interacting with the current nonseparable subset $X_i$, RDG first dichotomizes it and then investigates the interdependency between each of the resulting variable subset and $X_i$. This process is recursively conducted until all the variables interacting with $X_i$ are identified. In this way, RDG reduces its FE requirement to $O(n \log_2 n)$ for an $n$-dimensional problem.

Besides the detection mode, the decomposition threshold also affects much on the decompostion performance. As the original DG, XDG, gDG, and FII ask users to directly specify a threshold. This is a challenging task because the performance of these algorithms is very sensitive to the threshold [22], [26]. GDG and RDG weaken this sensitivity by setting the threshold to a value proportional to the minimum fitness value of several randomly generated solutions. However, this setting strategy still requires users to specify a value for the proportion coefficient and the resulting global threshold performs poorly on imbalanced problems [26]. To alleviate this issue, DG2 adaptively calculates a threshold for each pair of variables based on the greatest lower bound and the least upper bound on the corresponding computational roundoff error [26]. Sun *et al.* [44] recently simplified this threshold setting strategy and introduced it into RDG.

At present, decomposition remains one of the most important research topics of CC. Interested readers may refer to review papers [1], [3], and [5].

## III. Fast Differential Grouping

### A. Separability of LSBO Problems

Many real-world LSBO problems are difficult to tackle but possess an appealing characteristic, i.e., separability, where additive separability is the most common type and so is most extensively studied in the CC research field [13], [22]-[27]. The definition of additive separability can be described as follows.

**Definition 1.** *An n-dimensional function $f(\cdot)$ is said to be* additively separable *if it has the following general form:*

$$f(\boldsymbol{x}) = \sum_{i=1}^{k} f_i(\boldsymbol{x}_i), \; k = 2, \ldots, n, \tag{1}$$

*where $\boldsymbol{x} = (x_1, \ldots, x_n)$ is the decision vector, $\boldsymbol{x}_1, \ldots, \boldsymbol{x}_k$ are disjoint subcomponents of $\boldsymbol{x}$, $f_i(\cdot)$ denotes the subfunction of $\boldsymbol{x}_i$ ($i = 1, \ldots, k$), and $k$ is the number of subfunctions; otherwise, we say this function is* additively nonseparable.

This definition extends the Definition 2 in [13] which defines the function described by (1) to be *"partially additively separable"*, and makes the concept stricter by imposing a constraint on $k$. As this study focuses on additive separability, we omit the terms *"additive"* and *"additively"* hereinafter for the convenience of description.

**Definition 2.** *For a separable function $f(\boldsymbol{x})$ described by (1), let $X$ and $X_i$ ($i = 1, \ldots, k$) denote the variable sets of itself and its $i$th subfunction, respectively. We say that $\{X_1, \ldots, X_k\}$ is a* partition *of $X$ in the sense of function separability and any nonempty subset $X_i' \subseteq X_i$ is* separable *from any nonempty subset $X_{-i}' \subseteq (X \setminus X_i)$.*

If a function $f(\boldsymbol{x})$ is separable and can be formulized by (1), then the theoretical optimums of any decision variables in $X_i$ ($i = 1, \ldots, k$) is definitely not affected by any variables in $X_{-i}$ ($= X \setminus X_i$). From this perspective, there is no interdependency between any nonempty subset of $X_i$ and the one of $X_{-i}$, and accordingly we say they can be separated from each other. Especially, if a single variable $x_i$ ($i = 1, \ldots, n$) can be separated from $X \setminus \{x_i\}$, we simply say it is separable.

**Definition 3.** *For an n-dimensional separable function $f(\boldsymbol{x})$ with $X$ being its variable set, there necessarily exists a* maximum partition *of $X$ in the sense of function separability. Denote this partition as $\{X_1^*, \ldots X_k^*\}$, then this function is said to be* fully separable *if $k = n$; otherwise, it is said to be* partially separable.

As the maximum partition of $X$ in terms of function separability, $\{X_1^*, \ldots X_k^*\}$ implies that any variable subset $X_i^*$ ($i = 1, \ldots, k$) cannot be further partitioned. Besides, if each $X_i^*$ only contains a single element, then $k = n$ and the original function $f(\boldsymbol{x})$ is fully separable.

**Theorem 1.** *For an n-dimensional separable function $f(\boldsymbol{x})$ with $X$ being its variable set, if its two disjoint variable subsets $X_i$ and $X_j$ are separable from each other, then for $\forall \mathbf{cv} \in [\mathbf{lb}, \mathbf{ub}]^n$, $\forall \mathbf{x}_i^0, \mathbf{x}_i' \in [\mathbf{lb}_i, \mathbf{ub}_i]^{|X_i|}$, and $\forall \mathbf{x}_j^0, \mathbf{x}_j' \in [\mathbf{lb}_j, \mathbf{ub}_j]^{|X_j|}$, the following equation holds:*

$$\Delta(\mathbf{x}_i^0, \mathbf{x}_i' \mid \mathbf{x}_j^0) = \Delta(\mathbf{x}_i^0, \mathbf{x}_i' \mid \mathbf{x}_j'), \tag{2}$$

*where $\mathbf{lb}$ and $\mathbf{ub}$ denote the lower and the upper bounds of $\boldsymbol{x}$, respectively, $\mathbf{lb}_i = \mathbf{lb}(X_i)$ denotes the lower bound of the subcomponent $\boldsymbol{x}_i = \boldsymbol{x}(X_i)$, and $\Delta(\mathbf{x}_i^0, \mathbf{x}_i' \mid \mathbf{x}_j^0)$ is defined as*

$$\Delta(\mathbf{x}_i^0, \mathbf{x}_i' \mid \mathbf{x}_j^0) \triangleq f(\mathbf{cv} \leftarrow \mid \mathbf{x}_i', \mathbf{x}_j^0) - f(\mathbf{cv} \leftarrow \mid \mathbf{x}_i', \mathbf{x}_j^0) \tag{3}$$

*with $\mathbf{cv} \leftarrow \mid \mathbf{x}_i^0, \mathbf{x}_j^0$ denoting the complete solution obtained by inserting $\mathbf{x}_i^0$ and $\mathbf{x}_j^0$ into the corresponding positions in $\mathbf{cv}$.*

*Proof.* Please see the proof given in Section S-I of the supple-

mentary document. □

Theorem 1 indicates that (2) is just a necessary but not a sufficient condition for the separability between two variable subsets. This can be illustrated by the following example:

$$f(\boldsymbol{x}) = (x_1 - x_2)^2 + (x_2 - x_3)^2 + (x_3 - x_4)^2 + (x_4 - 1)^2.$$

For this function, it is obvious that the variable subsets $\{x_1, x_2\}$ and $\{x_4\}$ satisfy (2). However, they cannot be separated from each other since been linked together by the variable $x_3$. This kind of nonseparability can be referred to as indirect interdependency or indirect interaction [21], [23].

**Corollary 1.** *Let $X$ be the variable set of a function $f(\boldsymbol{x})$ and $X_i$ be a proper subset of $X$, then $X_i$ and $X_{-i}$ ($= X \setminus X_i$) are separable from each other if and only if they satisfy (2).*

*Proof.* Please see the proof given in Section S-I of the supplementary document. □

The definitions and theorems given above have the following merits compared with existing ones [13], [17], [22], [23]. First, they directly investigate separability from the perspective of variable subsets instead of variable individuals, which generalizes the concept of separability and facilitates reducing the separability detection times in the decomposition process. Second, the concept of separability is defined and proved based on the basic property of additive operation, and does not ask a function to be differentiable. Consequently, its application scope can be broadened.

### B. Interdependency Indicator and Its Distribution Analysis

Before defining the interdependency indicator, let us first study the influence of a variable subset $X_j$ on another subset $X_i$ by investigating

$$I(X_i, X_j) = \left| \Delta(\mathbf{x}_i^0, \mathbf{x}_i' \mid \mathbf{x}_j^0) - \Delta(\mathbf{x}_i^0, \mathbf{x}_i' \mid \mathbf{x}_j') \right|. \tag{4}$$

According to Theorem 1, $I(X_i, X_j)$ equals zero if $X_i$ and $X_j$ are separable. However, this does not strictly hold due to the computational roundoff error incurred by floating point operations [26]. This error cannot be accurately calculated, but may be roughly estimated. To weaken its negative influence, it is reasonable to subtract its estimate $\bar{e}$ from $I(X_i, X_j)$. The resulting difference $I(X_i, X_j) - \bar{e}$, however, may still differs much for different pairs of variable subsets since their contributions to the objective function may be rather imbalanced. Taking this phenomenon into account, we define the interdependency indicator as follows:

$$\varphi(X_i, X_j) = \max\left( \frac{I(X_i, X_j) - \bar{e}}{2 \max\left( \left| \Delta(\mathbf{x}_i^0, \mathbf{x}_i' \mid \mathbf{x}_j^0) \right|, \left| \Delta(\mathbf{x}_i^0, \mathbf{x}_i' \mid \mathbf{x}_j') \right| \right)}, \varepsilon_{\mathrm{M}} \right). \tag{5}$$

The denominator of the fractional expression in (5) reflects the maximum term for calculating $I(X_i, X_j)$ and generally normalizes $I(X_i, X_j) - \bar{e}$ into the range (0, 1). $\varepsilon_{\mathrm{M}}$ denotes the machine precision, which is $2^{-52}$ for the double precision floating point numbers of 64 bits. The introduction of $\varepsilon_{\mathrm{M}}$ ensures $\varphi(X_i, X_j) > 0$ even if $I(X_i, X_j) < \bar{e}$, which may be the case for a pair of separable $X_i$ and $X_j$. As for $\bar{e}$, this study employs the *"greatest lower bound"* $e_{\inf}$ developed in [26]; therefore, it can be formulized as

$$\overline{e} = \gamma_2 \left( \left| f(\mathbf{cv} \leftarrow |\mathbf{x}_i^0, \mathbf{x}_j^0) \right| + \left| f(\mathbf{cv} \leftarrow |\mathbf{x}_i^{'}, \mathbf{x}_j^0) \right| + \left| f(\mathbf{cv} \leftarrow |\mathbf{x}_i^0, \mathbf{x}_j^{'}) \right| + \left| f(\mathbf{cv} \leftarrow |\mathbf{x}_i^{'}, \mathbf{x}_j^{'}) \right| \right), \quad (6)$$

where $\gamma_2$ is an estimation coefficient and can be calculated by $\gamma_2 = \varepsilon_M / (1 - \varepsilon_M)$ [1].

From (3)-(6), it can be known that an indicator $\varphi(X_i, X_j)$ requires to sample and evaluate four solutions. As indicated in Theorem 1, any feasible $\mathbf{cv}$, $\mathbf{x}_i^0$, $\mathbf{x}_i^{'}$, $\mathbf{x}_j^0$, and $\mathbf{x}_j^{'}$ could be employed to generate these solutions. Nevertheless, $\mathbf{cv}$, which can be called context vector as in [9], is generally fixed to the lower bound (**lb**) of $\mathbf{x}$ for different indicators for the sake of reutilization, the subcomponent $\mathbf{x}_i$ is generally perturbed from its lower bound $\mathbf{x}_i^0 = \mathbf{lb}_i$ to the upper bound $\mathbf{x}_i^{'} = \mathbf{ub}_i$, and so is $\mathbf{x}_j$. For the convenience of description, we denote the solution $\mathbf{cv} \leftarrow |\mathbf{x}_i^0 = \mathbf{lb}_i, \mathbf{x}_j^0 = \mathbf{lb}_j$ as $\mathbf{x}^{l,l}$. Similar notations, including $\mathbf{x}^{u,l}$, $\mathbf{x}^{l,u}$, and $\mathbf{x}^{u,u}$, are given to the other three solutions. To illustrate the calculation of $\varphi(X_i, X_j)$ more clearly, Fig. S-1 in the supplementary document presents the relations among the four solutions and the intermediate items concerned therein.

This new interdependency indicator has three characteristics. First, it is oriented towards subsets rather than individuals of variables, although it is also applicable to individual variables. Second, it highlights the real interdependency among variable subsets by restraining the negative influence of the roundoff error. Finally, it fits well to different pairs of variable subsets of a LSBO problem due to its normalization operation. Based on these characteristics, we develop an interdependency indicator distribution analysis procedure (IDAP), which enables us to perform decomposition without specifying a threshold. For a set of $m$ indicator values contained in $\Phi$, IDAP analyzes their distribution according to the following three steps:

1) Sort the $m$ indicator values in an ascending order, and thus get a sequence $\varphi_{(1)}, \dots, \varphi_{(m)}$, where $\varphi_{(i)}$ $(i = 1, \dots, m)$ denotes the $i$th (smallest) element in $\Phi$. Here, the symbols of the two variable subsets concerned in each indicator are omitted for the convenience of description.

2) Compute the quotient between two adjacent indicator values in the sequence as follows:

$$\lambda_i = \frac{\varphi_{(i+1)}}{\varphi_{(i)}}, i = 1, \dots, m-1. \quad (7)$$

It is obvious that $\lambda_i \geq 1$ for $\forall \varphi_{(i)} \in \Phi$.

3) Find out the largest $\lambda_i$ $(i = 1, \dots, m-1)$ and denote it and the corresponding indicator in the denominator of (7) as $\lambda^*$ and $\varphi_{(i^*)}$, respectively, which means $\lambda^* = \varphi_{(i^*+1)} / \varphi_{(i^*)}$.

If $\lambda^*$ is much larger than all the other $\lambda^i$, it is reasonalbe to classify the indicator values in $\Phi$ into two classes, and the ones not greater than $\varphi_{(i^*)}$ (or not less than $\varphi_{(i^*+1)}$) can be considered to originate from separable variable subsets (or nonseparable variable subsets). On the contrary, we can deem that the indicator values in $\Phi$ are distributed rather uniformly and are only derived from separable or nonseparable variable subsets. For this case, if there exists an indicator value equal to $\varepsilon_M$, which

---

<sup>1</sup> $e_{inf}$ in [26] can be considered as a good estimate of the computational roundoff error since it is deduced according to the IEEE 754 standard [45]. However, strictly speaking, it is not a lower bound on the roundoff error because some real errors are less than their $e_{inf}$.

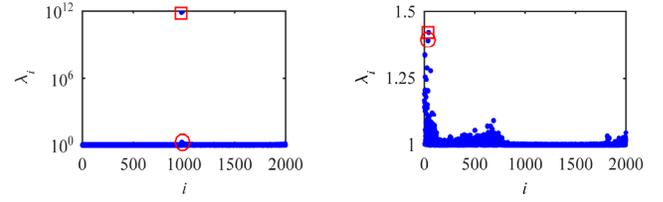

(a) The decision graph with an outlier (b) The decision graph without outliers
Fig. 1. Two examples of the decision graph.

---

**Algorithm 2:** $(type, \varphi_s, \varphi_n) \leftarrow \text{IDAP}(\Phi)$

**Input:** $\Phi$: a set of indicator values to be analyzed.
**Output:** $type$: the separability type of the indicator values in $\Phi$; $\varphi_s, \varphi_n$: the separability and the nonseparability thresholds.

1. Perform initialization: $\varphi_s \leftarrow null$, $\varphi_n \leftarrow null$, and $\varepsilon_M \leftarrow 1.0 \times 2^{-52}$;
2. Sort the elements in $\Phi$ in an ascending order and thus get a sequence: $\varphi_{(1)} \leq \varphi_{(2)} \leq \cdots \leq \varphi_{(m)}$, $m = |\Phi|$;
3. **for** $i = 1$ to $(m-1)$ **do**
4. | Calculate $\lambda_i$: $\lambda_i \leftarrow \varphi_{(i+1)} / \varphi_{(i)}$;
5. Find out the largest and the 2nd largest $\lambda$: $\lambda^* \leftarrow \varphi_{(i^*+1)} / \varphi_{(i^*)}$, and $\lambda^{**}$;
6. **if** $\lambda^* > 1000 \lambda^{**}$ **then** // the factor is set to 1000
7. | $type \leftarrow 'ps'$, $\varphi_s \leftarrow \varphi_{(i^*)}$, and $\varphi_n \leftarrow \varphi_{(i^*+1)}$; // partially separable
8. **elseif** for $\forall \varphi \in \Phi$, $\varphi > \varepsilon_M$ **then**
9. | $type \leftarrow 'ns'$; // nonseparable
10. **else**
11. | $type \leftarrow 'fs'$; // fully separable
12. **return** $type, \varphi_s, \varphi_n$;

---

implies that the corresponding influence quantity $I$ is less than or very close to its error estimate $\overline{e}$, then all pairs of variable subsets concerned in $\Phi$ can be regarded as separable. The reason behind this consists in that the influence quantity $I$ between two interdependent variable subsets is positively much larger than the corresponding $\overline{e}$.

At this moment, the key issue is how to measure the relative magnitude of $\lambda^*$ over the ones of all the other $\lambda_i$. Intuitively, we can construct a decision graph by taking $i = 1, \dots, m-1$ as abscissa and $\lambda_i$ as ordinate. If $\lambda^*$ is an outlier in the decision graph, it can be considered large enough. Fig. 1 shows two examples of this kind of decision graph. It can be obviously seen that the point corresponding to $\lambda^*$ in Fig. 1(a) is an outlier and the one in Fig. 1(b) is not. To automatically complete the decision process, we can directly judge whether $\lambda^*$ is larger than the second largest indicator value $\lambda^{**}$ by a large enough factor. As indicated by Fig. 1, this factor can be taken in a very wide range of values. Our preliminary experimental results show that it can be set to an arbitrary value within $[1.0 \times 10^3, 1.0 \times 10^8]$ for different LSBO problems without affecting their decomposition results.

To summarize, Algorithm 2 presents the pseudocode of IDAP. It is notable that if a significantly large $\lambda^*$ is found (step 6), the corresponding $\varphi_{(i^*)}$ and $\varphi_{(i^*+1)}$ are taken as the separability threshold $\varphi_s$ and the nonseparability threshold $\varphi_n$, respectively. They both are finally output together with the type of indicator values in $\Phi$, and can be employed to judge the separability of a pair of variable subsets.

IDAP provides an efficient way to detect the separability among variable subsets by analyzing the distribution of their indicator values. When applying this procedure, at least one of the following two conditions must be satisfied: 1) There are enough indicator values in $\Phi$ such that the distribution characteristics of all possible indicator values of the current LSBO instance can be reflected. 2) The set $\Phi$ contains a small number of indicator values, but it definitely includes at least an indicator value of a pair of separable variable subsets and the one of a pair of interdependent subsets if the current instance is partially separable. As a special case, if the indicator values of all pairs of variables are obtained, we can complete decomposition by applying Algorithm 2. However, this kind of direct application cannot reduce the FE consumption.

### C. Instance Type Identification Procedure

As stated in Section III-A, LSBO instances can be classified into three types: fully separable instances, nonseparable instances, and partially separable instances. If we know an instance is fully separable or nonseparable, there is no need to decompose it step by step, and thus many FEs can be saved. Different from existing decomposition algorithms which cope with these three types of instances in the same manner without considering the particularities of the first two types, FDG, the decomposition algorithm proposed in this study, explicitly identifies the type of an instance with a specially designed instance type identification procedure (ITIP) before decomposing it. The basic steps of ITIP can be described as follows:

    1) Find out at least a pair of separable variable subsets and a pair of nonseparable subsets under the assumption that the current instance is partially separable.

    2) Detect the separability of all pairs of variable subsets selected in step 1) by investigating their interdependency indicator values.

    3) Determine the type of the instance. If each pair of variable subsets are separable (or nonseparable), the instance is fully separable (or nonseparable); otherwise, it is really partially separable.

It is obvious that the key of ITIP lies in step 1. To achieve that goal, a sophisticated variable subset selection strategy, which involves two rules, is developed. The first rule aims to find out at least a pair of nonseparable variable subsets even for an almost fully separable instance. At each time of selection, it takes the two variable subsets generated by randomly halving the whole variable set. In contrast, the second rule randomly selects two variable individuals for the purpose of getting a pair of separable variables. To measure the success probabilities of these two selection rules, we give the following two theorems.

**Theorem 2.** *Assume that $f(\boldsymbol{x})$ is a partially separable function of $n$ dimensions, it involves $k$ disjoint sets of nonseparable variables with each set containing $s_i$ ($s_i \geq 2$, $i = 1, \ldots, k$) variables, and $\sum_{i=1}^{k} s_i \leq n$ holds. If we randomly partition the whole variable set into two (nearly) equal-sized subsets, then these two subsets are separable from each other with a probability of*

$$P_{s1} = \frac{C_{n_s}^{\lfloor n/2 \rfloor} + \sum_{i=1}^{k} C_{n_s}^{\lfloor n/2 \rfloor - s_i} + \sum_{i=1}^{k} \sum_{j=1, j \neq i}^{k} C_{n_s}^{\lfloor n/2 \rfloor - s_i - s_j} + \ldots + C_{n_s}^{\lfloor n/2 \rfloor - \sum_{i=1}^{k} s_i}}{C_n^{\lfloor n/2 \rfloor}}, \quad (8)$$

*where $\lfloor \cdot \rfloor$ denotes the round down operator, $n_s = n - \sum_{i=1}^{k} s_i$*

*denotes the number of separable variables, and the combination number $C_p^q$ is defined to be zero if $p < 0$, $q < 0$, or $p < q$. The probability of getting at least such a pair of nonseparable variable subsets within $l$ ($l \geq 1$) trials is*

$$P_{n1} = 1 - P_{s1}^{\ l}. \quad (9)$$

*Proof.* Please see the proof given in Section S-I of the supplementary document. □

**Corollary 2.** *If we further assume the $k$ nonseparable variable subsets of the function $f(\boldsymbol{x})$ considered in Theorem 2 have the same size of $s$, then any two variable subsets generated by randomly halving the original variable set are separable from each other with a probability of*

$$P_{s1}^{'} = \frac{\sum_{i=0}^{k} C_k^i C_{n_s}^{\lfloor n/2 \rfloor - is}}{C_n^{\lfloor n/2 \rfloor}}, \quad (10)$$

*Proof.* Please see the proof given in Section S-I of the supplementary document. □

**Theorem 3.** *For the function described in Theorem 2, its two randomly-selected variables are nonseparable from each other with a probability of*

$$P_{n2} = \frac{\sum_{i=1}^{k} C_{s_i}^{2}}{C_n^{2}}, \quad (11)$$

*and the probability of getting at least such a pair of separable variables within $l$ ($l \geq 1$) trials is*

$$P_{s2} = 1 - P_{n2}^{\ l}. \quad (12)$$

*Proof.* Please see the proof given in Section S-I of the supplementary document. □

Two examples are provided at the end of Section S-I of the supplementary document. They verify that the two selection rules can achieve their respective goals with a probability very close to 1 even when a few trials are allowed.

For each pair of selected variable subsets, ITIP calculates their interdependency indicator. After that, it detects the separability of all pairs of selected variable subsets by applying IDAP (Algorithm 2), and finally identifies the type of the current instance. Algorithm 3 presents the pseudocode of ITIP, where steps 2-6 and steps 7-11 fulfil the first and the second selection rules, respectively, and compute interdependency indicators by sampling and evaluating corresponding solutions. To be conservative, the trial number is set to 10 in steps 2 and 7, although a smaller number also works. It is notable that, besides $\mathbf{cv} = \mathbf{lb}$, the complete solution $\mathbf{ub}$ can also be reutilized in step 4 since it can be used as $\mathbf{x}^{u,u}$ required by different pairs of variable subsets.

The performance of ITIP is pretty appealing since it can identify the type of an instance with very a few FEs. Despite that, it may fail on some nonseparable instances involving indirect interdependencies because it regards two indirect interdependent variables as separable. Fortunately, this defect can be remedied by the following decomposition procedure developed for partially separable instances.

---

**Algorithm 3**: $(type, \varphi_s, \varphi_n, feNum) \leftarrow \text{ITIP}(\mathbf{lb}, \mathbf{ub}, y^{\mathbf{lb}}, y^{\mathbf{ub}})$

---

**Input**: $\mathbf{lb}$, $\mathbf{ub}$: the lower and the upper bounds of the decision vector; $y^{\mathbf{lb}}, y^{\mathbf{ub}}$: the fitness values of $\mathbf{lb}$ and $\mathbf{ub}$.

**Output**: $type$: the type of the current instance; $\varphi_s, \varphi_n$: the separability and the nonseparability thresholds; $feNum$: the number of FEs.

1. Perform initialization: $feNum \leftarrow 0$ and $\Phi \leftarrow \varnothing$ ;
2. **for** $i = 1$ to 10 **do** // $l = 10$ and $\Phi$ is employed to save indicator values
3.      Randomly partition $X$ into two equal-sized subsets: $X_1$ and $X_2$;
4.      Generate and evaluate four solutions for calculating $\varphi(X_1, X_2)$ :
$$\begin{cases} \mathbf{x}^{1,1} \leftarrow \mathbf{lb}, y^{1,1} \leftarrow y^{\mathbf{lb}}; \mathbf{x}^{u,1} \leftarrow \mathbf{x}^{1,1}, \mathbf{x}^{u,1}(X_1) \leftarrow \mathbf{ub}(X_1), y^{u,1} \leftarrow f(\mathbf{x}^{u,1}) \\ \mathbf{x}^{1,u} \leftarrow \mathbf{x}^{1,1}, \mathbf{x}^{1,u}(X_2) \leftarrow \mathbf{ub}(X_2), y^{1,u} \leftarrow f(\mathbf{x}^{1,u}); \mathbf{x}^{u,u} \leftarrow \mathbf{ub}, y^{u,u} \leftarrow y^{\mathbf{ub}} \end{cases};$$
5.      Calculate $\varphi(X_1, X_2)$ according to (3)-(6);
6.      Perform updating: $feNum \leftarrow feNum + 2$ and $\Phi \leftarrow \Phi \cup \varphi(X_1, X_2)$ ;
7. **for** $i = 1$ to 10 **do** // $l = 10$
8.      Randomly select two different variables, $x_1$ and $x_2$, from $X$;
9.      Generate and evaluate four solutions for calculating $\varphi(x_1, x_2)$ :
$$\begin{cases} \mathbf{x}^{1,1} \leftarrow \mathbf{lb}, y^{1,1} \leftarrow y^{\mathbf{lb}}; \mathbf{x}^{u,1} \leftarrow \mathbf{x}^{1,1}, \mathbf{x}^{u,1}(x_1) \leftarrow \mathbf{ub}(x_1), y^{u,1} \leftarrow f(\mathbf{x}^{u,1}) \\ \mathbf{x}^{1,u} \leftarrow \mathbf{x}^{1,1}, \mathbf{x}^{1,u}(x_2) \leftarrow \mathbf{ub}(x_2), y^{1,u} \leftarrow f(\mathbf{x}^{1,u}) \\ \mathbf{x}^{u,u} \leftarrow \mathbf{x}^{1,u}, \mathbf{x}^{u,u}(x_1) \leftarrow \mathbf{ub}(x_1), y^{u,u} \leftarrow f(\mathbf{x}^{u,u}) \end{cases};$$
10.    Calculate $\varphi(x_1, x_2)$ according to (3)-(6);
11.    Perform updating: $feNum \leftarrow feNum + 3$ and $\Phi \leftarrow \Phi \cup \varphi(x_1, x_2)$ ;
12. Analyze the distribution of indicator values in $\Phi$ by executing IDAP:
$$(type, \varphi_s, \varphi_n) \leftarrow \text{IDAP}(\Phi);$$
13. **return** $type, \varphi_s, \varphi_n, feNum$ ;

---

### D. Binary-Tree-Based Decomposition Procedure

If a LSBO instance is judged to be partially separable, it is necessary to further decompose it, the goal of which is to keep each group of nonseparable variables together and put each pair of separable variables into different groups. Here the key operation consists in that, for a variable subset $X_i$ nonseparable from another one $X_j$, how to efficiently capture the variables in $X_j$ that are really interdependent with it. A direct way is to check the interdependency between $X_i$ and each variable in $X_j$. However, this requires too many FEs. To alleviate this issue, this study proposes a binary-tree-based decomposition procedure (BTDP) which can significantly reduce separability detection times and effectively reutilize evaluated solutions by converting the key decomposition process into a search process in a binary tree.

Concretely, BTDP first initializes the variable subset $X_j$ as the root node of a binary tree, and then examines its separability from $X_i$. If they are interdependent with each other, BTDP divides $X_j$ into two (nearly) equal-sized subsets and takes them as child nodes. Next, it further detects the separability between $X_i$ and each child node. BTDP stops dividing a node if it is separable from $X_i$ or just involves a single variable. In this way, a binary tree is dynamically generated along with the decomposition process.

As mentioned in Section III-B, four evaluated solutions are required to calculate an interdependency indicator. Interestingly, if a certain sampling rule is obeyed, all the evaluated solutions can be reutilized during the search process in a binary tree. To illustrate this rule, let us take the schematic binary tree shown in Fig. 2 as an example. For the calculation of $\varphi(X_i, X_j)$, it is feasible to directly sample and evaluate four solutions according to the guideline given in Fig. S-1 in the supplementary document. However, when calculating $\varphi(X_i, X_{j1})$, its two solutions $\mathbf{x}^{1,1}$ and $\mathbf{x}^{u,1}$ can be inherited from the corresponding ones sampled for $\varphi(X_i, X_j)$, and its other two solutions $\mathbf{x}^{1,u}$ and $\mathbf{x}^{u,u}$ should be generated by perturbing $\mathbf{x}_i = \mathbf{x}(X_i)$ from $\mathbf{lb}_i$ to $\mathbf{ub}_i$ and keeping $\mathbf{x}_{j1}$ at $\mathbf{ub}_{j1}$. As for $\varphi(X_i, X_{j2})$, it is exciting to find that if we employ the solution $\mathbf{x}^{1,u}$ of $\varphi(X_i, X_j)$ as its context vector, its two solutions $\mathbf{x}^{1,1}$ and $\mathbf{x}^{u,1}$ can take the two ones $\mathbf{x}^{1,u}$ and $\mathbf{x}^{u,u}$ of $\varphi(X_i, X_{j1})$, respectively, and its other two solutions $\mathbf{x}^{1,u}$ and $\mathbf{x}^{u,u}$ can be inherited from the corresponding ones of $\varphi(X_i, X_j)$. This means that no new solution needs to be generated for the node $X_{j2}$ at all. Fig. 3 graphically demonstrates the relations among these solutions. It is notable that although $X_j$ corresponds to the root node of the binary tree shown in Fig. 2, no particularity of the root node is used in the solution sampling rule described above. Therefore, this rule applies to a general parent node and its two child nodes, such as $X_{j2}$, $X_{j21}$, and $X_{j22}$ in Fig. 2. On the other side, although this rule samples two new solutions for the left child node and enables the right one to inherit solutions from its brother and parent nodes, the reverse is also feasible.

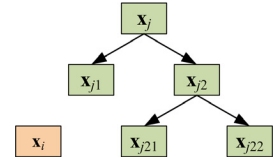

Fig. 2. A schematic example of the binary tree used for decomposition.

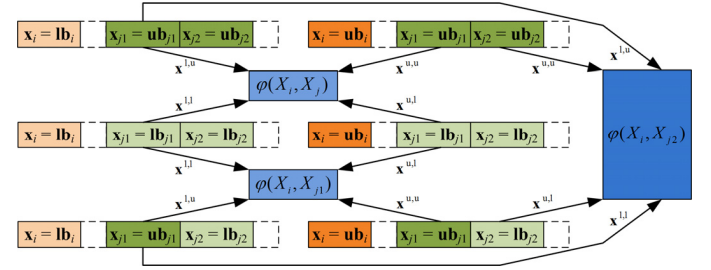

Fig. 3. The relations among the solutions sampled for a parent node and its two child nodes in a binary tree.

To traverse a binary tree, two strategies, i.e., the breadth-first strategy and the depth-first one, are available. For the application of the binary tree in this study, there is no essential difference between these two strategies. Algorithm 4 presents the pseudocode of BTDP, where the breadth-first traverse strategy is adopted. The *"rootNode"* in step 2 is defined as a struct and employed to record the whole variable subset $X_2$ and the four solutions used for calculating its interdependency indicator. Then in step 4, this node is pushed into *"nodeQueue"* which is a queue defined to store all possible nodes. Next, BTDP carries on a loop process in steps 5-16 until *"nodeQueue"* is empty. During the loop process, the first node in the front of *"nodeQueue"* is tackled in priority (step 6); once its two child nodes, i.e., *"leftNode"* and *"righNode"*, are generated (step 14), they are sequentially pushed into the rear of *"nodeQueue"* (step 16). By this means, the *"first-in, first-out"* mechanism required by the breadth-first traverse strategy can be realized.

Step 8 judges the separability between $X_1$ and the variable subset involved in the current node with the procedure given in Algorithm 5. It compares the corresponding interdependency

---

**Algorithm 4**: $(X_1, feNum, \varphi_s, \varphi_n) \leftarrow \text{BTDP}(X_1, X_2, \mathbf{lb}, \mathbf{ub}, y^{\mathbf{lb}}, \varphi_s, \varphi_n)$

---

**Input**: $X_1, X_2$: a pair of variable subsets to be detected; $\mathbf{lb}$, $\mathbf{ub}$: the lower and the upper bounds of the decision vector; $y^{\mathbf{lb}}$: the fitness value of $\mathbf{lb}$; $\varphi_s, \varphi_n$: the separability and the nonseparability thresholds.

**Output**: $X_1$: the variables in $X_2$ interdependent with $X_1$; $feNum$: the number of FEs; $\varphi_s, \varphi_n$: the updated thresholds.

1. Perform initialization: $X_1 \leftarrow \varnothing$, $feNum \leftarrow 0$, $\varphi_s \leftarrow \varphi_s$, and $\varphi_n \leftarrow \varphi_n$;
2. Set root node: $rootNode \begin{cases} .subset \leftarrow X_2 \\ \mathbf{x}^{\mathbf{l},1} \leftarrow \mathbf{lb}, \ .y^{\mathbf{l},1} \leftarrow y^{\mathbf{lb}} \\ \mathbf{x}^{\mathbf{u},1} \leftarrow \mathbf{x}^{\mathbf{l},1}, \ \mathbf{x}^{\mathbf{u},1}(X_1) \leftarrow \mathbf{ub}(X_1), \ .y^{\mathbf{u},1} \leftarrow f(\mathbf{x}^{\mathbf{u},1}) \\ \mathbf{x}^{\mathbf{l},u} \leftarrow \mathbf{x}^{\mathbf{l},1}, \ \mathbf{x}^{\mathbf{l},u}(X_2) \leftarrow \mathbf{ub}(X_2), \ .y^{\mathbf{l},u} \leftarrow f(\mathbf{x}^{\mathbf{l},u}) \\ \mathbf{x}^{\mathbf{u},u} \leftarrow \mathbf{x}^{\mathbf{l},u}, \ \mathbf{x}^{\mathbf{u},u}(X_1) \leftarrow \mathbf{ub}(X_1), \ .y^{\mathbf{u},u} \leftarrow f(\mathbf{x}^{\mathbf{u},u}) \end{cases}$;
3. Update $feNum$: $feNum \leftarrow feNum + 3$;
4. Push $rootNode$ into the rear of a queue $nodeQueue$;
5. **while** $nodeQueue \ != \varnothing$ **do**
6. | Pop up the 1st node in the front of $nodeQueue$ and save it into $curNode$;
7. | Calculate $\varphi(X_1, curNode.subset)$ according to (3)-(6);
8. | Judge the separability between $X_1$ and $curNode.subset$: $(isSep, \varphi_s, \varphi_n) \leftarrow \text{JudgeSep}(\varphi(X_1, curNode.subset), \varphi_s, \varphi_n)$;
9. | **if** $isSep = false$ **then**
10. | | **if** $|curNode.subset| == 1$ **then** //find a variable interdependent with $X_1$
11. | | | Update $X_1$: $X_1 \leftarrow X_1 \cup curNode.subset$;
12. | | **else**
13. | | | Randomly halve $curNode.subset$ into two subsets: $X_l$ and $X_r$;
14. | | | Generate two child nodes related to $X_l$ and $X_r$:
$leftNode \begin{cases} leftNode \leftarrow curNode \\ leftNode \begin{cases} .subset \leftarrow X_l \\ \mathbf{x}^{\mathbf{l},u}(X_r) \leftarrow \mathbf{lb}(X_r), \ .y^{\mathbf{l},u} \leftarrow f(\mathbf{x}^{\mathbf{l},u}) \\ \mathbf{x}^{\mathbf{u},u}(X_r) \leftarrow \mathbf{lb}(X_r), \ .y^{\mathbf{u},u} \leftarrow f(\mathbf{x}^{\mathbf{u},u}) \end{cases} \end{cases}$,
$rightNode \begin{cases} rightNode \leftarrow curNode \\ rightNode \begin{cases} .subset \leftarrow X_r \\ \mathbf{x}^{\mathbf{l},1} \leftarrow leftNode.\mathbf{x}^{\mathbf{l},u}, \ .y^{\mathbf{l},1} \leftarrow leftNode.y^{\mathbf{l},u} \\ \mathbf{x}^{\mathbf{u},1} \leftarrow leftNode.\mathbf{x}^{\mathbf{u},u}, \ .y^{\mathbf{u},1} \leftarrow leftNode.y^{\mathbf{u},u} \end{cases} \end{cases}$;
15. | | | Update $feNum$: $feNum \leftarrow feNum + 2$;
16. | | | Sequentially push $leftNode$, $rightNode$ into the rear of $nodeQueue$;
17. **return** $X_1, feNum, \varphi_s, \varphi_n$;

---

**Algorithm 5**: $(isSep, \varphi_s, \varphi_n) \leftarrow \text{JudgeSep}(\varphi, \varphi_s, \varphi_n)$

---

**Input**: $\varphi$: the indicator value to be judged; $\varphi_s, \varphi_n$: the separability and the nonseparability thresholds.

**Output**: $isSep$: the separability flag; $\varphi_s, \varphi_n$: the updated thresholds.

1. Perform initialization: $\varphi_s \leftarrow \varphi_s$ and $\varphi_n \leftarrow \varphi_n$;
2. **if** $\varphi / \varphi_s < \varphi_n / \varphi$ **then** // involve the case of $\varphi < \varphi_s$
3. | $isSep \leftarrow true$ and $\varphi_s \leftarrow \max(\varphi, \varphi_s)$;
4. **else** // involve the case of $\varphi > \varphi_n$
5. | $isSep \leftarrow false$ and $\varphi_n \leftarrow \min(\varphi, \varphi_n)$;
6. **return** $isSep, \varphi_s, \varphi_n$;

---

indicator with two thresholds $\varphi_s'$ and $\varphi_n'$ which can be generally initialized by the ones returned by ITIP. Note that, to tackle the case of $\varphi \in (\varphi_s, \varphi_n)$, Algorithm 5 compares $\varphi / \varphi_s$ with $\varphi_n / \varphi$. If the former is less (or larger) than the latter, it judges the two variable subsets to be separable (or nonseparable) and updates the threshold $\varphi_s$ (or $\varphi_n$) with $\varphi$.

### E. Decomposition of Partially Separable Instances

Equipped with BTDP, FDG decomposes a partially separable instance as follows: It starts from an arbitrary variable $X_i = \{x_i \in X\}$ and tries to find variables interdependent with $X_i$ from the remaining ones $X_c = X \setminus X_i$ by applying BTDP.

If some variables are really found, FDG combines them into $X_i$ and further detects the interdependency between the updated $X_i$ and $X_c = X \setminus X_i$ for the purpose of capturing the variables indirectly interdependent with the original $X_i$. This iterative process continues until no new interdependent variable can be found. Then FDG takes $X_i$ as a nonseparable variable group or a single separable variable according to the number of elements it contains, and moves on to an untreated variable in $X_c$. The whole decomposition process terminates when all the variables are grouped.

From the above process, it can be known that if some separable variables are involved in a partially separable instance, they will reside in the candidate variable subset $X_c$ until being directly treated. This may increase many interdependency detection times. A way to avoid this issue is to exclude this kind of variables in advance by detecting the separability of each of them as done by FII [25]. However, this requires $n$ detections for an $n$-dimensional instance and takes no effect if the instance does not involve any separable variable at all. To balance this contradiction, we develop a new separable variable exclusion procedure (SVEP). It first randomly selects $l$ ($l \ll n$) variables, then detects their separability. If a separable variable is really found, it further investigates the remaining variables and finally excludes all the separable ones from the candidate variable subset $X_c$; otherwise, it directly quits this trial procedure. It can be easily proved that, for the general partially separable function described in Theorem 2, SVEP can get at least a separable variable with a probability of

$$P_{s3} = 1 - \frac{C_{n_n}^l}{C_n^l} > 1 - \left(\frac{n_n}{n}\right)^l, \tag{13}$$

where $n_n = \sum_{i=1}^k s_i$ denotes the number of nonseparable variables. Obviously, this probability exponentially increases with the increase of $l$. Algorithm 6 shows the pseudocode of SVEP, where steps 2-11 carry on a loop which may be discontinued if all of the first $l$ variables are found nonseparable (steps 3-4).

---

**Algorithm 6**: $(seps, feNum, \varphi_s, \varphi_n) \leftarrow \text{SVEP}(\mathbf{lb}, \mathbf{ub}, y^{\mathbf{lb}}, \varphi_s, \varphi_n)$

---

**Input**: $\mathbf{lb}$, $\mathbf{ub}$: the lower and the upper bounds of the decision vector; $y^{\mathbf{lb}}, y^{\mathbf{ub}}$: the fitness values of $\mathbf{lb}$ and $\mathbf{ub}$. $\varphi_s, \varphi_n$: the separability and the nonseparability thresholds.

**Output**: $seps$: a set of separable variables; $feNum$: the number of FEs; $\varphi_s, \varphi_n$: the updated thresholds.

1. Perform initialization: $seps \leftarrow \varnothing$, $feNum \leftarrow 0$, $\varphi_s \leftarrow \varphi_s$, $\varphi_n \leftarrow \varphi_n$, and $flag \leftarrow false$;
2. **for** $i = 1$ to $|X|$ **do**
3. | **if** $i > 10$ and $flag == false$ **then** //As in ITIP, $l$ is set to 10
4. | | **break**;
5. | Randomly select a variable $x$ from $X$ without repetition;
6. | Generate and evaluate four solutions for calculating $\varphi(\{x\}, X \setminus \{x\})$:
$\begin{cases} \mathbf{x}^{\mathbf{l},1} \leftarrow \mathbf{lb}, \ .y^{\mathbf{l},1} \leftarrow y^{\mathbf{lb}}; \mathbf{x}^{\mathbf{u},1} \leftarrow \mathbf{x}^{\mathbf{l},1}, \mathbf{x}^{\mathbf{u},1}(x) \leftarrow \mathbf{ub}(x), y^{\mathbf{u},1} \leftarrow f(\mathbf{x}^{\mathbf{u},1}) \\ \mathbf{x}^{\mathbf{l},u} \leftarrow \mathbf{x}^{\mathbf{u},u}, \mathbf{x}^{\mathbf{l},u}(x) \leftarrow \mathbf{lb}(x), y^{\mathbf{l},u} \leftarrow f(\mathbf{x}^{\mathbf{l},u}); \mathbf{x}^{\mathbf{u},u} \leftarrow \mathbf{ub}, y^{\mathbf{u},u} \leftarrow y^{\mathbf{ub}} \end{cases}$;
7. | Update $feNum$: $feNum \leftarrow feNum + 2$;
8. | Calculate $\varphi(\{x\}, X \setminus \{x\})$ according to (3)-(6);
9. | Judge the separability of $x$: $(isSep, \varphi_s, \varphi_n) \leftarrow \text{JudgeSep}(\varphi(\{x\}, X \setminus \{x\}), \varphi_s, \varphi_n)$;
10. | **if** $isSep$ **then**
11. | | Update updating: $seps \leftarrow seps \cup \{x\}$ and $flag \leftarrow true$;
12. **return** $seps$ and $feNum$;

---

---

**Algorithm 7:** $(nonseps, seps, feNum) \leftarrow \text{PSDP}(\mathbf{lb}, \mathbf{ub}, y^{\mathbf{lb}}, y^{\mathbf{ub}}, \varphi_s, \varphi_n)$

---

**Input:** $\mathbf{lb}$, $\mathbf{ub}$: the lower and the upper bounds of the decision vector;
$\quad y^{\mathbf{lb}}, y^{\mathbf{ub}}$: the fitness values of $\mathbf{lb}$ and $\mathbf{ub}$;
$\quad \varphi_s, \varphi_n$: the separability and the nonseparability thresholds.
**Output:** $nonseps$: a group of nonseparable variable subsets; $seps$: a set of separable variables; $feNum$: the number of FEs.

---

1. Perform initialization: $nonseps \leftarrow \varnothing$ ;
2. Exclude separable variables:
$\quad (seps, feNum, \varphi_s, \varphi_n) \leftarrow \text{SVEP}(\mathbf{lb}, \mathbf{ub}, y^{\mathbf{lb}}, y^{\mathbf{ub}}, \varphi_s, \varphi_n)$ ;
3. Get remaining variables: $X_c \leftarrow X \setminus seps$ ;
4. Initialize $X_1$ with a variable $x$ randomly selected from $X_c$ ;
5. **while** $|X_1| < |X_c|$ **do** // imply a nonempty $X_2$
6. | Set $X_2$ : $X_2 \leftarrow X_c \setminus X_1$ ;
7. | Find the variables interdependent with $X_1$ from $X_2$ :
$\quad (X_1', feNum', \varphi_s, \varphi_n) \leftarrow \text{BTDP}(X_1, X_2, \mathbf{lb}, \mathbf{ub}, \varphi_s, \varphi_n)$ ;
8. | Update $feNum$: $feNum \leftarrow feNum + feNum'$ ;
9. | **if** $X_1' == \varnothing$ **then**
10. | | **if** $|X_1| == 1$ **then**
11. | | | Update $seps$: $seps \leftarrow seps \cup X_1$ ;
12. | | **else**
13. | | | Update $nonseps$: $nonseps \leftarrow nonseps \cup \{X_1\}$ ;
14. | | Update $X_c$: $X_c \leftarrow X_2$ ;
15. | | Reinitialize $X_1$ with a variable $x$ randomly selected from $X_c$ ;
16. | **else**
17. | | Update $X_1$: $X_1 \leftarrow X_1 \cup X_1'$ ;
18. **if** $|X_1| == 1$ **then**
19. | Update $seps$: $seps \leftarrow seps \cup X_1$ ;
20. **else**
21. | Update $nonseps$: $nonseps \leftarrow nonseps \cup \{X_1\}$ ;
22. **return** $nonseps$, $seps$, and $feNum$;

---

**Algorithm 8:** $(nonseps, seps, feNum) \leftarrow \text{FDG}(\mathbf{lb}, \mathbf{ub})$

---

**Input:** $\mathbf{lb}$, $\mathbf{ub}$: the lower and the upper bounds of the decision vector.
**Output:** $nonseps$: a group of nonseparable variable subsets; $seps$: a set of separable variables; $feNum$: the number of FEs.

---

1. Perform initialization: $nonseps \leftarrow \varnothing$ , $seps \leftarrow \varnothing$ , and $feNum \leftarrow 0$ ;
2. Evaluate $\mathbf{lb}$ and $\mathbf{ub}$: $y^{\mathbf{lb}} \leftarrow f(\mathbf{lb})$ and $y^{\mathbf{ub}} \leftarrow f(\mathbf{ub})$ , and update $feNum$: $feNum \leftarrow feNum + 2$ ;
3. Identify instance type: $(type, \varphi_s, \varphi_n, feNum') \leftarrow \text{ITIP}(\mathbf{lb}, \mathbf{ub}, y^{\mathbf{lb}}, y^{\mathbf{ub}})$ ;
4. Update $feNum$: $feNum \leftarrow feNum + feNum'$ ;
5. **if** $type == 'ns'$ **then**
6. | Update $nonseps$: $nonseps \leftarrow \{X\}$ ;
7. **elseif** $type == 'fs'$ **then**
8. | Update $seps$: $seps \leftarrow X$ ;
9. **else**
10. | Decompose the identified partially separable instance:
$\quad (nonseps, seps, feNum') \leftarrow \text{PSDP}(\mathbf{lb}, \mathbf{ub}, y^{\mathbf{lb}}, y^{\mathbf{ub}}, \varphi_s, \varphi_n)$ ;
11. | Update $feNum$: $feNum \leftarrow feNum + feNum'$ ;
12. **return** $nonseps$, $seps$, and $feNum$;

---

By integrating SVEP and BTDP together, Algorithm 7 presents the pseudocode of the whole decomposition procedure for partially separable instances (PSDP). After attempting to exclude separable variables with SVEP in step 2, PSDP groups the remaining variables by iteratively applying BTDP in steps 5-17. Note that although PSDP is mainly designed for partially separable instances, it can also be applied to fully separable and nonseparable ones. For the former, PSDP achieves decomposition by just applying SVEP; while for the latter, it mainly relies on BTDP.

### F. Framework and Computational Complexity of FDG

The whole framework of FDG is presented in Algorithm 8, where step 3 identifies the type of the current instance with ITIP. If the instance is judged to be nonseparable or fully separable, its decomposition result is directly given in step 6 or 8, respectively; otherwise, it is achived by invoking PSDP in step 10.

As for the computational complexity of a decomposition algorithm, it is generally analyzed in terms of the number of FEs since FE is much more time-consuming than other algorithmic operations. To highlight main characteristics of FDG and facilitate the analysis, the case of indirect interdependency is not considered below.

The ITIP in FDG detects separability for $2l$ pairs of variable subsets. As indicated at the end of Section III-C, each pair of variable subsets sampled by the first (or second) rule needs 2 (or 3) new FEs and shares the other 2 (or 1) FEs with other pairs. Therefore, a total number of $5l + 2$ FEs are required by ITIP, and then the computational complexity of FDG on fully separable and nonseparable instances is $O(5l + 2) = O(l)$ , where $l$ is generally a constant much less than $n$.

For an $n$-dimensional partially separable instance which has no separable variable but involves $k$ groups of nonseparable variables with each group containing $s$ elements, the SVEP in FDG only checks $l$ variables with $2l$ FEs being consumed, and the BTDP in FDG constructs a binary tree for each group of nonseparable variables. In the worst case, BTDP uniformly assigns the $(s-1)$ nonseparable partners of a selected variable to all the nodes in each layer of the binary tree because this tends not to exclude any node and results in the most number of detections, and the total computational complexity of FDG in this case is $O(n \log_2 k)$ . In the best case, BTDP concentratively assigns the nonseparable partners of a selected variable to one of the nodes in each layer such that those variables separable from the selected one can be quickly excluded and the number of detections can be reduced to the greatest extent possible. The corresponding computational complexity becomes $O(\max(n, k \log_2 k))$ . For the detailed deduction, please refer to Section S-III in the supplementary document.

For an $n$-dimensional partially separable instance which involves $k$ nonseparable variable subsets of size $s$ and $(n - ks)$ separable variables, FDG will invoke SVEP to detect the separability of all variables and consumes $2n$ FEs. After this operation, the original instance can be considered as a $ks$-dimensional partially separable instance without any separable variable. Then the computational complexity of FDG on this type of instances is $O(\max(n, k \log_2 k))$ in the worst case and $O(\max(n, k \log_2 k))$ in the best case.

Table I summarizes the computational complexity of FDG and five existing decomposition algorithms, including DG [13], GDG [22], FII [25], DG2 [26], and RDG [17], on LSBO problems involving subcomponents of different sizes. It can be seen that FDG possesses the lowest computational complexity on all types of problems.

FDG shares some similarities with RDG since RDG performs its basic decomposition operation in a recursive manner which can also be described by a binary tree. Nevertheless, FDG still differs distinctly from RDG in the following four aspects:

1) RDG tackles different types of LSBO instances in the same way, while FDG performs decomposition by first identifying the type of an instance and then employing

TABLE I

Computational complexity of DG, GDG, FII, DG2, RDG, and FDG on different types of problems. $n : k \times s - n_s$ denotes an $n$-dimensional problem which involves $k$ nonseparable variable subsets of size $s$ and $n_s$ separable variables.

| Dimensions of subcomponents | DG | GDG | FII | DG2 | RDG | FDG |
|---|---|---|---|---|---|---|
| $n : 0 \times 0 - n$ | $O(n^2)$ | $O(n^2)$ | $O(n)$ | $O(n^2)$ | $O(n)$ | $O(l)$ |
| $n : k \times s - (n-ks)$ | $O(n^2)$ | $O(n^2)$ | $O(\max(n,k^2s))$ | $O(n^2)$ | $O(\max(n,ks\log_2 n))$ | $O(\max(n,ks\log_2 k))$ |
| $n : k \times s - 0$ | $O(kn)$ | $O(n^2)$ | $O(kn)$ | $O(n^2)$ | $O(n\log_2(n))$ | $O(n\log_2 k)$ |
| $n : 1 \times n - 0$ | $O(n)$ | $O(n^2)$ | $O(n)$ | $O(n^2)$ | $O(n)$ | $O(l)$ |

different strategies according to the identification result. As a consequence, FDG can decompose fully separable and nonseparable instances at negligible cost.

2) Due to the lack of instance-type information, RDG does not attempt to exclude separable variables involved in partially separable instances in advance, although this operation is likely to effectively reduce the number of variables to be grouped and thus interdependency detection times.

3) When detecting the interdependency between a variable subset $X_i$ and each of the two subsets $X_{j1}$ and $X_{j2}$ obtained by partitioning $X_j$, RDG needs to sample and evaluate three solutions [17]. By contrast, FDG just requires two FEs for $X_{j1}$ or $X_{j2}$, and does not require any FE for the other one, which essentially reduces its computational complexity on partially separable instances.

4) The two algorithms also adopt different interdependency indicators. RDG inherits the one developed by GDG [22] and requires users to specify a control coefficient which is somewhat sensitive, while FDG employs a normalized interdependency indicator whose threshold can be adaptively generated by analyzing its distribution.

## IV. Experimental Studies

To evaluate the performance of FDG, we first tested its decomposition efficiency and accuracy by comparing it with several popular learning-based decomposition algorithms, then further investigated its ability to improve the final optimization performance of CC. The CEC'2010 and the CEC'2013 benchmark suites were employed in our experiments. These two suites contain 20 and 15 large scale optimization functions, respectively, all of which are minimization functions of 1000 dimensions except that the functions $f_{13}$ and $f_{14}$ in the CEC'2013 suite have 905 decision variables [46], [47]. Table II presents the types of these functions. It is necessary to point out that, according to the nomenclature in the CEC'2013 suite, the functions $f_{12}$-$f_{14}$ therein and $f_{20}$ in the CEC'2010 suite are overlapping functions. Now it is still not clear what is the best decomposition for this type of functions, and most learning-based decomposition algorithms tend to place all overlapped variable subsets into the same group [17], [22]-[26]. For this reason, we classify this type of functions into nonseparable functions in Table II.

### A. Comparison of Decomposition Performance

In this experiment, DG [13], GDG [22], FII [25], DG2 [26], and RDG [17] were employed as competitors of FDG. As a cla-ssi learning-based decomposition algorithm, DG provides the original decomposition criterion for the other four algorithms. GDG and DG2 can be considered as the same class of decomposition algorithms since they both detect all pairwise interdependencies. FII and RDG share some similarities with FDG because they all perform decomposition from the perspective of variable subsets without collecting the full interdependency information. The main ideas of these five competitors are described in Section II, and their parameters in this experiment were all set according to the suggestions given in the respective original papers.

TABLE II
Types of benchmark functions

| Function type | | CEC'2010 suite | CEC'2013 suite |
|---|---|---|---|
| Fully separable (T$_1$) | | $f_1$ -$f_3$ | $f_1$ -$f_3$ |
| Partially separable (T$_2$) | with separable vars. (T$_{2\text{-}1}$) | $f_4$ -$f_{13}$ | $f_4$ -$f_7$ |
| | without separable vars. (T$_{2\text{-}2}$) | $f_{14}$ -$f_{18}$ | $f_8$ -$f_{11}$ |
| Nonseparable (T$_3$) | with overlapping vars. (T$_{3\text{-}1}$) | $f_{20}$ | $f_{12}$ -$f_{14}$ |
| | without overlapping vars. (T$_{3\text{-}2}$) | $f_{19}$ | $f_{15}$ |

To measure the accuracy of a decomposition, we investigated its consistency with the corresponding ideal decomposition (ID) by an indicator called normalized mutual information (NMI) [48]. NMI can felicitously quantify the similarity between two partitions of a set. For two decompositions $D^1$ and $D^2$, a confusion matrix $M$ can be generated by setting each of its element $M_{ij}$ to the number of common variables in the subset $X_i$ of $D^1$ and the subset $X_j$ of $D^2$, and the NMI between $D^1$ and $D^2$ is defined as

$$\text{NMI}(D^1,D^2) =$$

$$\frac{-2\sum_{i=1}^{k_1}\sum_{j=1}^{k_2} M_{ij}\log_2\left(nM_{ij}\Big/ M_i^1 M_j^2\right)}{\sum_{i=1}^{k_1} M_i^1 \log_2\left(M_i^1\Big/ n\right)+\sum_{j=1}^{k_2} M_j^2 \log_2\left(M_j^2\Big/ n\right)}\times 100\%, \quad (14)$$

where $n$, $k_1$, and $M_i^1$ denote the total number of decision variables in the current instance, the number of variable subsets in $D^1$, and the number of variables in the $i$th subset of $D^1$, respectively. Obviously, the equation $M_i^1 = \sum_{j=1}^{k_2} M_{ij}$ holds. For any two decompositions $D^1$ and $D^2$ of a LSBO instance, $\text{NMI}(D^1,D^2)$ varies in the range $[0,1]$, and the more consistent they are, the larger $\text{NMI}(D^1,D^2)$ is.

Besides the accuracy on all the variables, we also severally calculated the accuracies of a decomposition algorithm on nonseparable variables and separable ones, which were carried out by excluding separable variables and nonseparable ones from a decomposition, respectively. As for the decomposition

TABLE III

Decomposition accuracy and efficiency of DG, GDG, FII, DG2, RDG, and FDG on different types of benchmark functions in the CEC'2010 and the CEC'2013 suites. $\rho_1$, $\rho_2$, and $\rho_3$ denote the NMI indicators on all the variables, nonseparable ones, and separable ones, respectively. If a type of functions do not contain nonseparable or separable variables, the corresponding entries in the table are marked with "–". "Success No." indicates the number of functions on which a decomposition algorithm achieves ideal decomposition.

| Function type (No.) | Statistics | DG $\rho_1$ | DG $\rho_2$ | DG $\rho_3$ | DG feNum | GDG $\rho_1$ | GDG $\rho_2$ | GDG $\rho_3$ | FII $\rho_1$ | FII $\rho_2$ | FII $\rho_3$ | FII feNum | DG2 $\rho_1$ | DG2 $\rho_2$ | DG2 $\rho_3$ | RDG $\rho_1$ | RDG $\rho_2$ | RDG $\rho_3$ | RDG feNum | FDG $\rho_1$ | FDG $\rho_2$ | FDG $\rho_3$ | FDG feNum |
|---|---|---|---|---|---|---|---|---|---|---|---|---|---|---|---|---|---|---|---|---|---|---|---|
| $T_1$ (6) | Median | 100 | – | 100 | 1.00e+06 | 100 | – | 100 | 100 | – | 100 | 3.00e+03 | 100 | – | 100 | 100 | – | 100 | 3.00e+03 | 100 | – | 100 | 5.20e+01 |
| | Mean | 100 | – | 100 | 1.00e+06 | 66.67 | – | 66.67 | 100 | – | 100 | 3.00e+03 | 66.67 | – | 66.67 | 66.67 | – | 66.67 | | 66.67 | – | 66.67 | 5.20e+01 |
| | Std. | 100 | – | 100 | | 47.14 | – | 47.14 | 100 | – | 100 | | 47.17 | – | 47.14 | 47.14 | – | 47.14 | 1.41e+04 | 47.14 | – | 47.14 | 0 |
| | Success No. | 6 | | | | 4 | | | 6 | | | | 4 | | | 4 | | | | 4 | | | |
| $T_{2-1}$ (14) | Median | 98.60 | 96.48 | 100 | 2.70e+05 | 100 | 100 | 100 | 99.66 | 100 | 100 | 4.18e+03 | 100 | 100 | 100 | 100 | 100 | 100 | 1.17e+04 | 100 | 100 | 100 | 3.35e+03 |
| | Mean | 79.05 | 75.05 | 77.03 | 2.99e+05 | 96.40 | 97.16 | 92.86 | 94.70 | 86.35 | 95.62 | 1.17e+04 | 86.36 | 99.30 | 79.85 | 92.14 | 97.92 | 86.28 | 1.40e+04 | 85.82 | 100 | 79.20 | 3.84e+03 |
| | Std. | 24.37 | 36.43 | 27.24 | 3.05e+05 | 11.72 | 7.43 | 25.75 | 8.14 | 26.73 | 8.61 | 2.35e+04 | 26.69 | 2.51 | 38.78 | 19.08 | 7.50 | 33.65 | 1.18e+04 | 28.20 | 0 | 39.87 | 1.92e+03 |
| | Success No. | 5 | | | | 11 | | | 7 | | | | 10 | | | 11 | | | | 11 | | | |
| $T_{2-2}$ (9) | Median | 99.67 | 99.57 | – | 2.10e+04 | 100 | 100 | – | 100 | 100 | – | 2.30e+04 | 100 | 100 | – | 100 | 100 | – | 2.05e+04 | 100 | 100 | – | 6.96e+03 |
| | Mean | 87.00 | 87.00 | – | 2.46e+04 | 84.81 | 84.81 | – | 82.02 | 82.02 | – | 5.82e+04 | 98.48 | 98.48 | – | 86.26 | 86.26 | – | 2.23e+04 | 97.92 | 97.92 | – | 8.16e+03 |
| | Std. | 20.29 | 20.29 | – | 1.13e+04 | 31.03 | 31.03 | – | 32.22 | 32.22 | – | 1.17e+05 | 4.31 | 4.31 | – | 30.92 | 30.92 | – | 1.01e+04 | 5.87 | 5.87 | – | 3.31e+03 |
| | Success No. | 4 | | | | 5 | | | 5 | | | | 8 | | | 6 | | | | 8 | | | |
| $T_{3-1}$ (4) | Median | 0.00 | 0.00 | – | 8.14e+04 | 100 | 100 | – | 100 | 100 | – | 2.54e+05 | 100 | 100 | – | 100 | 100 | – | 3.02e+04 | 100 | 100 | – | 1.17e+04 |
| | Mean | 0.00 | 0.00 | – | 8.09e+04 | 100 | 100 | – | 100 | 100 | – | 2.54e+05 | 100 | 100 | – | 75 | 75 | – | 2.99e+04 | 100 | 100 | – | 1.17e+04 |
| | Std. | 0 | 0 | – | 7.11e+04 | 0 | 0 | – | 0 | 0 | – | 2.49e+05 | 0 | 0 | – | 43.30 | 43.30 | – | 2.07e+04 | 0 | 0 | – | 6.32e+03 |
| | Success No. | 0 | | | | 4 | | | 4 | | | | 4 | | | 3 | | | | 4 | | | |
| $T_{3-2}$ (2) | Median | 100 | 100 | – | 2.00e+03 | 100 | 100 | – | 100 | 100 | – | 4.00e+03 | 100 | 100 | – | 100 | 100 | – | 6.08e+03 | 100 | 100 | – | 5.20e+01 |
| | Mean | 100 | 100 | – | 2.00e+03 | 100 | 100 | – | 100 | 100 | – | 4.00e+03 | 100 | 100 | – | 100 | 100 | – | 6.08e+03 | 100 | 100 | – | 5.20e+01 |
| | Std. | 0 | 0 | – | 0 | 0 | 0 | – | 0 | 0 | – | 0 | 0 | 0 | – | 0 | 0 | – | 8.00e+01 | 0 | 0 | – | 0 |
| | Success No. | 2 | | | | 2 | | | 2 | | | | 2 | | | 2 | | | | 2 | | | |
| | Total Success No. | 17 | | | | 26 | | | 24 | | | | 28 | | | 26 | | | | 29 | | | |

efficiency, it was directly measured by the number of FEs consumed during the decomposition process. Table S-I in the supplementary document provides detailed results obtained by the six decomposition algorithms, and Table III summarizes these results in terms of function type. It is worth mentioning that the FEs of GDG and DG2 are not listed out due to the page limit and the fact that these two algorithms invariably consume 501501 and 500501 FEs, respectively, on different types of functions of 1000 dimensions [26]. From Tables S-I and III, it can be observed that FDG holds significant superiority over its five competitors since it achieves the highest decomposition accuracy with much fewer FEs on most types of functions.

For each of the six fully separable functions ($f_1$-$f_3$ in both CEC'2010 and CEC'2013 suites), FDG completes decomposition by just executing ITIP and consumes only 52 FEs, which are negligible compared with the ones consumed by other algorithms. It is notable that FDG seems to get wrong decompositions for $f_3$ in both suites. In fact, both functions are shifted Ackley function and are actually additively nonseparable [46]. In this sense, FDG makes correct identifications. The same case is encountered by GDG, DG2 and RDG, while DG and FII avoid this by taking much larger thresholds.

The functions $f_4$-$f_{18}$ in the CEC'2010 suite and $f_4$-$f_{11}$ in the CEC'2013 suite are partially separable functions. On this type of functions, FDG also shows satisfying performance. It achieves 100% decomposition accuracy on 19 out of total 23 functions with the other four ones being $f_6$ and $f_{11}$ in the CEC'2010 suite and $f_6$ and $f_8$ in the CEC'2013 suite. The reason for its imperfect accuracy on the first three functions lies in that these functions all take Ackley function as a separable subfunction, while FDG tends to regard it as additively nonseparable. This can be revealed by the corresponding 100% accuracy on nonseparable variables and extremely low accuracy on separable ones. The function $f_8$ in the CEC'2013 suite is the only function on which FDG really makes a mistake. By investigating its composition, it can be found that this function involves 20 nonseparable variable subsets and most subsets have sharply different contribution weights to the original objective function. Then the roundoff errors related to the variables of small contributions tend to be overestimated. As a consequence, the SVEP in FDG wrongly judges some nonseparable variables of this kind to be separable and tends to generate a very small separability threshold, which further leads to the wrong combination of some actually separable variables.

Among the other five decomposition algorithms, DG2 achieves almost the same accuracy with FDG, but fails on one more function, i.e., $f_7$ in the CEC'2013 suite. FII improves the accuracy of DG on some functions, such as $f_{13}$ and $f_{18}$ in the CEC'2010 suite, by taking indirect interdependency into account, but still fails on 11 out of 23 partially separable functions due to the same decomposition threshold on different functions. GDG and RDG enhance the adaptivity of DG and FII to a certain extent by generating different thresholds for different functions. However, their performance is obviously inferior to that of DG2 and FDG on the imbalanced functions, such as $f_{10}$ and $f_{11}$ in the CEC'2013 suite.

As for the decomposition efficiency, Table III clearly demonstrates that FDG saves about two-thirds of FEs on average on partially separable functions when compared with the fastest existing algorithms, i.e., FII for functions of type $T_{2-1}$ and RDG for functions of type $T_{2-2}$. In comparison with GDG and DG2, FDG consumes less than 2% FEs. The high efficiency of FDG on partially separable functions profits from its two algorithmic components, i.e., SVEP and BTDP. To further verify the effectiveness of these two components, we specially did an experiment by removing SVEP from FDG. Fig. S-4 in the supplementary document presents the variation of the number of FEs consumed by FDG. It can be found that without SVEP, the FE consumption of FDG increases to a certain extent, but is still much less than those of FII and RDG.

When it comes to the four overlapping functions ($f_{20}$ in the CEC'2010 suite and $f_{12}$-$f_{14}$ in the CEC'2013 suite), the ITIP in FDG judges them to be partially separable since it does not consider the particularity of overlap. Despite of this, following procedures in FDG still enable it to capture right decompositions with much fewer FEs than existing algorithms. As for $f_{19}$ in the CEC'2010 suite and $f_{15}$ in the CEC'2013 suite, FDG correctly identifies them as nonseparable functions, and thus completes the decomposition with negligible FEs.

### B. Comparison of Optimization Performance

*1) Comparison with other decomposition algorithms under DECC framework:* To investigate FDG's capability in enhancing the performance of CC, we embedded it into the canonical CC framework, DECC [12], [13], and compared its final optimization results with the ones obtained by DG [13], GDG [22], FII [25], DG2 [26], RDG [17], and ID. DECC employs a well known variant of differential evolution called SaNSDE [49] as optimizer and optimizes all the subproblems in a round-robin fashion. In this experiment, the parameters of DECC and SaNSDE were strictly set according to their original papers. As suggested by [46], each CC algorithm was allowed to exhaust a maximum number of $3.0 \times 10^6$ FEs in each run, and the result on each function was calculated in terms of the median, the mean, and the standard deviation of the best solutions obtained over 25 independent runs. For each decomposition algorithm, the number of FEs consumed during the decomposition process was also counted into the total number of allowed FEs.

Table S-II in the supplementary document presents the final experimental results obtained by the seven DECC algorithms, where the differences between the results of FDG and the ones of the other six algorithms are detected using a two-tailed Wilcoxon rank-sum test at a significance level of 0.05. Note that we consider there is no difference between two solutions if the magnitude orders of their fitness values are not greater than −10. Table IV summarizes the numbers of wins, ties, and losses of FDG against its six competitors based on the results in Table S-II. From these results, it can be observed that FDG exhibits very satisfying performance. It performs no worse than ID on the total 35 functions. It is rather unexpected to see that FDG even achieves higher solution quality in terms of median than ID on $f_6$ in the CEC'2010 suite and $f_8$ in the CEC'2013 suite. Compared with DG, GDG, FII, DG2, and RDG, FDG yields much better (or similar) solutions for 18, 25, 12, 24, and 5 (or 17, 10, 22, 10, and 30) functions, respectively.

As the seven decomposition algorithms are integrated into the same CC framework, the differences among their final optimization results mainly stem from their differences in the decomposition accuracy and efficiency. It is understandable that the higher the decomposition accuracy is and the more FEs left for optimization there are, the better the final optimization result will be. This general conclusion can be verified by most of the results shown in Tables S-I and S-II.

TABLE IV
THE NUMBERS OF WINS, TIES, AND LOSSES (W/T/L) OF FDG AGAINST DG, GDG, FII, DG2, RDG, AND ID ON DIFFERENT TYPES OF FUNCTIONS.

| Fun. type (No.) | DG | GDG | FII | DG2 | RDG | ID |
|---|---|---|---|---|---|---|
| $T_1$ (6) | 3/ 3/0 | 3/ 3/0 | 0/ 6/0 | 3/ 3/0 | 0/ 6/0 | 0/ 6/0 |
| $T_{2\text{-}1}$ (14) | 8/ 6/0 | 10/ 4/0 | 6/ 7/1 | 8/ 5/1 | 2/12/0 | 1/13/0 |
| $T_{2\text{-}2}$ (9) | 3/ 6/0 | 6/ 3/0 | 4/ 5/0 | 7/ 2/0 | 2/ 7/0 | 1/ 8/0 |
| $T_{3\text{-}1}$ (4) | 4/ 0/0 | 4/ 0/0 | 2/ 2/0 | 4/ 0/0 | 1/ 3/0 | 0/ 4/0 |
| $T_{3\text{-}2}$ (2) | 0/ 2/0 | 2/ 0/0 | 0/ 2/0 | 2/ 0/0 | 0/ 2/0 | 0/ 2/0 |
| Total | 18/17/0 | 25/10/0 | 12/22/1 | 24/10/1 | 5/30/0 | 2/33/0 |

Nevertheless, there are some exceptions which need to be discussed. As indicated in the last subsection, GDG, DG2, RDG, and FDG judge $f_3$ in both suites to be nonseparable, and thus achieve a decomposition accuracy of 0%. However, they obtain almost the same optimization result with DG and FII, which correctly identify the two functions to be fully separable. The reason is that DECC treats all the separable variables as a whole since the optimal division of separable variables remains so far an open problem [26], [30]. For $f_6$ in the CEC'2010 suite and $f_8$ in the CEC'2013 suite, FDG achieves lower decomposition accuracies than most of its competitors, but helps DECC find better solutions, which is somewhat counterintuitive. Undoubtedly, the lower FE cost contributes a lot to the success of FDG on these two functions. However, there are some other reasons. When it comes to $f_6$, a closer observation reveals that FDG divides its 950 separable variables into three variable subsets, which exactly reduces the scales of its subfunctions and thus its solving difficulty. As for $f_8$, the analysis in the last subsection indicates that FDG merges many variables of smaller contributions together. Consequently, the total number of subfunctions decreases and the subfunctions of larger contributions can be assigned with more FEs, which is beneficial to finding better solutions.

It should be pointed out that although FDG consumes much fewer FEs than all the other decomposition algorithms (except ID which is implemented manually without FE), it fails to achieve better solutions for part of functions, such as $f_5$ in the CEC'2010 suite and $f_9$ in the CEC'2013 suite. This is because that the allowed FEs ($3.0 \times 10^6$) are much more than the ones required by decomposition algorithms and the optimizer of DECC, namely SaNSDE, tends to prematurely converge before exhausting remaining FEs. To reveal the performance differences between FDG and its six competitors more comprehensively, we also recorded their intermediate results obtained with $6.0 \times 10^5$ and $1.2 \times 10^5$ FEs according to the suggestions in [46]. The detailed results are reported in Tables S-III and S-IV in the supplementary document. Fig. S-5 further presents the numbers of wins, ties, and losses of FDG against its competitors under the cases of three different quantities of FEs. It can be seen that, with the reduction of the number of available FEs, the advantage of FDG becomes more obvious. When $6.0 \times 10^5$ FEs are allowed, FDG performs no worse than all of its six competitors except for being surpassed by ID and FII on $f_{18}$ in the CEC'2010 suite and $f_8$ in the CEC'2013 suite, respectively. When the number of FEs reduces to $1.2 \times 10^5$, it outperforms DG and RDG on more functions.

*2) Comparison with state-of-the-art algorithms for LSBO:* To deeply test the effectiveness of FDG, we integrated it into the CC framework with a fine-grained computation resource allocation strategy (FCRACC) [14]. As a contribution-based CC framework, FCRACC allocates available FEs among all the subproblems according to their real-time contributions to the overall fitness improvement rather than in a round-robin fashion. Besides, it takes an excellent variant of differential evolution called SHADE [50] as optimizer. We compared the performance of the resulting algorithm, FCRACC$_{FDG+SHADE}$, with that of four state-of-the-art algorithms for LSBO, including the memetic algorithm that creates local search chains with Solis Wets' method (MA-SW-Chains) [51], the multiple offspring sampling framework for the CEC'2013 competition (MOS-CEC2013) [52], the segment-based predominant learning swarm optimizer with dynamic segment number (DSPLSO) [53], and the level-based learning swarm optimizer with dynamic level number (DLLSO) [54]. MA-SW-Chains was ranked first on the CEC'2010 competition on large scale global optimization, while MOS-CEC2013 won the CEC'2013 and CEC'2015 competitions. As two excellent algorithms recently developed for LSBO, DSPLSO and DLLSO expands the ideas of the competitive swarm optimizer [55] and the social learning PSO [56], respectively, and were shown to be more effective than some CC algorithms.

Table S-V in the supplementary document reports the experimental results of the five algorithms. To ensure the comparison fair, the results of MA-SW-Chains and MOS-CEC2013 are directly taken from [57], which is a review paper published by the proposer of MOS-CEC2013, and the results of DSPLSO and DLLSO are taken from their original papers. It is notable that all these results were generated according to the default experimental settings suggested by [46] except that the ones of DSPLSO and DLLSO were yielded based on 30 independent runs. To measure the performance of FCRACC$_{FDG+SHADE}$ against each of its four competitors, we quantified the difference between their fitness means with Cohen's $d$ effect size [58]. Compared with significance test such as $t$-test, Cohen's $d$ effect size is affected less by the sample size. It is generally considered that there is no significant difference between two means if the absolute value of their Cohen's $d$ effect size is less than 0.2. Table V summarizes the numbers of wins, ties, and losses of FDG against its four competitors based on the results in Table S-V. Besides, it also reports the rankings of the five algorithms obtained through the Friedman test.

These results clearly show that FCRACC$_{FDG+SHADE}$ not only outperforms MA-SW-Chains, DSPLSO, and DLLSO, but also achieves competitive performance with MOS-CEC2013. It dominates these four algorithms on 20, 23, 22, and 16 out of total 35 functions, respectively. Moreover, it also shows similar performance on the other 4, 4, 2, and 4 functions in comparison with the above four algorithms, respectively. Especially on functions such as $f_4$, $f_8$, $f_9$, $f_{13}$ in the CEC'2010 suite and $f_2$, $f_4$, $f_8$ in the CEC'2013 suite, FCRACC$_{FDG+SHADE}$ yields better solution quality than its four competitors by several orders of magnitude in terms of fitness mean. The Friedman test results listed in Table V indicate that FCRACC$_{FDG+SHADE}$ can be ranked almost the same with MOS-CEC2013, followed by MA-SW-Chains, DLLSO, and DSPLSO.

TABLE V
STATISTICAL RESULTS OF MA-SW-CHAINS, MOS-CEC2013, DSPLSO, DLLSO, AND FCRACC$_{FDG+SHADE}$ ON TWO BENCHMARK SUITES.

| Functions | Indicator | MA-SW-Chains | MOS-CEC2013 | DSPLSO | DLLSO | FCRA$_{FDG+SHADE}$ |
|---|---|---|---|---|---|---|
| All the 35 functions | w/t/l | 20/4/11 | 16/4/15 | 23/4/8 | 22/2/11 | – |
| | Ranking | 3.0143 | 4.2286 | 3.5517 | 3.5143 | 2.5000 |
| 29 (partially) separable functions | w/t/l | 17/3/9 | 15/4/10 | 20/3/6 | 18/2/9 | – |
| | Ranking | 3.0345 | 2.7241 | 3.3448 | 3.4655 | 2.4310 |

A closer observation on Table S-V further reveals that FCRACC$_{FDG+SHADE}$ tends to lose its superiority on nonseparable functions (such as $f_{20}$ in the CEC'2010 suite and $f_{12}$ in the CEC'2013 suite) and the functions where it gets improper decomposition (such as $f_6$ in both benchmark suites). Both kinds of functions are of large scales or contain large scale subfunctions. This indicates the optimizer of FCRACC$_{FDG+SHADE}$ is not so efficient on LSBO problems. If we exclude the six nonseparable functions, the superiority of FCRACC$_{FDG+SHADE}$ can be further highlighted. It surpasses MOS-CEC2013 as well as its other three competitors, which can be verified by the summarized results given in the last two rows of Table V. From above analysis, it can be concluded that the superiority of FCRACC$_{FDG+SHADE}$ mainly benefits from its CC framework and decomposition algorithm.

The other four algorithms also demonstrate their specialities. MA-SW-Chains performs best on $f_3$, $f_6$, and $f_{10}$ in the CEC'2013 suite, and MOS-CEC2013 outperforms its competitors on nine out of total 35 functions, over a half of which are nonseparable ones. DSPLSO yields best solution quality on $f_5$, $f_6$, $f_{11}$ in the CEC'2010 suite and $f_5$, $f_9$ in the CEC'2013 suite. As for DLLSO, it has an edge over the other four algorithms when used to solve $f_{10}$ and $f_{15}$ in the CEC'2010 suite.

## V. CONCLUSION AND FUTURE WORK

This paper presents a fast learning-based decomposition algorithm named FDG for LSBO problems. Different from existing decomposition algorithms, FDG employs a normalized interdependency indicator, which enables it to analyze the distribution of the indicator values for different pairs of variable subsets. Benefiting from this characteristic and equipped with a sophisticated variable subset selection strategy, FDG can identify the type of an LSBO instance at a negligible computation cost and adaptively generate decomposition thresholds for partially separable instances. Furthermore, FDG adopts a binary-tree-based decomposition procedure, which can significantly reduce separability detection times and greatly enhance the reutilization degree of evaluated solutions. Consequently, FDG is able to decompose an $n$-dimensional partially separable instance with $O(n\log_2 k)$ FEs in the worst case, where $k$ is the number of nonseparable variable subsets. Experimental results on two benchmark suites show that FDG can generate highly accurate decompositions with much fewer FFs than five popular learning-based algorithms. When embedded into CC frameworks, it can outperform these existing algorithms in terms of solution quality, and even shows clear superiority over four state-of-the-art optimization algorithms for LSBO.

At present, FDG only applies to additively separable LSBO problems like most of the other learning-based decomposition algorithms. In the future, we will extend the idea to generally separable problems. Besides, it is necessary to develop decomposition algorithms that can adaptively determine its decomposition granularity according to available computation resources, since a more fine-grained decomposition may not necessarily lead to a better optimization result.

REFERENCES

[1] S. Mahdavi, M. E. Shiri, and S. Rahnamayan, "Metaheuristics in large-scale global continues optimization: A survey," *Inf. Sci.*, vol. 295, pp. 407–428, Feb. 2015.

[2] M. Bhattacharya, R. Islam, and J. Abawajy. "Evolutionary optimization: A big data perspective," *Journal of Network & Computer Applications*, vol. 59, pp. 416–426, Jan. 2016.

[3] G. A. Trunfio, "Metaheuristics for continuous optimization of high dimensional problems: State of the art and perspectives," in *Big Data Optimization: Recent Developments and Challenges*, A. Emrouznejad, Ed. Cham, Switzerland: Springer, 2016, pp. 437–460.

[4] N. R. Sabar, J. Abawajy, and J. Yearwood, "Heterogeneous cooperative co-evolution memetic differential evolution algorithm for big data optimization problems," *IEEE Trans. Evol. Comput.*, vol. 21, no. 2, pp. 315–327, Apr., 2017.

[5] X. Ma, X. Li, Q. Zhang, K. Tang, Z. Liang, W. Xie, and Z. Zhu, "A survey on cooperative co-evolutionary algorithms," *IEEE Trans. Evol. Comput.*, 2018. in press, doi: 10.1109/TEVC.2018. 2868770.

[6] T. Chai, Y. Jin, and B. Sendhoff, "Evolutionary complex engineering optimization: Opportunities and challenges," *IEEE Comput. Intell. Mag.*, vol. 8, no. 3, pp. 12–15, Aug. 2013.

[7] S. A. Thomas and Y. Jin, "Reconstructing biological gene regulatory networks: Where optimization meets big data," *Evol. Intell.*, vol. 7, no. 1, pp. 29–47, Apr. 2014.

[8] Y. Liu, X. Yao, Q. Zhao, and T. Higuchi, "Scaling up fast evolutionary programming with cooperative coevolution," in *Proc. IEEE Congr. Evol.*, 2001, pp. 1101–1108.

[9] F. van den Bergh and A. P. Engelbrecht, "A cooperative approach to particle swarm optimization," *IEEE Trans. Evol. Comput.*, vol. 8, no. 3, pp. 225–239, 2004.

[10] M. A. Potter and K. A. De Jong, "A cooperative coevolutionary approach to function optimization," in *Proc. Int. Conf. Parallel Problem Solving Nat.*, 1994, pp. 249–257.

[11] M. A. Potter and K. A. De Jong, "Cooperative coevolution: An architecture for evolving coadapted subcomponents," *Evol. Comput.*, vol. 8, no. 1, pp. 1–29, Mar. 2000.

[12] Z. Yang, K. Tang, and X. Yao, "Large scale evolutionary optimization using cooperative coevolution," *Inf. Sci.*, vol. 178, pp. 2986–2999, Aug. 2008.

[13] M. N. Omidvar, X. Li, Y. Mei, and X. Yao, "Cooperative co-evolution with differential grouping for large scale optimization," *IEEE Trans. Evol. Comput.*, vol. 18, no. 3, pp. 378–393, 2014.

[14] Z. Ren, Y. Liang, A. Zhang, Y. Yang, Z. Feng, and L. Wang, "Boosting cooperative coevolution for large scale optimization with a fine-grained computation resource allocation strategy," *IEEE Trans. Cybern.*, 2018. in press, doi: 10.1109/TCYB.2018.2859635.

[15] M. Yang, M. N. Omidvar, and C. Li et al., "Efficient resource allocation in cooperative co-evolution for large-scale global optimization," *IEEE Trans. Evol. Comput.*, vol. 21, no. 4, pp. 493–505, Aug. 2017.

[16] P. Yang, K. Tang, and X. Yao. "Turning high-dimensional optimization into computationally expensive optimization," *IEEE Trans. Evol. Comput.*, vol. 22, no. 1: 143–156, Feb. 2018.

[17] Y. Sun, M. Kirley, and S.K. Halgamuge, "A recursive decomposition method for large scale continuous optimization," *IEEE Trans. Evol. Comput.*, vol. 22, no. 5: 647–661, Oct. 2018.

[18] X. Peng, Y. Jin, and H. Wang, "Multimodal optimization enhanced cooperative coevolution for large-scale optimization," *IEEE Trans. Cybern.*, 2018. in press, doi: 10.1109/TCYB.2018.2846179.

[19] M. N. Omidvar, X. Li, Z. Yang, and X. Yao, "Cooperative co-evolution for large scale optimization through more frequent random grouping," in *Proc. IEEE Congr. Evol.*, 2010, pp. 1754–1761.

[20] W. Chen, T. Weise, Z. Yang, and K. Tang, "Large-scale global optimization using cooperative coevolution with variable interaction learning," In *Proc. Int. Conf. Parallel Problem Solving Nat.*, 2010, pp. 300–309.

[21] Y. Sun, M. Kirley, and S. K. Halgamuge, "Quantifying variable interactions in continuous optimization problems," *IEEE Trans. Evol. Comput.*, vol. 21, no. 2, pp. 249–264, Apr. 2017.

[22] Y. Mei, M. Omidvar, X. Li, and X. Yao, "A competitive divide-and-conquer algorithm for unconstrained large-scale black-box optimization," *ACM Trans. Math. Softw.*, vol. 42, no. 2, p. 13, Jun. 2016.

[23] Y. Sun, M. Kirley, and S.K. Halgamuge, "Extended differential grouping for large scale global optimization with direct and indirect variable interactions," in *Proc. Genet. Evol. Comput. Conf.*, 2015, pp. 313–320.

[24] Y. Ling, H. Li, and B. Cao, "Cooperative co-evolution with graph-based differential grouping for large scale global optimization," in *Proc. Int. Conf. Natural Computation, Fuzzy Systems and Knowledge Discovery*, Aug. 2016, pp. 95–102.

[25] X.-M. Hu, F.-L. He, W.-N. Chen, and J. Zhang, "Cooperation coevolution with fast interdependency identification for large scale optimization," *Inf. Sci.*, vol. 381, pp. 142–160, Mar. 2017.

[26] M. N. Omidvar, M. Yang, Y. Mei, X. Li, and X. Yao, "DG2: A faster and more accurate differential grouping for large-scale black-box optimization," *IEEE Trans. Evol. Comput.*, vol. 21, no. 6: 929–942, Dec. 2017.

[27] Z. Ren, B. Pang, M. Wang, Z. Feng, Y. Liang, and Y. Zhang, "Surrogate model assisted cooperative coevolution for large scale optimization," *Appl. Intell.*, vol. 49, no. 2: 513–531, Feb. 2019.

[28] Y. Shi, H. Teng, and Z. Li, "Cooperative co-evolutionary differential evolution for function optimization," in *Proc. Int. Conf. Nat. Comput.*, 2005, pp. 1080–1088.

[29] Z. Yang, K. Tang, and X. Yao, "Multilevel cooperative coevolution for large scale optimization," in *Proc. IEEE Congr. Evol.*, June 2008, pp. 1663–1670.

[30] M. N. Omidvar, Y. Mei, and X. Li. "Effective decomposition of large scale separable continuous functions for cooperative co-evolutionary algorithms," In *Proc. IEEE Congr. Evol.*, 2014, pp. 1305–1312.

[31] X. Li and X. Yao, "Cooperatively coevolving particle swarms for large scale optimization," *IEEE Trans. Evol. Comput.*, vol. 16, no. 2, pp. 210–224, Apr. 2012.

[32] G. A. Trunfio. A cooperative coevolutionary differential evolution algorithm with adaptive subcomponents. *Procedia Computer Science*, vol. 51, no. 1, pp. 834–844, 2015.

[33] G. A. Trunfio, P. Topa, and J. Was. A new algorithm for adapting the configuration of subcomponents in large-scale optimization with cooperative coevolution. *Inf. Sci.*, vol. 372, pp. 773–795, 2016.

[34] T. Ray and X. Yao, "A cooperative coevolutionary algorithm with correlation based adaptive variable partitioning," in *Proc. IEEE Congr. Evol.*, May 2009, pp. 983–989.

[35] W. Dong, T. Chen, P. Tino, and X. Yao, "Scaling up estimation of distribution algorithms for continuous optimization," *IEEE Trans. Evol. Comput.*, vol. 17, no. 6, pp. 797–822, 2013.

[36] Q. Xu, M. L. Sanyang, and A. Kaban, "Large scale continuous EDA using mutual information," in *Proc. IEEE Congr. Evol.*, 2016, pp. 3718–3725.

[37] K. Weicker and N. Weicker, "On the improvement of coevolutionary optimizers by learning variable interdependencies," in *Proc. IEEE Congr. Evol.*, 1999, pp. 1627–1632.

[38] L. Sun, S. Yoshida, X. Cheng, and Y. Liang, "A cooperative particle swarm optimizer with statistical variable interdependence learning," *Inf. Sci.*, vol. 186, pp. 20–39, Mar. 2012.

[39] H. Ge, L. Sun, X. Yang, S. Yoshida, and Y. Liang, "Cooperative differential evolution with fast variable interdependence learning and cross-cluster mutation," *Appl. Soft Comput.*, vol. 36, pp. 300–314, Nov. 2015.

[40] H. Ge, L. Sun, G. Tan, Z. Chen, and C. L. Philip Chen, "Cooperative hierarchical PSO with two stage variable interaction reconstruction for large scale optimization," *IEEE Trans. Cybern.*, vol. 47, no. 9, pp. 2809–2823, Sep. 2017.

[41] M. N. Omidvar, X. Li, and X. Yao, "Cooperative co-evolution with delta grouping for large scale non-separable function optimization," in *Proc. IEEE Congr. Evol.*, 2010, pp. 1762–1769.

[42] Y. Rojas and R. Landa, "Towards the use of statistical information and differential evolution for large scale global optimization," In *Proc. Int. Conf. Electrical Engineering Computing Science and Automatic Control*, 2011, pp. 1–6.

[43] M. Tezuka, M. Munetomo, and K. Akama, "Linkage identification by nonlinearity check for real-coded genetic algorithms," In *Proc. Genet. Evol. Comput. Conf.*, 2004, pp. 222–233.

[44] Y. Sun, M. N. Omidvar, M. Kirley, and X. Li, "Adaptive threshold parameter estimation with recursive differential grouping for problem decomposition," in *Proc. Genet. Evol. Comput. Conf.*, 2018, pp. 889–896.

[45] "IEEE standard for floating-point arithmetic, IEEE std. 754-2008," Aug. 2008.

[46] K. Tang, X. Li, P. N. Suganthan, Z. Yang, and T. Weise, "Benchmark functions for the CEC'2010 special session and competition on large scale global optimization," Nat. Inspired Comput. Appl. Lab., Univ. Sci. Technol. China, Hefei, China, Rep., 2010.

[47] X. Li, K. Tang, M. N. Omidvar, Z. Yang, and K. Qin, "Benchmark functions for the CEC'2013 special session and competition on large scale global optimization," Evol. Comput. Mach. Learn. Group, RMIT Univ., Melbourne, VIC, Australia, Rep., 2013.

[48] L. Danon, J. Duch, A. Diaz-Guilera, and A. Arenas, "Comparing community structure identification," *Journal of Statistical Mechanics: Theory and Experiment*, vol. 2005, p. 09008, Sep. 2005.

[49] Z. Yang, K. Tang, and X. Yao, "Self-adaptive differential evolution with neighborhood search," in *Proc. IEEE Congr. Evol. Comput.*, 2008, pp. 1110–1116.

[50] R. Tanabe and A. Fukunaga, "Success-history based parameter adaptation for differential evolution," in *Proc. IEEE Congr. Evol. Comput.*, Jun. 2013, pp. 71–78.

[51] D. Molina, M. Lozano, and F. Herrera, "MA-SW-Chains: Memetic algorithm based on local search chains for large scale continuous global optimization," in *Proc. IEEE Congr. Evol. Comput.*, 2010, pp. 1–8.

[52] A. LaTorre, S. Muelas, and J.-M. Peña, "Large scale global optimization: Experimental results with MOS-based hybrid algorithms," in *Proc. IEEE Congr. Evol. Comput.*, 2013, pp. 2742–2749.

[53] Q. Yang, W. N. Chen, T. Gu, H. Zhang, J. D. Deng, Y. Li, and J. Zhang, "Segment-Based Predominant Learning Swarm Optimizer for Large-Scale Optimization," *IEEE Trans. Cybern.,* vol. 47, no. 9, pp. 2896–2910, Sep. 2017.

[54] Q. Yang, W. N. Chen, J. D. Deng, Y. Li, T. Gu, and J. Zhang, "A level-based learning swarm optimizer for large-scale optimization," *IEEE Trans. Evol. Comput.,* vol. 22, no. 4, pp. 578–594, Aug. 2018.

[55] R. Cheng and Y. Jin, "A competitive swarm optimizer for large scale optimization," *IEEE Trans. Cybern.,* vol. 45, no. 2, pp. 191–204, Feb. 2015.

[56] R. Cheng and Y. Jin, "A social learning particle swarm optimization algorithm for scalable optimization," *Inf. Sci.*, vol. 291, pp. 43–60, Jan. 2015.

[57] A. LaTorre, S. Muelas, and J.-M. Peña, "A comprehensive comparison of large scale global optimizers," *Inf. Sci.*, vol. 316, pp. 517–549, Sep. 2015.

[58] J. Cohen, *Statistical Power Analysis for the Behavioral Sciences*, 2nd ed. Hillsdale, NJ, USA: Lawrence Erlbaum, 1988.

# Supplementary Materials for "A Fast Differential Grouping Algorithm for Large Scale Black-Box Optimization"

Zhigang Ren, *Member, IEEE*, An Chen, Yaochu Jin, *Fellow, IEEE*, Wenhua Guo, Yongsheng Liang, *Student Member, IEEE*, and Zuren Feng, *Member, IEEE*

## S-I Theorems and Proofs

**Theorem 1.** *For an n-dimensional separable function $f(\mathbf{x})$ with $X$ being its variable set, if its two disjoint variable subsets $X_i$ and $X_j$ are separable from each other, then for $\forall \mathbf{cv} \in [\mathbf{lb}, \mathbf{ub}]^n$, $\forall \mathbf{x}_i^0, \mathbf{x}_i^{'} \in [\mathbf{lb}_i, \mathbf{ub}_i]^{|X_i|}$, and $\forall \mathbf{x}_j^0, \mathbf{x}_j^{'} \in [\mathbf{lb}_j, \mathbf{ub}_j]^{|X_j|}$, the following equation holds:*

$$\Delta(\mathbf{x}_i^0, \mathbf{x}_i^{'} \mid \mathbf{x}_j^0) = \Delta(\mathbf{x}_i^0, \mathbf{x}_i^{'} \mid \mathbf{x}_j^{'}) , \qquad (S\text{-}1)$$

*where $\mathbf{lb}$ and $\mathbf{ub}$ denote the lower and the upper bounds of $\mathbf{x}$, respectively, $\mathbf{lb}_i = \mathbf{lb}(X_i)$ denotes the lower bound of the subcomponent $\mathbf{x}_i = \mathbf{x}(X_i)$, and $\Delta(\mathbf{x}_i^0, \mathbf{x}_i^{'} \mid \mathbf{x}_j^0)$ is defined as*

$$\Delta(\mathbf{x}_i^0, \mathbf{x}_i^{'} \mid \mathbf{x}_j^0) \triangleq f(\mathbf{cv} \leftarrow \mid \mathbf{x}_i^0, \mathbf{x}_j^0) - f(\mathbf{cv} \leftarrow \mid \mathbf{x}_i^{'}, \mathbf{x}_j^0) \qquad (S\text{-}2)$$

*with $\mathbf{cv} \leftarrow \mid \mathbf{x}_i^0, \mathbf{x}_j^0$ denoting the complete solution obtained by inserting $\mathbf{x}_i^0$ and $\mathbf{x}_j^0$ into the corresponding positions in $\mathbf{cv}$.*

*Proof.* If $X_i$ and $X_j$ are separable from each other, then according to **Definition 2** given in the main body of this paper, there necessarily exist two disjoint subsets $X_i^+ \supseteq X_i$ and $X_j^+ \supseteq X_j$ satisfying

$$f(\mathbf{x}) = f_i(\mathbf{x}_i^+) + f_j(\mathbf{x}_j^+) ,$$

where $\mathbf{x}_i^+ = \mathbf{x}(X_i^+)$ is a $|X_i^+|$-dimensional subvector composed by the decision variables in the subset $X_i^+$ and $f_i(\cdot)$ is the corresponding subfunction. Then

$$\Delta(\mathbf{x}_i^0, \mathbf{x}_i^{'} \mid \mathbf{x}_j^0)$$
$$= f(\mathbf{cv} \leftarrow \mid \mathbf{x}_i^0, \mathbf{x}_j^0) - f(\mathbf{cv} \leftarrow \mid \mathbf{x}_i^{'}, \mathbf{x}_j^0)$$
$$= f(\dots, \mathbf{x}_i^0, \dots, \mathbf{x}_j^0, \dots) - f(\dots, \mathbf{x}_i^{'}, \dots, \mathbf{x}_j^0, \dots)$$
$$= f_i(\dots, \mathbf{x}_i^0, \dots) + f_j(\dots, \mathbf{x}_j^0, \dots) - f_i(\dots, \mathbf{x}_i^{'}, \dots) - f_j(\dots, \mathbf{x}_j^0, \dots, )$$

$$= f_i(\dots, \mathbf{x}_i^0, \dots) + f_j(\dots, \mathbf{x}_j^{'}, \dots) - f_i(\dots, \mathbf{x}_i^{'}, \dots) - f_j(\dots, \mathbf{x}_j^{'}, \dots, )$$
$$= f(\dots, \mathbf{x}_i^0, \dots, \mathbf{x}_j^{'}, \dots) - f(\dots, \mathbf{x}_i^{'}, \dots, \mathbf{x}_j^{'}, \dots)$$
$$= f(\mathbf{cv} \leftarrow \mid \mathbf{x}_i^0, \mathbf{x}_j^{'}) - f(\mathbf{cv} \leftarrow \mid \mathbf{x}_i^{'}, \mathbf{x}_j^{'})$$
$$= \Delta(\mathbf{x}_i^0, \mathbf{x}_i^{'} \mid \mathbf{x}_j^{'}) . \qquad \square$$

**Corollary 1.** *Let $X$ be the variable set of a function $f(\mathbf{x})$ and $X_i$ be a proper subset of $X$, then $X_i$ and $X_{-i} (= X \setminus X_i)$ are separable from each other if and only if they satisfy (S-1).*

*Proof.* The necessity of (S-1) for the separability between $X_i$ and $X_{-i}$ has been verified in Theorem 1. Here we focus on the sufficiency of (S-1).

If $X_i$ and $X_{-i}$ satisfy (S-1), then any variable in $X_{-i}$ does not affect the variation of $f(\mathbf{x})$ with respect to any variable in $X_i$. This implies that there are no other operations except addition or subtraction between the variables in $X_i$ and the ones in $X_{-i}$. On the other side, there do not exist other variables that can link the two complementary variable subsets together. Then we can independently define subfunctions for $X_i$ and $X_{-i}$, and the operation between the two subfunctions is additive. According to this conclusion and Definition 2 given in the main body of this paper, $X_i$ and $X_{-i}$ can be separated from each other. $\square$

**Theorem 2.** *Assume that $f(\mathbf{x})$ is a partially separable function of $n$ dimensions, it involves $k$ disjoint sets of nonseparable variables with each set containing $s_i$ ($s_i \geq 2$, $i = 1, \dots, k$) variables, and $\sum_{i=1}^{k} s_i \leq n$ holds. If we randomly partition the whole variable set into two (nearly) equal-sized subsets, then these two subsets are separable from each other with a probability of*

$$P_{s1} = \frac{C_{n_s}^{\lfloor n/2 \rfloor} + \sum_{i=1}^{k} C_{n_s}^{\lfloor n/2 \rfloor - s_i} + \sum_{i=1}^{k} \sum_{j=1, j \neq i}^{k} C_{n_s}^{\lfloor n/2 \rfloor - s_i - s_j} + \dots + C_{n_s}^{\lfloor n/2 \rfloor - \sum_{i=1}^{k} s_i}}{C_n^{\lfloor n/2 \rfloor}} , \qquad (S\text{-}3)$$

*where $\lfloor \cdot \rfloor$ denotes the round down operator, $n_s = n - \sum_{i=1}^{k} s_i$ denotes the number of separable variables, and the combination number $C_p^q$ is defined to be zero if $p < 0$, $q < 0$, or $p < q$. The probability of getting at least such a pair of nonseparable variable subsets within $l$ ($l \geq 1$) trials is*

$$P_{n1} = 1 - P_{s1}^{l} . \qquad (S\text{-}4)$$

This work was supported in part by the National Nature Science Foundation of China under Grant 61873199 and in part by China Postdoctoral Science Foundation under Grants 2014M560784 and 2016T90922. (*Corresponding Authors: Zhigang Ren and Yaochu Jin*)

Z. Ren, A. Chen, and Y. Liang are with the Autocontrol Institute, School of Electronic and Information Engineering, Xi'an Jiaotong University, Xi'an 710049, China (e-mail: renzg@mail.xjtu.edu.cn; chenan123@stu.xjtu.edu.cn; liangyongsheng@stu.xjtu.edu.cn).

Y. Jin is with the Department of Computer Science, University of Surrey, Guildford GU2 7XH, U.K. (e-mail: yaochu.jin@surrey.ac.uk).

W. Guo is with the State Key Laboratory for Manufacturing Systems Engineering, School of Mechanical Engineering, Xi'an Jiaotong University, Xi'an 710049, China (e-mail: markguo@mail.xjtu.edu.cn).

Z. Feng is with the State Key Laboratory for Manufacturing Systems Engineering, School of Electronic and Information Engineering, Xi'an Jiaotong University, Xi'an 710049, China (e-mail: fzr9910@mail.xjtu.edu.cn).

*Proof.* To randomly partition the whole variable set into two (nearly) equal-sized subsets, we just need to randomly select $\lfloor n/2 \rfloor$ variables from $n$ ones. Here we perform selection with repetition because this can greatly simplify the selection process and rarely generate two same partitions. Then the total number of possible selections is $C_n^{\lfloor n/2 \rfloor}$.

To ensure the separability between the two generated variable subsets, the only way is to keep all the variables belonging to each of the $k$ nonseparable variable groups in the same subset. This can be achieved as follows:

1) Select $\lfloor n/2 \rfloor$ variables from $n_s$ separable ones and keep all the nonseparable ones together. The number of this kind of selections is $C_{n_s}^{\lfloor n/2 \rfloor}$.

2) Select all the $s_i$ variables in the $i$th nonseparable variable group and the other $\lfloor n/2 \rfloor - s_i$ variables from $n_s$ separable ones. The number of this kind of selections is $C_{n_s}^{\lfloor n/2 \rfloor - s_i}$. Considering the $k$ nonseparable variable groups, the total number of this kind of selections is $\sum_{i=1}^{k} C_{n_s}^{\lfloor n/2 \rfloor - s_i}$.

3) Select all the $s_i$ and $s_j$ variables from the $i$th and the $j$th nonseparable variable groups, respectively, and the other $\lfloor n/2 \rfloor - s_i - s_j$ variables from $n_s$ separable ones. The number of this kind of selections is $C_{n_s}^{\lfloor n/2 \rfloor - s_i - s_j}$. Considering different combinations of two nonseparable variable groups, the total number of this kind of selections is $\sum_{i=1}^{k} \sum_{j=1, j \neq i}^{k} C_{n_s}^{\lfloor n/2 \rfloor - s_i - s_j}$.

4) Similarly, we can select variables from 3, 4, ..., or $k$ nonseparable variable groups and the remaining variables from $n_s$ separable ones. For the last case, the number of possible selections is $C_{n_s}^{\lfloor n/2 \rfloor - \sum_{i=1}^{k} s_i}$.

Then the total number of possible selections is

$$C_{n_s}^{\lfloor n/2 \rfloor} + \sum_{i=1}^{k} C_{n_s}^{\lfloor n/2 \rfloor - s_i} + \sum_{i=1}^{k} \sum_{j=1, j \neq i}^{k} C_{n_s}^{\lfloor n/2 \rfloor - s_i - s_j} + \ldots + C_{n_s}^{\lfloor n/2 \rfloor - \sum_{i=1}^{k} s_i},$$

and the probability that the two generated variable subsets can be separated from each other is $P_{s1}$ as shown in (S-3).

As we perform selection with repetition, the selections among different trials are independent. Then each pair of variable subsets obtained in $l$ ($l \geq 1$) trials can be separated with a probability of $P_{s1}^{l}$, and the probability of getting at least a pair of nonseparable variable subsets within $l$ trials is $P_{n1}$. □

**Corollary 2.** *If we further assume the $k$ nonseparable variable subsets of the function $f(\mathbf{x})$ considered in* Theorem 2 *have the same size of $s$, then any two variable subsets generated by randomly halving the original variable set are separable from each other with a probability of*

$$P_{s1}' = \frac{\sum_{i=0}^{k} C_k^i C_{n_s}^{\lfloor n/2 \rfloor - is}}{C_n^{\lfloor n/2 \rfloor}}, \tag{S-5}$$

*Proof.* This is a special case of the one considered in Theorem 2. Equation (S-5) can be obtained by simply replacing the subset sizes $s_i$ and $s_j$ in (S-3) with $s$ and counting the number of possibilities of selecting $i$ nonseparable variable groups from $k$ ones with $C_k^i$. □

**Theorem 3.** *For the function described in* Theorem 2*, its two randomly-selected variables are nonseparable from each other*

*with a probability of*

$$P_{n2} = \frac{\sum_{i=1}^{k} C_{s_i}^2}{C_n^2}, \tag{S-6}$$

*and the probability of getting at least such a pair of separable variables within $l$ ($l \geq 1$) trials is*

$$P_{s2} = 1 - P_{n2}^{l}. \tag{S-7}$$

*Proof.* The number of possibilities of selecting two variables from $n$ ones is $C_n^2$. To ensure the nonseparability between two variables, they must be selected from the same nonseparable variable group. If we select them from the $i$th group, the number of possible selections is $C_{s_i}^2$. Considering that there are $k$ nonseparable variable groups, the total number of this kind of selections becomes $\sum_{i=1}^{k} C_{s_i}^2$, and the probability of getting two nonseparable variables is $P_{n2}$ as shown in (S-6). According to the independence among different trials, we can get at least a pair of separable variables within $l$ ($l \geq 1$) trials with a probability of $P_{s2}$. □

Theorem 2 and Corollary 2 reveal that, even for an almost fully separable function, the first selection rule can easily capture a pair of nonseparable variable subsets. To verify this conclusion, let us take the function $f_4$ in the CEC'2010 benchmark suite [1] as an example. This function has a single group of 50 nonseparable variables with all the other 950 variables being separable. According to Theorem 2 or Corollary 2, it can be deduced that, when $l$ trials are allowed, the first rule is able to get at least a pair of nonseparable variable subsets for $f_4$ with a probability of

$$P_{n1} = 1 - \left( \frac{C_{950}^{500} + C_{950}^{450}}{C_{1000}^{500}} \right)^{l} = 1 - \left( 2 \times \frac{500! 950!}{450! 1000!} \right)^{l} \geq 1 - 2^{-49l}.$$

This probability is so large that the rule can almost get a pair of nonseparable variable subsets at each trial.

As a contrast, Theorem 3 reveals that the second selection rule is very likely to find a pair of separable variables even for a function containing a small number of separable variables. For instance, for the function $f_{14}$ (involving 20 nonseparable variable subsets with each including 50 elements) in the CEC'2010 benchmark suite, this rule can capture at least a pair of separable variables within $l$ trials with a probability of

$$P_{s2} = 1 - \left( \frac{20 C_{50}^2}{C_{1000}^2} \right)^{l} = 1 - 0.049^{l},$$

which approaches to 1 even when a few trials are allowed.

## S-II  Calculation of Interdependency Indicator

As indicated in Theorem 1, any feasible $\mathbf{cv}$, $\mathbf{x}_i^0$, $\mathbf{x}_i'$, $\mathbf{x}_j^0$, and $\mathbf{x}_j'$ could be employed to generate the four complete solutions required by the calculation of an indicator $\varphi(X_i, X_j)$. Nevertheless, $\mathbf{cv}$ is generally fixed to the lower bound ($\mathbf{lb}$) of $\mathbf{x}$ for different indicators for the sake of reutilization, $\mathbf{x}_i^0$ and $\mathbf{x}_i'$ are generally set to the lower bound ($\mathbf{lb}_i$) and the upper bound ($\mathbf{ub}_i$) of $\mathbf{x}_i$, respectively, and similar settings are provided for $\mathbf{x}_j^0$ and $\mathbf{x}_j'$.

Fig. S-1 graphically presents the relations among the four

solutions and the intermediate items concerned in the calculation of $\varphi(X_i, X_j)$. For the convenience of the description, the solution $\mathbf{cv} \leftarrow |\mathbf{x}_i^0 = \mathbf{lb}_i, \mathbf{x}_j^0 = \mathbf{lb}_j$ is simply denoted as $\mathbf{x}^{l,l}$. Similar notations, including $\mathbf{x}^{u,l}$, $\mathbf{x}^{l,u}$, and $\mathbf{x}^{u,u}$, are given for the other three solutions.

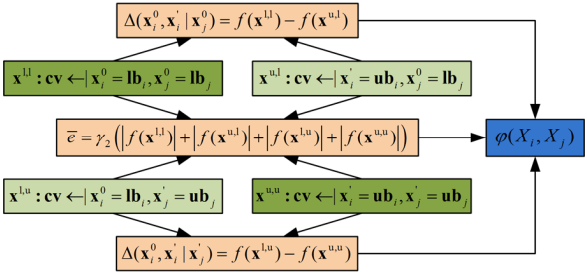

Fig. S-1. The relations among the four solutions and the intermediate items required by the calculation of an interdependency indicator

## S-III Computational Complexity of FDG on Partially Separable Instances

When decomposing an $n$-dimensional partially separable instance which has no separable variable but involves $k$ groups of nonseparable variables with each group containing $s$ elements, the BTDP in FDG constructs a binary tree for each variable group. In the worst case, it uniformly assigns the $(s-1)$ nonseparable partners of a selected variable to all the nodes in each layer of the binary tree because this tends not to exclude any node and results in the most number of detections. To make the analysis convenient, we assume each variable has $n$ partners (or nonseparable partners) instead of $(n-1)$ (or $(s-1)$) ones and both $n$ and $s$ are integral powers of 2. Then for the first selected variable, BTDP will construct a binary tree that can be schematically described by Fig. S-2. It includes $(1 + \log_2 s + \log_2 k)$ layers. The $i$th ($i = 0, 1, \ldots, \log_2 s$) layer has $2^i$ nodes and each node involves $n/2^i$ variables, $s/2^i$ of which are nonseparable from the selected one. As a consequence, no variable is excluded in the first $(1 + \log_2 s)$ layers. Each of the last $\log_2 k$ layers in the binary tree has $2s$ nodes. Half of them are excluded due to their separability from the selected variable, and each of the other $s$ nodes includes a variable nonseparable from the selected one and is further halved until just containing a single variable.

Considering that the root node requires 3 new FEs, each left child node requires 2 FEs, and each right child node does not require any FE, it can be deduced that BTDP needs

$$1 \times 3 + \sum_{i=1}^{\log_2 s} 2^i \Big/ 2 \times 2 + 2s \times \log_2(k)/2 \times 2 = 2s\log_2(k) + 2s + 1$$

FEs to get the first group of nonseparable variables. Then the total number of FEs required by the $k$ groups of variables is

$$2s \times \sum_{i=1}^{k} \log_2(i) + (2s+1) \times k ,$$
$$< 2ks\log_2(k) + 2ks + k = 2n\log_2(k) + 2n + k$$

and the total computational complexity of FDG on this type of instances is $O(5l + 2 + 2l + 2n\log_2(k) + 2n + k) = O(n\log_2 k)$ in the worst case.

In the best case, BTDP concentratively assigns the $(s-1)$ nonseparable partners of a selected variable to one of the nodes in each layer such that those variables separable from the selected one can be quickly excluded and the number of detections can be reduced to the greatest extent possible. Fig. S-3 presents the schematic binary tree for the first selected variable. It also includes $(1 + \log_2 k + \log_2 s)$ layers. Except the root node which corresponds to the 0th layer, each of the first $(1 + \log_2 k)$ layers has 2 nodes, one of which is separable from the selected one and thus is directly excluded. Each of the last $\log_2 s$ layers has $2^i$ ($i = 1, \ldots, \log_2 s$) nodes. As these nodes are nonseparable from the selected variable, they are persistently halved until only a single variable is involved in each of them. With the similar analysis method shown above, we can know that BTDP needs

$$1 \times 3 + \log_2(k) \times 2 + \sum_{i=1}^{\log_2 s} 2^i \Big/ 2 \times 2 = 2\log_2(k) + 2s + 1 ,$$

FEs to construct the binary tree for the first group of nonseparable variables. Then the total number of FEs required by the $k$ groups of variables is

$$2 \times \sum_{i=1}^{k} \log_2(i) + k \times (2s+1) < 2k\log_2(k) + 2n + k ,$$

and the computational complexity of FDG on this type of instances becomes $O(5l + 2 + 2l + 2k\log_2(k) + 2n + k) = O(\max(n, k\log_2 k))$ in the best case.

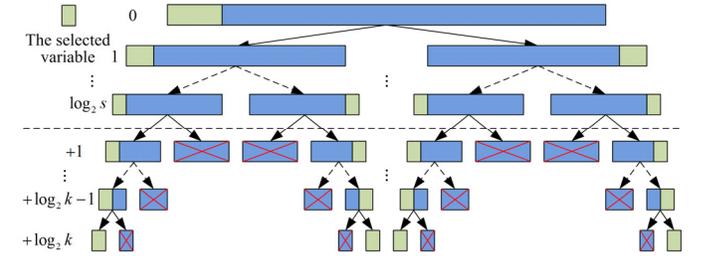

Fig. S-2. The schematic binary tree constructed by BTDP for the first selected variable in a $k \times s$ partially separable instance in the worst case.

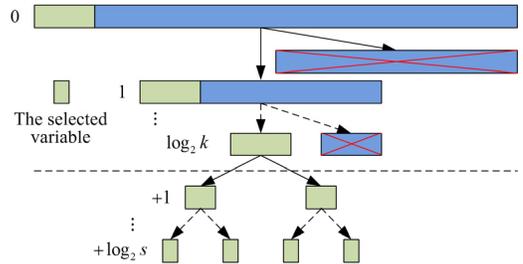

Fig. S-3. The schematic binary tree constructed by BTDP for the first selected variable in a $k \times s$ partially separable instance in the best case.

## S-IV  Detailed Experimental Results

### TABLE S-I

Decomposition accuracy and efficiency of DG, GDG, FII, DG2, RDG, and FDG on each benchmark function in the CEC'2010 and the CEC'2013 suites. $\rho_1$, $\rho_2$, and $\rho_3$ denote the NMI indicators on all the variables, nonseparable ones, and separable ones, respectively. If a function does not contain nonseparable or separable variables, the corresponding entry in the table is marked with "–". The result in boldface indicates that it is of the highest accuracy and is obtained with the fewest FEs.

| Typ. | Fun. | DG | | | | GDG | | | FII | | | | DG2 | | | RDG | | | | FDG | | | |
| --- | --- | --- | --- | --- | --- | --- | --- | --- | --- | --- | --- | --- | --- | --- | --- | --- | --- | --- | --- | --- | --- | --- | --- |
| | | $\rho_1$ | $\rho_2$ | $\rho_3$ | feNum | $\rho_1$ | $\rho_2$ | $\rho_3$ | $\rho_1$ | $\rho_2$ | $\rho_3$ | feNum | $\rho_1$ | $\rho_2$ | $\rho_3$ | $\rho_1$ | $\rho_2$ | $\rho_3$ | feNum | $\rho_1$ | $\rho_2$ | $\rho_3$ | feNum |
| T1 | CEC'2010 $f_1$ | 100 | – | 100 | 1.00e+06 | 100 | – | 100 | 100 | – | 100 | 3.00e+03 | 100 | – | 100 | 100 | – | 100 | 3.00e+03 | **100** | – | **100** | **5.20e+01** |
| | CEC'2010 $f_2$ | 100 | – | 100 | 1.00e+06 | 100 | – | 100 | 100 | – | 100 | 3.00e+03 | 100 | – | 100 | 100 | – | 100 | 3.00e+03 | **100** | – | **100** | **5.20e+01** |
| | CEC'2010 $f_3$ | 100 | – | 100 | 1.00e+06 | 0.00 | – | 0.00 | 100 | – | 100 | 3.00e+03 | 0.00 | – | 0.00 | 0.00 | – | 0.00 | 6.00e+03 | 0.00 | – | 0.00 | 5.20e+01 |
| | CEC'2013 $f_1$ | 100 | – | 100 | 1.00e+06 | 100 | – | 100 | 100 | – | 100 | 3.00e+03 | 100 | – | 100 | 100 | – | 100 | 3.00e+03 | **100** | – | **100** | **5.20e+01** |
| | CEC'2013 $f_2$ | 100 | – | 100 | 1.00e+06 | 100 | – | 100 | 100 | – | 100 | 3.00e+03 | 100 | – | 100 | 100 | – | 100 | 3.00e+03 | **100** | – | **100** | **5.20e+01** |
| | CEC'2013 $f_3$ | 100 | – | 100 | 1.00e+06 | 0.00 | – | 0.00 | 100 | – | 100 | 3.00e+03 | 0.00 | – | 0.00 | 0.00 | – | 0.00 | 6.00e+03 | 0.00 | – | 0.00 | 5.20e+01 |
| T2-1 | CEC'2010 $f_4$ | 48.75 | 100 | 46.35 | 1.45e+04 | 100 | 100 | 100 | 75.49 | 100 | 74.55 | 3.73e+03 | 100 | 100 | 100 | 100 | 100 | 100 | 4.20e+03 | **100** | **100** | **100** | **2.15e+03** |
| | CEC'2010 $f_5$ | 100 | 100 | 100 | 9.05e+05 | 100 | 100 | 100 | 100 | 100 | 100 | 3.05e+03 | 100 | 100 | 100 | 100 | 100 | 100 | 4.15e+03 | **100** | **100** | **100** | **2.15e+03** |
| | CEC'2010 $f_6$ | 100 | 100 | 100 | 9.06e+04 | 100 | 100 | 100 | 100 | 100 | 100 | 3.05e+03 | 22.22 | 100 | 17.91 | 100 | 100 | 100 | 5.00e+04 | 13.85 | 100 | 8.83 | 2.44e+03 |
| | CEC'2010 $f_7$ | 58.30 | 100 | 58.11 | 6.77e+04 | 100 | 100 | 100 | 100 | 100 | 100 | 3.05e+03 | 100 | 100 | 100 | 100 | 100 | 100 | 4.23e+03 | **100** | **100** | **100** | **2.15e+03** |
| | CEC'2010 $f_8$ | 47.98 | 100 | 45.64 | 2.32e+04 | 100 | 100 | 100 | 77.41 | 100 | 76.56 | 3.58e+03 | 100 | 100 | 100 | 100 | 100 | 100 | 5.60e+03 | **100** | **100** | **100** | **2.46e+03** |
| | CEC'2010 $f_9$ | 100 | 100 | 100 | 2.70e+05 | 100 | 100 | 100 | 100 | 100 | 100 | 8.01e+03 | 100 | 100 | 100 | 100 | 100 | 100 | 1.40e+04 | **100** | **100** | **100** | **4.74e+03** |
| | CEC'2010 $f_{10}$ | 100 | 100 | 100 | 2.72e+05 | 100 | 100 | 100 | 100 | 100 | 100 | 8.01e+03 | 100 | 100 | 100 | 100 | 100 | 100 | 1.40e+04 | **100** | **100** | **100** | **4.75e+03** |
| | CEC'2010 $f_{11}$ | 99.95 | 99.79 | 100 | 2.70e+05 | 54.28 | 0.00 | 0.00 | 99.31 | 97.11 | 100 | 9.49e+03 | 54.28 | 0.00 | 0.00 | 54.28 | 0.00 | 0.00 | 1.36e+04 | 54.28 | 100 | 0.00 | 4.65e+03 |
| | CEC'2010 $f_{12}$ | 100 | 100 | 100 | 2.71e+05 | 100 | 100 | 100 | 100 | 100 | 100 | 8.01e+03 | 100 | 100 | 100 | 100 | 100 | 100 | 1.43e+04 | **100** | **100** | **100** | **4.80e+03** |
| | CEC'2010 $f_{13}$ | 63.04 | 51.49 | 58.25 | 5.03e+04 | 100 | 100 | 100 | 100 | 100 | 100 | 9.61e+04 | 100 | 100 | 100 | 100 | 100 | 100 | 2.92e+04 | **100** | **100** | **100** | **9.75e+03** |
| | CEC'2013 $f_4$ | 53.58 | 93.17 | 36.47 | 1.56e+04 | 100 | 100 | 100 | 91.79 | 0.00 | 97.88 | 3.69e+03 | 100 | 100 | 100 | 100 | 100 | 100 | 9.84e+03 | **100** | **100** | **100** | **3.31e+03** |
| | CEC'2013 $f_5$ | 98.67 | 87.45 | 100 | 5.27e+05 | 96.69 | 73.25 | 100 | 96.59 | 72.66 | 100 | 4.64e+03 | 100 | 100 | 100 | 100 | 100 | 100 | 1.01e+04 | **100** | **100** | **100** | **3.31e+03** |
| | CEC'2013 $f_6$ | 98.53 | 86.29 | 100 | 5.79e+05 | 98.62 | 86.98 | 100 | 95.82 | 68.24 | 100 | 3.58e+03 | 33.37 | 100 | 0.00 | 37.81 | 100 | 7.89 | 1.32e+04 | 33.37 | 100 | 0.00 | 3.73e+03 |
| | CEC'2013 $f_7$ | 37.83 | 32.52 | 33.62 | 1.14e+04 | 100 | 100 | 100 | 89.42 | 70.86 | 89.72 | 6.23e+03 | 99.16 | 90.25 | 100 | 97.84 | 70.86 | 100 | 9.02e+03 | **100** | **100** | **100** | **3.39e+03** |
| T2-2 | CEC'2010 $f_{14}$ | 100 | 100 | – | 2.10e+04 | 100 | 100 | – | 100 | 100 | – | 2.30e+04 | 100 | 100 | – | 100 | 100 | – | 2.05e+04 | **100** | **100** | – | **6.96e+03** |
| | CEC'2010 $f_{15}$ | 100 | 100 | – | 2.10e+04 | 100 | 100 | – | 100 | 100 | – | 2.30e+04 | 100 | 100 | – | 100 | 100 | – | 2.05e+04 | **100** | **100** | – | **6.94e+03** |
| | CEC'2010 $f_{16}$ | 99.67 | 99.57 | – | 2.11e+04 | 100 | 100 | – | 96.32 | 96.32 | – | 2.88e+04 | 100 | 100 | – | 100 | 100 | – | 2.09e+04 | **100** | **100** | – | **7.08e+03** |
| | CEC'2010 $f_{17}$ | 100 | 100 | – | 2.10e+04 | 100 | 100 | – | 100 | 100 | – | 2.30e+04 | 100 | 100 | – | 100 | 100 | – | 2.07e+04 | **100** | **100** | – | **7.03e+03** |
| | CEC'2010 $f_{18}$ | 45.39 | 45.39 | – | 3.96e+04 | 100 | 100 | – | 100 | 100 | – | 3.70e+05 | 100 | 100 | – | 100 | 100 | – | 4.98e+04 | **100** | **100** | – | **1.74e+04** |
| | CEC'2013 $f_8$ | 94.03 | 94.03 | – | 2.26e+04 | 74.72 | 74.72 | – | 54.38 | 54.38 | – | 1.31e+04 | 86.29 | 86.29 | – | 84.83 | 84.83 | – | 1.99e+04 | 81.32 | 81.32 | – | 8.40e+03 |
| | CEC'2013 $f_9$ | 100 | 100 | – | 1.76e+04 | 98.51 | 98.51 | – | 100 | 100 | – | 2.06e+04 | 100 | 100 | – | 100 | 100 | – | 1.92e+04 | **100** | **100** | – | **6.49e+03** |
| | CEC'2013 $f_{10}$ | 89.81 | 89.81 | – | 4.86e+04 | 90.08 | 90.08 | – | 87.49 | 87.49 | – | 1.70e+04 | 100 | 100 | – | 91.53 | 91.53 | – | 1.91e+04 | **100** | **100** | – | **6.53e+03** |
| | CEC'2013 $f_{11}$ | 54.14 | 54.14 | – | 9.10e+03 | 100 | 100 | – | 100 | 100 | – | 5.00e+03 | 100 | 100 | – | 100 | 100 | – | 1.04e+04 | **100** | **100** | – | **6.60e+03** |
| T3-1 | CEC'2010 $f_{20}$ | 0.00 | 0.00 | – | 1.55e+05 | 100 | 100 | – | 100 | 100 | – | 5.03e+05 | 100 | 100 | – | 100 | 100 | – | 5.08e+04 | **100** | **100** | – | **1.80e+04** |
| | CEC'2013 $f_{12}$ | 0.00 | 0.00 | – | 1.49e+05 | 100 | 100 | – | 100 | 100 | – | 5.03e+05 | 100 | 100 | – | 100 | 100 | – | 5.08e+04 | **100** | **100** | – | **1.80e+04** |
| | CEC'2013 $f_{13}$ | 0.00 | 0.00 | – | 5.86e+03 | 100 | 100 | – | **100** | **100** | – | **4.59e+03** | 100 | 100 | – | 0.00 | 0.00 | – | 8.31e+03 | 100 | 100 | – | 5.35e+03 |
| | CEC'2013 $f_{54}$ | 0.00 | 0.00 | – | 1.39e+04 | 100 | 100 | – | 100 | 100 | – | 6.08e+03 | 100 | 100 | – | 100 | 100 | – | 9.54e+03 | **100** | **100** | – | **5.36e+03** |
| T3-2 | CEC'2010 $f_{19}$ | 100 | 100 | – | 2.00e+03 | 100 | 100 | – | 100 | 100 | – | 4.00e+03 | 100 | 100 | – | 100 | 100 | – | 6.00e+03 | **100** | **100** | – | **5.20e+01** |
| | CEC'2013 $f_{15}$ | 100 | 100 | – | 2.00e+03 | 100 | 100 | – | 100 | 100 | – | 4.00e+03 | 100 | 100 | – | 100 | 100 | – | 6.16e+03 | **100** | **100** | – | **5.20e+01** |

TABLE S-II

The median, the mean, and the standard deviation of the best solutions obtained by FDG on each benchmark function in the CEC'2010 and the CEC'2013 suites over 25 independent runs with each run allowing a maximum number of $3.0\times10^6$ FEs. FDG is embedded into the DECC framework and is compared with DG, GDG, FII, DG2, RDG, and ID. The superscripts "−", "≈", and "+" indicate that the corresponding result is worse than, similar to, and better than that of FDG, respectively.

**CEC'2010 suite**

| DG | GDG | FII | DG2 | RDG | ID | FDG | Fun. |
|---|---|---|---|---|---|---|---|
| 9.51e+05$^-$ | 2.57e+05$^-$ | 7.28e+04$^+$ | 2.57e+05$^-$ | 7.28e+04$^+$ | 7.27e+04$^+$ | 7.27e+04 | $f_1$ |
| 3.08e+06 | 5.36e+05 | 2.44e+05 | 5.35e+05 | 2.44e+05 | 2.44e+05 | 2.44e+05 | |
| 6.06e+06 | 6.55e+05 | 5.41e+05 | 6.55e+05 | 5.41e+05 | 5.41e+05 | 5.41e+05 | |
| 4.28e+03$^+$ | 4.28e+03$^+$ | 4.28e+03$^-$ | 4.28e+03$^+$ | 4.28e+03$^+$ | 4.28e+03$^+$ | 4.28e+03 | $f_2$ |
| 4.30e+03 | 4.30e+03 | 4.30e+03 | 4.30e+03 | 4.30e+03 | 4.30e+03 | 4.30e+03 | |
| 4.10e+02 | 4.08e+02 | 4.09e+02 | 4.08e+02 | 4.09e+02 | 4.09e+02 | 4.09e+02 | |
| 1.13e+01$^+$ | 1.13e+01$^+$ | 1.13e+01$^+$ | 1.13e+01$^+$ | 1.13e+01$^-$ | 1.13e+01$^+$ | 1.13e+01 | $f_3$ |
| 1.11e+01 | 1.11e+01 | 1.11e+01 | 1.11e+01 | 1.11e+01 | 1.11e+01 | 1.11e+01 | |
| 1.10e+00 | 1.10e+00 | 1.10e+00 | 1.10e+00 | 1.10e+00 | 1.10e+00 | 1.10e+00 | |
| 5.37e+11$^-$ | 3.75e+10$^+$ | 5.61e+10$^-$ | 3.75e+10$^+$ | 3.09e+10$^+$ | 3.09e+10$^+$ | 3.09e+10 | $f_4$ |
| 6.02e+11 | 4.33e+10 | 7.03e+10 | 4.32e+10 | 3.26e+10 | 3.26e+10 | 3.26e+10 | |
| 3.94e+11 | 2.97e+10 | 6.51e+10 | 2.97e+10 | 2.55e+10 | 2.55e+10 | 2.55e+10 | |
| 7.25e+07$^+$ | 7.05e+07$^+$ | 6.97e+07$^+$ | 7.05e+07$^+$ | 6.97e+07$^+$ | 6.97e+07$^+$ | 6.59e+07 | $f_5$ |
| 6.89e+07 | 6.76e+07 | 6.76e+07 | 6.76e+07 | 6.59e+07 | 6.59e+07 | 6.59e+07 | |
| 1.25e+07 | 1.23e+07 | 1.24e+07 | 1.23e+07 | 1.24e+07 | 1.24e+07 | 1.24e+07 | |
| 1.57e+01$^-$ | 1.57e+01$^-$ | 1.57e+01$^+$ | 1.03e+01$^-$ | 1.57e+01$^-$ | 1.57e+01$^-$ | 1.19e+01 | $f_6$ |
| 1.58e+01 | 1.57e+01 | 1.57e+01 | 1.01e+01 | 1.57e+01 | 1.57e+01 | 1.18e+01 | |
| 1.31e+00 | 1.32e+00 | 1.32e+00 | 1.56e+00 | 1.32e+00 | 1.32e+00 | 1.35e+00 | |
| 2.13e+04$^-$ | 1.94e+04 | 1.06e+04$^+$ | 1.94e+04$^-$ | 1.06e+04$^+$ | 1.06e+04$^+$ | 1.06e+04 | $f_7$ |
| 2.61e+04 | 1.95e+04 | 1.95e+04 | 1.95e+04 | 1.15e+04 | 1.12e+04 | 1.15e+04 | |
| 2.22e+04 | 8.10e+03 | 5.90e+03 | 8.10e+03 | 5.90e+03 | 5.57e+03 | 5.89e+03 | |
| 3.15e+07$^-$ | 1.85e+04 | 3.38e+01$^-$ | 1.85e+04$^-$ | 7.90e+03$^+$ | 7.86e+03$^+$ | 7.87e+03 | $f_8$ |
| 3.40e+07 | 3.34e+05 | 6.38e+05 | 3.34e+05 | 3.28e+05 | 3.28e+05 | 3.28e+05 | |
| 3.02e+07 | 1.10e+06 | 1.49e+06 | 1.10e+06 | 1.10e+06 | 1.10e+06 | 1.10e+06 | |
| 3.76e+07$^-$ | 4.59e+07$^-$ | 2.89e+07$^+$ | 4.58e+07$^-$ | 2.89e+07$^+$ | 2.86e+07$^+$ | 2.89e+07 | $f_9$ |
| 4.92e+07 | 6.47e+07 | 3.72e+07 | 6.41e+07 | 3.72e+07 | 3.59e+07 | 3.70e+07 | |
| 3.75e+07 | 4.54e+07 | 2.79e+07 | 4.55e+07 | 2.79e+07 | 2.78e+07 | 2.79e+07 | |
| 3.40e+03$^-$ | 3.63e+03$^-$ | 3.21e+03$^+$ | 3.63e+03$^-$ | 3.21e+03$^+$ | 3.20e+03$^-$ | 3.20e+03 | $f_{10}$ |
| 3.39e+03 | 3.61e+03 | 3.20e+03 | 3.60e+03 | 3.21e+03 | 3.17e+03 | 3.18e+03 | |
| 1.58e+02 | 1.73e+02 | 1.36e+02 | 1.75e+02 | 1.37e+02 | 1.27e+02 | 1.36e+02 | |
| 2.57e+01$^+$ | 2.46e+01$^+$ | 3.05e+01$^-$ | 2.46e+01$^+$ | 2.71e+01$^-$ | 2.50e+01$^-$ | 2.51e+01 | $f_{11}$ |
| 2.60e+01 | 2.53e+01 | 3.01e+01 | 2.53e+01 | 2.63e+01 | 2.50e+01 | 2.53e+01 | |
| 2.48e+00 | 3.36e+00 | 2.59e+00 | 2.99e+00 | 2.24e+00 | 2.47e+00 | 2.25e+00 | |
| 3.21e+04$^+$ | 3.83e+04$^-$ | 2.76e+04$^+$ | 3.78e+04$^-$ | 2.76e+04$^+$ | 2.90e+04$^+$ | 2.72e+04 | $f_{12}$ |
| 3.11e+04 | 3.74e+04 | 2.66e+04 | 3.71e+04 | 2.66e+04 | 2.56e+04 | 2.66e+04 | |
| 1.12e+04 | 1.27e+04 | 1.08e+04 | 1.25e+04 | 1.08e+04 | 1.07e+04 | 1.08e+04 | |
| 1.44e+04$^+$ | 1.81e+04$^-$ | 1.38e+04$^+$ | 1.81e+04$^-$ | 1.34e+04$^+$ | 1.26e+04$^+$ | 1.33e+04 | $f_{13}$ |
| 1.52e+04 | 1.95e+04 | 1.38e+04 | 1.94e+04 | 1.33e+04 | 1.28e+04 | 1.33e+04 | |
| 7.68e+03 | 6.43e+03 | 5.26e+03 | 6.42e+03 | 5.24e+03 | 5.24e+03 | 5.24e+03 | |
| 2.16e+07$^+$ | 2.76e+07$^-$ | 2.16e+07$^+$ | 2.76e+07$^-$ | 2.16e+07$^+$ | 2.14e+07$^+$ | 2.15e+07 | $f_{14}$ |
| 2.16e+07 | 2.79e+07 | 2.16e+07 | 2.78e+07 | 2.16e+07 | 2.14e+07 | 2.15e+07 | |
| 1.66e+06 | 2.12e+06 | 1.66e+06 | 2.11e+06 | 1.66e+06 | 1.63e+06 | 1.64e+06 | |
| 2.67e+03$^+$ | 2.91e+03$^-$ | 2.67e+03$^+$ | 2.91e+03$^-$ | 2.67e+03$^+$ | 2.65e+03$^+$ | 2.65e+03 | $f_{15}$ |
| 2.73e+03 | 3.01e+03 | 2.73e+03 | 3.01e+03 | 2.73e+03 | 2.72e+03 | 2.72e+03 | |
| 2.47e+02 | 2.55e+02 | 2.47e+02 | 2.55e+02 | 2.47e+02 | 2.47e+02 | 2.44e+02 | |
| 2.09e+01$^+$ | 1.91e+01$^-$ | 2.65e+01$^-$ | 1.91e+01$^+$ | 1.91e+01$^+$ | 1.91e+01$^+$ | 1.91e+01 | $f_{16}$ |
| 1.96e+01 | 1.84e+01 | 2.59e+01 | 1.84e+01 | 1.84e+01 | 1.84e+01 | 1.84e+01 | |
| 3.75e+00 | 3.08e+00 | 3.81e+00 | 3.08e+00 | 3.08e+00 | 3.08e+00 | 3.08e+00 | |
| 8.20e+00$^+$ | 4.99e+01$^-$ | 8.22e+00$^+$ | 4.97e+01$^-$ | 8.20e+00$^+$ | 7.67e+00$^+$ | 7.93e+00 | $f_{17}$ |
| 8.81e+00 | 5.26e+01 | 8.86e+00 | 5.23e+01 | 8.81e+00 | 8.04e+00 | 8.30e+00 | |
| 3.71e+00 | 1.48e+01 | 3.72e+00 | 1.48e+01 | 3.71e+00 | 2.95e+00 | 3.00e+00 | |
| 2.38e+10$^-$ | 1.39e+03$^-$ | 1.32e+03$^-$ | 1.39e+03$^-$ | 1.21e+03$^+$ | 1.18e+03$^+$ | 1.19e+03 | $f_{18}$ |
| 2.42e+10 | 1.38e+03 | 1.32e+03 | 1.38e+03 | 1.18e+03 | 1.16e+03 | 1.16e+03 | |
| 4.87e+09 | 1.61e+02 | 1.55e+02 | 1.61e+02 | 1.67e+02 | 1.62e+02 | 1.62e+02 | |
| 4.45e+11$^-$ | 1.86e+07$^-$ | 1.86e+07$^-$ | 1.86e+07$^-$ | 5.09e+06$^+$ | 5.01e+06$^+$ | 5.00e+06 | $f_{20}$ |
| 4.45e+10 | 1.42e+08 | 1.42e+08 | 1.42e+08 | 7.48e+07 | 7.42e+07 | 7.43e+07 | |
| 5.93e+09 | 3.19e+08 | 3.19e+08 | 3.19e+08 | 2.29e+08 | 2.29e+08 | 2.29e+08 | |
| 9.01e+05$^-$ | 1.01e+06$^-$ | 9.01e+05$^-$ | 1.01e+06$^-$ | 9.01e+05$^-$ | 9.01e+05$^-$ | 9.01e+05 | $f_{19}$ |
| 9.04e+05 | 1.02e+06 | 9.04e+05 | 1.02e+06 | 9.05e+05 | 9.03e+05 | 9.04e+05 | |
| 4.99e+04 | 5.14e+04 | 4.99e+04 | 5.14e+04 | 4.99e+04 | 4.99e+04 | 4.99e+04 | |

**CEC'2013 suite**

| Fun. | Type | FDG | ID | RDG | DG2 | FII | GDG | DG |
|---|---|---|---|---|---|---|---|---|
| $f_1$ | $T_1$ | 1.41e+05 | 1.41e+05$^+$ | 1.41e+05$^+$ | 6.17e+05$^-$ | 1.41e+05$^+$ | 6.17e+05$^-$ | 1.56e+06$^-$ |
| | | 1.02e+06 | 1.02e+06 | 1.02e+06 | 1.53e+06 | 1.02e+06 | 1.53e+06 | 2.96e+06 |
| | | 2.25e+06 | 2.25e+06 | 2.25e+06 | 2.30e+06 | 2.25e+06 | 2.30e+06 | 3.42e+06 |
| $f_2$ | | 1.48e+04 | 1.48e+04$^+$ | 1.48e+04$^+$ | 1.48e+04$^+$ | 1.48e+04$^+$ | 1.48e+04$^+$ | 1.48e+04$^+$ |
| | | 1.44e+04 | 1.44e+04 | 1.44e+04 | 1.44e+04 | 1.44e+04 | 1.44e+04 | 1.44e+04 |
| | | 1.65e+03 | 1.65e+03 | 1.65e+03 | 1.65e+03 | 1.65e+03 | 1.65e+03 | 1.65e+03 |
| $f_3$ | | 2.06e+01 | 2.06e+01$^+$ | 2.06e+01$^+$ | 2.06e+01$^+$ | 2.06e+01$^+$ | 2.06e+01 | 2.07e+01$^-$ |
| | | 2.06e+01 | 2.06e+01 | 2.06e+01 | 2.06e+01 | 2.06e+01 | 2.06e+01 | 2.07e+01 |
| | | 9.07e-03 | 9.97e-03 | 8.84e-03 | 9.14e-03 | 9.14e-03 | 9.36e-03 | 1.12e-02 |
| $f_4$ | $T_{2\text{-}1}$ | 2.55e+08 | 2.55e+08$^-$ | 2.55e+08$^-$ | 3.55e+08$^-$ | 4.05e+09$^-$ | 3.55e+08$^-$ | 1.08e+11$^-$ |
| | | 2.97e+08 | 2.97e+08 | 2.97e+08 | 4.23e+08 | 4.17e+09 | 4.23e+08 | 1.16e+11 |
| | | 1.51e+08 | 1.51e+08 | 1.51e+08 | 1.90e+08 | 1.44e+09 | 1.90e+08 | 5.94e+10 |
| $f_5$ | | 2.26e+06 | 2.26e+06$^+$ | 2.26e+06$^-$ | 2.26e+06$^-$ | 2.46e+06$^-$ | 2.40e+06$^-$ | 2.77e+06$^-$ |
| | | 2.27e+06 | 2.27e+06 | 2.27e+06 | 2.35e+06 | 2.41e+06 | 2.59e+06 | 2.76e+06 |
| | | 3.06e+05 | 3.06e+05 | 3.06e+05 | 3.68e+05 | 4.66e+05 | 7.82e+05 | 8.64e+05 |
| $f_6$ | | 1.06e+06 | 1.06e+06$^+$ | 1.06e+06$^-$ | 1.06e+06$^-$ | 1.06e+06$^-$ | 1.06e+06$^-$ | 1.06e+06$^-$ |
| | | 1.06e+06 | 1.06e+06 | 1.06e+06 | 1.06e+06 | 1.06e+06 | 1.06e+06 | 1.06e+06 |
| | | 1.51e+03 | 1.93e+03 | 1.53e+03 | 1.91e+03 | 1.97e+03 | 1.75e+03 | 1.24e+03 |
| $f_7$ | | 7.09e+04 | 7.09e+04$^+$ | 1.23e+07$^-$ | 1.46e+06$^-$ | 5.33e+06$^-$ | 1.06e+05$^-$ | 1.03e+08$^-$ |
| | | 6.10e+05 | 6.10e+05 | 1.38e+07 | 1.91e+06 | 6.48e+06 | 4.55e+06 | 1.18e+08 |
| | | 2.70e+06 | 2.70e+06 | 9.06e+06 | 1.43e+06 | 4.03e+06 | 2.21e+07 | 4.59e+07 |
| $f_8$ | $T_{2\text{-}2}$ | 3.58e+13 | 8.21e+13$^-$ | 6.62e+13$^-$ | 9.66e+13$^-$ | 4.32e+13$^-$ | 5.20e+13$^-$ | 1.50e+15$^-$ |
| | | 4.01e+13 | 8.67e+13 | 7.45e+13 | 1.05e+14 | 5.24e+13 | 5.76e+13 | 2.35e+15 |
| | | 4.14e+13 | 4.73e+13 | 4.32e+13 | 7.75e+13 | 2.59e+13 | 3.67e+13 | 2.17e+15 |
| $f_9$ | | 2.88e+08 | 2.88e+08$^-$ | 2.88e+08$^-$ | 3.25e+08$^-$ | 2.88e+08$^-$ | 3.15e+08$^-$ | 2.88e+08$^-$ |
| | | 2.67e+08 | 2.67e+08 | 2.67e+08 | 2.85e+08 | 2.67e+08 | 2.99e+08 | 2.67e+08 |
| | | 7.54e+07 | 7.54e+07 | 7.54e+07 | 8.02e+07 | 7.54e+07 | 8.80e+07 | 7.54e+07 |
| $f_{10}$ | | 9.43e+07 | 9.43e+07$^-$ | 9.43e+07$^-$ | 9.45e+07$^-$ | 9.42e+07$^-$ | 9.41e+07$^-$ | 9.44e+07$^-$ |
| | | 9.43e+07 | 9.43e+07 | 9.43e+07 | 9.44e+07 | 9.42e+07 | 9.41e+07 | 9.44e+07 |
| | | 2.71e+06 | 2.71e+06 | 3.81e+06 | 2.99e+06 | 2.63e+06 | 3.83e+06 | 2.83e+06 |
| $f_{11}$ | | 1.77e+08 | 1.77e+08$^+$ | 4.29e+08$^-$ | 2.35e+08$^-$ | 3.63e+08$^-$ | 5.28e+08$^-$ | 5.05e+09$^-$ |
| | | 1.58e+10 | 1.58e+10 | 7.09e+08 | 1.78e+10 | 4.00e+08 | 1.08e+09 | 2.78e+10 |
| | | 3.60e+10 | 3.60e+10 | 8.51e+08 | 4.03e+10 | 1.36e+08 | 2.04e+09 | 7.06e+10 |
| $f_{12}$ | $T_{3\text{-}1}$ | 4.15e+06 | 3.76e+06$^+$ | 6.74e+06$^-$ | 2.74e+07$^-$ | 2.74e+07$^-$ | 2.74e+07$^-$ | 1.48e+11$^-$ |
| | | 4.44e+07 | 4.41e+07 | 4.57e+07 | 1.46e+08 | 1.46e+08 | 1.46e+08 | 1.46e+11 |
| | | 1.83e+08 | 1.83e+08 | 1.82e+08 | 3.19e+08 | 3.19e+08 | 3.19e+08 | 1.60e+10 |
| $f_{13}$ | | 4.50e+08 | 4.47e+08$^+$ | 2.50e+09$^-$ | 7.60e+08$^-$ | 4.50e+08$^-$ | 7.60e+08$^-$ | 7.09e+09$^-$ |
| | | 5.29e+08 | 5.27e+08 | 2.65e+09 | 7.50e+08 | 5.29e+08 | 7.51e+08 | 7.11e+09 |
| | | 2.18e+08 | 2.18e+08 | 8.22e+08 | 2.66e+08 | 2.18e+08 | 2.67e+08 | 2.16e+09 |
| $f_{14}$ | | 4.03e+08 | 4.00e+08$^+$ | 4.07e+08$^-$ | 9.81e+08$^-$ | 4.04e+08$^-$ | 9.82e+08$^-$ | 9.34e+09$^-$ |
| | | 6.32e+08 | 6.26e+08 | 6.38e+08 | 1.48e+09 | 6.33e+08 | 1.48e+09 | 1.10e+10 |
| | | 7.70e+08 | 7.68e+08 | 7.71e+08 | 1.40e+09 | 7.70e+08 | 1.40e+09 | 8.43e+09 |
| $f_{15}$ | $T_{3\text{-}2}$ | 4.62e+06 | 4.62e+06$^-$ | 4.63e+06$^-$ | 5.34e+06$^-$ | 4.63e+06$^-$ | 5.34e+06$^-$ | 4.63e+06$^-$ |
| | | 5.27e+06 | 5.27e+06 | 5.28e+06 | 6.05e+06 | 5.28e+06 | 6.05e+06 | 5.28e+06 |
| | | 1.70e+06 | 1.70e+06 | 1.70e+06 | 1.76e+06 | 1.70e+06 | 1.76e+06 | 1.70e+06 |

TABLE S-III

The median, the mean, and the standard deviation of the best solutions obtained by FDG on each benchmark function in the CEC'2010 and the CEC'2013 suites over 25 independent runs with each run allowing a maximum number of $6.0\times10^5$ FEs. FDG is embedded into the DECC framework and is compared with DG, GDG, FII, DG2, RDG, and ID. DG does not support the optimization on some functions because it spends all the FEs in the decomposition, and the corresponding entries in the table are marked with "–". The superscripts "$-$," "$\approx$," and "$+$" indicate that the corresponding result is worse than, similar to, and better than that of FDG, respectively.

**CEC'2010 suite**

| DG | GDG | FII | DG2 | RDG | ID | FDG | Fun. |
|---|---|---|---|---|---|---|---|
| –$^-$ | 1.82e+09$^-$ | 5.71e+07$^+$ | 1.78e+09$^-$ | 5.71e+07$^+$ | 5.47e+07$^+$ | 5.47e+07 | $f_1$ |
| – | 1.79e+09 | 6.23e+07 | 1.77e+09 | 6.23e+07 | 6.16e+07 | 6.16e+07 | |
| – | 3.27e+08 | 3.80e+07 | 3.24e+08 | 3.80e+07 | 3.78e+07 | 3.78e+07 | |
| –$^-$ | 1.02e+04$^-$ | 4.45e+03$^+$ | 1.01e+04$^-$ | 4.45e+03$^+$ | 4.45e+03$^+$ | 4.45e+03 | $f_2$ |
| – | 1.02e+04 | 4.48e+03 | 1.01e+04 | 4.48e+03 | 4.47e+03 | 4.47e+03 | |
| – | 4.03e+02 | 4.00e+02 | 4.00e+02 | 4.00e+02 | 4.00e+02 | 4.00e+02 | |
| –$^-$ | 1.77e+01$^-$ | 1.17e+01$^+$ | 1.76e+01$^-$ | 1.17e+01$^+$ | 1.17e+01$^+$ | 1.17e+01 | $f_3$ |
| – | 1.78e+01 | 1.15e+01 | 1.77e+01 | 1.15e+01 | 1.15e+01 | 1.15e+01 | |
| – | 4.36e-01 | 1.07e+00 | 4.40e-01 | 1.07e+00 | 1.07e+00 | 1.07e+00 | |
| 6.21e+12$^-$ | 6.94e+12$^-$ | 8.32e+11$^-$ | 6.94e+12$^-$ | 3.31e+11$^+$ | 3.26e+11$^+$ | 3.29e+11 | $f_4$ |
| 6.48e+12 | 7.06e+12 | 9.69e+11 | 7.06e+12 | 4.10e+11 | 4.00e+11 | 4.07e+11 | |
| 3.05e+12 | 2.72e+12 | 4.74e+11 | 2.72e+12 | 2.05e+11 | 2.01e+11 | 2.04e+11 | |
| –$^-$ | 2.83e+08 | 1.03e+08$^-$ | 2.83e+08$^-$ | 1.03e+08$^+$ | 1.01e+08$^+$ | 1.03e+08 | $f_5$ |
| – | 2.91e+08 | 1.02e+08 | 2.91e+08 | 1.02e+08 | 1.01e+08 | 1.02e+08 | |
| – | 3.13e+07 | 1.40e+07 | 3.13e+07 | 1.40e+07 | 1.43e+07 | 1.40e+07 | |
| –$^-$ | 1.75e+02 | 1.68e+01$^-$ | 2.74e+04$^-$ | 1.69e+01$^-$ | 1.68e+01$^-$ | 1.33e+01 | $f_6$ |
| – | 2.26e+02 | 1.70e+01 | 4.22e+04 | 1.71e+01 | 1.70e+01 | 1.33e+01 | |
| – | 2.09e+02 | 8.34e-01 | 2.21e+04 | 7.92e-01 | 8.34e-01 | 1.08e+00 | |
| 7.69e+08$^-$ | 1.18e+09$^-$ | 2.42e+05$^+$ | 1.18e+09$^-$ | 2.43e+05$^+$ | 2.42e+05$^+$ | 2.42e+05 | $f_7$ |
| 1.00e+09 | 1.27e+09 | 2.46e+05 | 1.27e+09 | 2.46e+05 | 2.46e+05 | 2.46e+05 | |
| 8.58e+08 | 3.80e+08 | 3.19e+04 | 3.80e+08 | 3.20e+04 | 3.19e+04 | 3.19e+04 | |
| 9.18e+07$^-$ | 1.02e+08 | 2.81e+07$^-$ | 1.02e+08$^-$ | 5.18e+06$^+$ | 5.17e+06$^+$ | 5.37e+06 | $f_8$ |
| 9.69e+07 | 1.15e+08 | 4.74e+07 | 1.15e+08 | 1.39e+07 | 1.28e+07 | 1.37e+07 | |
| 3.81e+07 | 5.18e+07 | 3.67e+07 | 5.18e+07 | 2.07e+07 | 2.22e+07 | | |
| 4.03e+09$^-$ | 2.96e+10 | 1.34e+09$^+$ | 2.96e+10$^-$ | 1.34e+09$^+$ | 1.34e+09$^+$ | 1.34e+09 | $f_9$ |
| 4.00e+09 | 3.00e+10 | 1.36e+09 | 2.99e+10 | 1.36e+09 | 1.36e+09 | 1.36e+09 | |
| 6.87e+08 | 3.18e+09 | 2.79e+08 | 3.18e+09 | 2.79e+08 | 2.79e+08 | 2.79e+08 | |
| 9.13e+03$^-$ | 1.24e+04$^-$ | 7.51e+03$^+$ | 1.24e+04$^-$ | 7.53e+03$^+$ | 7.49e+03$^+$ | 7.49e+03 | $f_{10}$ |
| 9.20e+03 | 1.24e+04 | 7.44e+03 | 1.24e+04 | 7.46e+03 | 7.43e+03 | 7.44e+03 | |
| 2.50e+02 | 1.30e+02 | 2.43e+02 | 1.31e+02 | 2.44e+02 | 2.41e+02 | 2.42e+02 | |
| 2.95e+01$^-$ | 9.90e-01$^-$ | 3.17e+01$^-$ | 9.73e+01$^-$ | 2.83e+01$^+$ | 2.69e+01$^+$ | 2.66e+01 | $f_{11}$ |
| 2.95e+01 | 9.91e-01 | 3.12e+01 | 9.75e+01 | 2.79e+01 | 2.68e+01 | 2.70e+01 | |
| 2.50e+00 | 2.43e+00 | 2.40e+00 | 3.16e+00 | 2.24e+00 | 2.35e+00 | 2.24e+00 | |
| 5.95e+05$^-$ | 1.69e+06$^-$ | 3.15e+05$^+$ | 1.69e+06$^-$ | 3.15e+05$^+$ | 3.15e+05$^+$ | 3.08e+05 | $f_{12}$ |
| 5.92e+05 | 1.70e+06 | 3.08e+05 | 1.69e+06 | 3.09e+05 | 3.08e+05 | 3.08e+05 | |
| 2.67e+04 | 5.73e+04 | 5.62e+04 | 2.43e+04 | 2.43e+04 | 2.43e+04 | | |
| 1.43e+09$^-$ | 7.29e+08$^-$ | 2.36e+05$^-$ | 6.83e+08$^-$ | 2.10e+05$^+$ | 2.10e+05$^+$ | 1.73e+05 | $f_{13}$ |
| 1.54e+09 | 7.49e+08 | 2.41e+05 | 7.08e+08 | 2.12e+05 | 2.12e+05 | 2.12e+05 | |
| 6.65e+08 | 2.04e+08 | 1.86e+04 | 2.00e+08 | 1.73e+04 | 1.73e+04 | 1.73e+04 | |
| 2.47e+08$^-$ | 5.80e+09$^-$ | 2.50e+08$^-$ | 5.54e+09$^-$ | 2.46e+08$^+$ | 2.33e+08$^+$ | 2.36e+08 | $f_{14}$ |
| 2.47e+08 | 5.90e+09 | 2.49e+08 | 5.67e+09 | 2.46e+08 | 2.31e+08 | 2.38e+08 | |
| 2.08e+07 | 6.59e+08 | 2.12e+07 | 6.80e+08 | 2.08e+07 | 1.85e+07 | 1.91e+07 | |
| 6.18e+03$^-$ | 1.01e+04 | 6.20e+03$^+$ | 1.00e+04$^-$ | 6.18e+03$^+$ | 6.09e+03$^+$ | 6.10e+03 | $f_{15}$ |
| 6.33e+03 | 1.01e+04 | 6.34e+03 | 1.00e+04 | 6.33e+03 | 6.27e+03 | 6.29e+03 | |
| 4.30e+02 | 8.77e+01 | 4.34e+02 | 8.69e+01 | 4.31e+02 | 4.16e+02 | 4.27e+02 | |
| 2.39e+01$^-$ | 2.74e+02 | 2.82e+01$^-$ | 2.72e+02$^-$ | 2.15e+01$^+$ | 2.09e+01$^+$ | 2.08e+01 | $f_{16}$ |
| 2.27e+01 | 2.75e+02 | 2.83e+01 | 2.73e+02 | 2.11e+01 | 2.06e+01 | 2.08e+01 | |
| 3.35e+00 | 4.53e+00 | 3.34e+00 | 4.64e+00 | 2.76e+00 | 2.84e+00 | 2.81e+00 | |
| 9.90e+04$^-$ | 1.46e+06 | 9.97e+04$^-$ | 1.43e+06$^-$ | 9.89e+04$^+$ | 9.01e+04$^+$ | 9.26e+04 | $f_{17}$ |
| 9.88e+04 | 1.51e+06 | 9.98e+04 | 1.49e+06 | 9.87e+04 | 8.97e+04 | 9.26e+04 | |
| 6.74e+03 | 1.84e+05 | 6.83e+03 | 1.88e+05 | 6.74e+03 | 6.42e+03 | 6.52e+03 | |
| 1.78e+11$^-$ | 1.73e+10$^-$ | 2.11e+08$^-$ | 1.43e+10$^-$ | 7.36e+04$^+$ | 3.92e+04$^+$ | 4.89e+04 | $f_{18}$ |
| 1.75e+11 | 1.74e+10 | 2.09e+08 | 1.42e+10 | 7.54e+04 | 3.73e+04 | 4.91e+04 | |
| 2.26e+10 | 1.67e+09 | 4.01e+07 | 1.23e+09 | 1.06e+04 | 1.08e+04 | 1.36e+04 | |
| 3.11e+12$^-$ | 1.99e+11$^-$ | 2.03e+11$^-$ | 1.97e+11$^-$ | 1.51e+10$^+$ | 1.27e+10$^+$ | 1.31e+10 | $f_{20}$ |
| 3.02e+12 | 2.02e+11 | 2.07e+11 | 2.00e+11 | 1.65e+10 | 1.34e+10 | 1.41e+10 | |
| 2.72e+12 | 2.52e+10 | 2.58e+10 | 2.50e+10 | 5.86e+09 | 4.71e+09 | 4.94e+09 | |
| 2.13e+06$^-$ | 5.11e+06 | 2.14e+06$^-$ | 5.06e+06 | 2.14e+06$^+$ | 2.13e+06$^+$ | 2.13e+06 | $f_{19}$ |
| 2.16e+06 | 5.12e+06 | 2.16e+06 | 5.08e+06 | 2.16e+06 | 2.15e+06 | 2.15e+06 | |
| 1.09e+05 | 3.43e+05 | 1.09e+05 | 3.37e+05 | 1.09e+05 | 1.09e+05 | 1.09e+05 | |

**CEC'2013 suite**

| Type | Fun. | FDG | ID | RDG | DG2 | FII | GDG | DG |
|---|---|---|---|---|---|---|---|---|
| $T_1$ | $f_1$ | 8.11e+07 | 8.11e+07$^-$ | 8.15e+07$^-$ | 2.36e+09$^-$ | 8.15e+07$^-$ | 2.41e+09$^-$ | –$^-$ |
| | | 9.46e+07 | 9.46e+07 | 9.58e+07 | 2.47e+09 | 9.58e+07 | 2.50e+09 | – |
| | | 5.44e+07 | 5.43e+07 | 5.52e+07 | 4.31e+08 | 5.52e+07 | 4.35e+08 | – |
| | $f_2$ | 1.51e+04 | 1.51e+04$^-$ | 1.51e+04$^-$ | 1.78e+04$^-$ | 1.51e+04$^+$ | 1.79e+04$^-$ | –$^-$ |
| | | 1.46e+04 | 1.46e+04 | 1.47e+04 | 1.80e+04 | 1.47e+04 | 1.81e+04 | – |
| | | 1.70e+03 | 1.70e+03 | 1.70e+03 | 9.29e+02 | 1.70e+03 | 9.20e+02 | – |
| | $f_3$ | 2.10e+01 | 2.10e+01$^-$ | 2.10e+01$^-$ | 2.13e+01$^-$ | 2.10e+01$^-$ | 2.13e+01$^-$ | – |
| | | 2.10e+01 | 2.10e+01 | 2.10e+01 | 2.13e+01 | 2.10e+01 | 2.13e+01 | – |
| | | 1.01e-02 | 1.01e-02 | 1.07e-02 | 1.14e-02 | 1.09e-02 | 8.52e-03 | – |
| $T_{2-1}$ | $f_4$ | 6.29e+09 | 6.29e+09$^-$ | 6.29e+09$^-$ | 1.01e+11$^-$ | 1.97e+10$^-$ | 1.01e+11$^-$ | 2.22e+11$^-$ |
| | | 7.34e+09 | 7.34e+09 | 7.34e+09 | 1.15e+11 | 2.12e+10 | 1.15e+11 | 2.47e+11 |
| | | 2.61e+09 | 2.61e+09 | 2.61e+09 | 4.44e+10 | 8.66e+09 | 4.44e+10 | 1.44e+11 |
| | $f_5$ | 4.11e+06 | 4.11e+06$^-$ | 4.23e+06$^-$ | 9.99e+06$^-$ | 4.88e+06$^-$ | 9.60e+06$^-$ | 1.10e+07$^-$ |
| | | 3.85e+06 | 3.85e+06 | 3.93e+06 | 9.98e+06 | 4.58e+06 | 9.57e+06 | 1.08e+07 |
| | | 1.13e+06 | 1.13e+06 | 1.17e+06 | 5.97e+05 | 1.14e+06 | 3.95e+05 | 6.19e+05 |
| | $f_6$ | 1.06e+06 | 1.06e+06$^-$ | 1.06e+06$^+$ | 1.07e+06$^-$ | 1.06e+06$^-$ | 1.06e+06$^-$ | 1.07e+06$^-$ |
| | | 1.06e+06 | 1.06e+06 | 1.06e+06 | 1.07e+06 | 1.06e+06 | 1.06e+06 | 1.07e+06 |
| | | 1.27e+03 | 1.55e+03 | 1.39e+03 | 1.67e+03 | 1.30e+03 | 1.39e+03 | |
| | $f_7$ | 4.15e+07 | 4.15e+07$^-$ | 3.93e+08$^-$ | 1.47e+09$^-$ | 9.88e+07$^-$ | 1.89e+09$^-$ | 1.07e+09$^-$ |
| | | 2.26e+08 | 2.26e+08 | 6.94e+08 | 1.76e+09 | 1.55e+08 | 2.01e+09 | 1.02e+09 |
| | | 3.62e+08 | 3.62e+08 | 8.81e+08 | 5.84e+08 | 1.85e+08 | 7.07e+08 | 4.24e+08 |
| $T_{2-2}$ | $f_8$ | 5.54e+14 | 1.21e+15$^-$ | 7.05e+14$^-$ | 7.39e+16$^-$ | 2.39e+14$^+$ | 7.33e+16$^-$ | 3.11e+15$^-$ |
| | | 5.45e+14 | 1.23e+15 | 8.35e+14 | 7.20e+16 | 2.71e+14 | 8.19e+16 | 5.42e+15 |
| | | 4.30e+14 | 6.76e+14 | 4.66e+14 | 2.17e+16 | 1.45e+14 | 3.56e+16 | 5.47e+15 |
| | $f_9$ | 5.53e+08 | 5.53e+08$^-$ | 5.53e+08$^-$ | 9.88e+08$^-$ | 5.53e+08$^-$ | 9.78e+08$^-$ | 5.53e+08$^-$ |
| | | 5.53e+08 | 5.53e+08 | 5.53e+08 | 9.79e+08 | 5.53e+08 | 9.74e+08 | 5.53e+08 |
| | | 4.95e+07 | 4.95e+07 | 4.95e+07 | 7.17e+07 | 4.95e+07 | 7.79e+07 | 4.95e+07 |
| | $f_{11}$ | 9.48e+07 | 9.48e+07$^-$ | 9.49e+07$^+$ | 9.52e+07$^-$ | 9.47e+07$^+$ | 9.55e+07$^-$ | 9.49e+07$^+$ |
| | | 9.47e+07 | 9.47e+07 | 9.48e+07 | 9.52e+07 | 9.47e+07 | 9.54e+07 | 9.49e+07 |
| | | 2.62e+05 | 2.62e+05 | 2.73e+05 | 2.72e+05 | 2.81e+05 | 3.42e+05 | 2.80e+05 |
| | $f_{17}$ | 1.17e+10 | 1.16e+10$^-$ | 1.73e+10$^-$ | 5.86e+11$^-$ | 9.08e+09$^+$ | 1.33e+11$^-$ | 9.55e+10$^-$ |
| | | 5.66e+10 | 5.65e+10 | 1.96e+10 | 7.28e+11 | 9.56e+09 | 1.30e+11 | 1.42e+11 |
| | | 8.72e+10 | 8.72e+10 | 1.41e+10 | 3.58e+11 | 5.33e+09 | 6.80e+10 | 1.21e+11 |
| $T_{3-1}$ | $f_{12}$ | 1.49e+10 | 1.35e+10$^-$ | 1.62e+10$^-$ | 2.10e+11$^-$ | 2.17e+11$^-$ | 2.12e+11$^-$ | 3.21e+12$^-$ |
| | | 1.53e+10 | 1.44e+10 | 1.71e+10 | 2.11e+11 | 2.17e+11 | 2.12e+11 | 3.22e+12 |
| | | 4.63e+09 | 4.47e+09 | 5.05e+09 | 2.28e+10 | 2.32e+10 | 2.29e+10 | 2.84e+11 |
| | $f_{13}$ | 5.32e+09 | 5.29e+09$^-$ | 9.57e+09$^-$ | 1.30e+10$^-$ | 5.32e+09$^-$ | 1.31e+10$^-$ | 2.92e+10$^-$ |
| | | 5.41e+09 | 5.37e+09 | 1.07e+10 | 1.28e+10 | 5.40e+09 | 1.28e+10 | 2.75e+10 |
| | | 1.34e+09 | 1.33e+09 | 2.76e+09 | 2.70e+09 | 1.34e+09 | 2.71e+09 | 7.27e+09 |
| | $f_{14}$ | 5.99e+10 | 5.91e+10$^-$ | 6.03e+10$^+$ | 1.43e+11$^-$ | 6.00e+10$^+$ | 1.44e+11$^-$ | 1.19e+11$^-$ |
| | | 5.34e+10 | 5.29e+10 | 5.38e+10 | 1.51e+11 | 5.34e+10 | 1.52e+11 | 1.21e+11 |
| | | 2.06e+10 | 2.05e+10 | 2.07e+10 | 4.65e+10 | 2.06e+10 | 4.68e+10 | 3.66e+10 |
| $T_{3-2}$ | $f_{15}$ | 1.76e+07 | 1.76e+07$^-$ | 1.77e+07$^-$ | 3.29e+08$^-$ | 1.77e+07$^+$ | 3.42e+08$^-$ | 1.76e+07$^+$ |
| | | 1.86e+07 | 1.86e+07 | 1.87e+07 | 3.97e+08 | 1.87e+07 | 4.14e+08 | 1.86e+07 |
| | | 4.29e+06 | 4.29e+06 | 4.34e+06 | 2.21e+08 | 4.32e+06 | 2.33e+08 | 4.31e+06 |

TABLE S-IV

The median, the mean, and the standard deviation of the best solutions obtained by FDG on each benchmark function in the CEC'2010 and the CEC'2013 suites over 25 independent runs with each run allowing a maximum number of $1.2\times10^5$ FEs. FDG is embedded into the DECC framework and is compared with DG, GDG, FII, DG2, RDG, and ID. Some algorithms do not support the optimization on part of functions because they spend all the FEs in the decomposition, and the corresponding entries are marked with "–". The superscripts "$-$," "$\approx$," and "$+$" indicate that the corresponding result is worse than, similar to, and better than that of FDG, respectively.

**CEC'2010 suite**

| DG | GDG | FII | DG2 | RDG | ID | FDG | Fun. |
|---|---|---|---|---|---|---|---|
| – | – | 1.40e+09$^-$ | – | 1.40e+09$^-$ | 1.35e+09$^-$ | 1.35e+09 | $f_1$ |
| – | – | 1.38e+09 | – | 1.38e+09 | 1.33e+09 | 1.33e+09 | |
| – | – | 2.87e+08 | – | 2.87e+08 | 2.85e+08 | 2.84e+08 | |
| – | – | 9.47e+03$^-$ | – | 9.47e+03$^-$ | 9.28e+03$^-$ | 9.28e+03 | $f_2$ |
| – | – | 9.38e+03 | – | 9.38e+03 | 9.25e+03 | 9.25e+03 | |
| – | – | 4.74e+02 | – | 4.74e+02 | 5.05e+02 | 5.03e+02 | |
| – | – | 1.70e+01$^-$ | – | 1.71e+01$^-$ | 1.68e+01$^-$ | 1.68e+01 | $f_3$ |
| – | – | 1.71e+01 | – | 1.72e+01 | 1.70e+01 | 1.70e+01 | |
| – | – | 4.93e-01 | – | 4.89e-01 | 4.80e-01 | 4.81e-01 | |
| 1.37e+14$^-$ | – | 9.86e+12$^-$ | – | 5.03e+12$^-$ | 5.02e+12$^-$ | 5.52e+12 | $f_4$ |
| 1.38e+14 | – | 1.00e+13 | – | 5.39e+12 | 5.31e+12 | 5.31e+12 | |
| 4.10e+13 | – | 2.73e+12 | – | 2.18e+12 | 2.14e+12 | 2.14e+12 | |
| – | – | 2.62e+08$^-$ | – | 2.62e+08$^-$ | 2.62e+08$^-$ | 2.62e+08 | $f_5$ |
| – | – | 2.67e+08 | – | 2.67e+08 | 2.67e+08 | 2.67e+08 | |
| – | – | 3.22e+07 | – | 3.21e+07 | 3.22e+07 | 3.22e+07 | |
| – | – | 2.93e+01$^-$ | – | 8.25e+03$^-$ | 2.92e+01$^-$ | 7.92e+00 | $f_6$ |
| – | – | 3.56e+01 | – | 1.24e+04 | 3.55e+01 | 1.29e+05 | |
| – | – | 1.80e+01 | – | 1.33e+04 | 1.80e+01 | 1.42e+05 | |
| 2.94e+10$^-$ | – | 5.59e+08$^-$ | – | 5.90e+08$^-$ | 5.59e+08$^-$ | 5.59e+08 | $f_7$ |
| 3.41e+10 | – | 5.77e+08 | – | 6.02e+08 | 5.77e+08 | 5.77e+08 | |
| 1.82e+10 | – | 2.22e+08 | – | 2.24e+08 | 2.22e+08 | 2.22e+08 | |
| 5.13e+10$^-$ | – | 1.53e+08$^-$ | – | 9.84e+07$^-$ | 9.53e+07$^-$ | 9.53e+07 | $f_8$ |
| 7.23e+10 | – | 1.98e+08 | – | 1.00e+08 | 9.89e+07 | 9.89e+07 | |
| 7.70e+10 | – | 1.61e+08 | – | 4.90e+07 | 4.77e+07 | 4.77e+07 | |
| – | – | 1.45e+10$^-$ | – | 2.94e+10$^-$ | 1.44e+10$^-$ | 1.45e+10 | $f_9$ |
| – | – | 1.47e+10 | – | 2.97e+10 | 1.46e+10 | 1.46e+10 | |
| – | – | 1.59e+09 | – | 3.17e+09 | 1.59e+09 | 1.59e+09 | |
| – | – | 1.17e+04$^-$ | – | 1.23e+04$^-$ | 1.17e+04$^-$ | 1.17e+04 | $f_{10}$ |
| – | – | 1.17e+04 | – | 1.24e+04 | 1.16e+04 | 1.17e+04 | |
| – | – | 1.48e+02 | – | 1.29e+02 | 1.53e+02 | 1.48e+02 | |
| – | – | 8.16e+01$^-$ | – | 8.76e+01$^-$ | 7.59e+01$^-$ | 7.84e+01 | $f_{11}$ |
| – | – | 8.14e+01 | – | 8.77e+01 | 7.63e+01 | 7.81e+01 | |
| – | – | 2.70e+00 | – | 2.18e+00 | 1.96e+00 | 2.58e+00 | |
| – | – | 1.32e+06$^-$ | – | 1.63e+06$^-$ | 1.29e+06$^-$ | 1.31e+06 | $f_{12}$ |
| – | – | 1.34e+06 | – | 1.63e+06 | 1.30e+06 | 1.32e+06 | |
| – | – | 4.01e+04 | – | 4.69e+04 | 4.09e+04 | 1.42e+04 | |
| 1.19e+12$^-$ | – | 1.76e+12$^-$ | – | 1.92e+09$^-$ | 8.69e+07$^-$ | 1.04e+08 | $f_{13}$ |
| 1.19e+12 | – | 1.76e+12 | – | 1.89e+09 | 9.51e+07 | 1.11e+08 | |
| 1.50e+11 | – | 1.58e+11 | – | 3.49e+08 | 2.21e+07 | 2.25e+07 | |
| 5.56e+09$^-$ | – | 7.04e+09$^-$ | – | 5.54e+09 | 4.06e+09$^-$ | 4.57e+09 | $f_{14}$ |
| 5.76e+09 | – | 6.91e+09 | – | 5.68e+09 | 4.11e+09 | 4.80e+09 | |
| 6.78e+08 | – | 8.23e+08 | – | 6.77e+08 | 5.20e+08 | 5.88e+08 | |
| 1.00e+04$^-$ | – | 1.02e+04$^-$ | – | 1.00e+04$^-$ | 9.72e+03$^-$ | 9.81e+03 | $f_{15}$ |
| 1.00e+04 | – | 1.02e+04 | – | 1.00e+04 | 9.74e+03 | 9.83e+03 | |
| 9.23e+01 | – | 8.22e+01 | – | 8.69e+01 | 1.04e+02 | 1.01e+02 | |
| 2.80e+02$^-$ | – | 3.24e+02$^-$ | – | 2.73e+02$^-$ | 2.37e+02$^-$ | 2.49e+02 | $f_{16}$ |
| 2.81e+02 | – | 3.23e+02 | – | 2.74e+02 | 2.38e+02 | 2.49e+02 | |
| 5.39e+00 | – | 6.83e+00 | – | 4.48e+00 | 3.78e+00 | 4.18e+00 | |
| 1.45e+06$^-$ | – | 1.55e+06$^-$ | – | 1.43e+06$^-$ | 1.15e+06$^-$ | 1.23e+06 | $f_{17}$ |
| 1.50e+06 | – | 1.60e+06 | – | 1.49e+06 | 1.21e+06 | 1.29e+06 | |
| 1.87e+05 | – | 1.79e+05 | – | 1.87e+05 | 1.59e+05 | 1.71e+05 | |
| 3.21e+12$^-$ | – | 2.05e+12$^-$ | – | 9.85e+09$^-$ | 1.21e+09 | 1.21e+10 | $f_{18}$ |
| 3.17e+12 | – | 2.08e+12 | – | 9.81e+09 | | 1.22e+10 | |
| 2.80e+11 | – | 2.15e+11 | – | 1.02e+09 | | 1.09e+09 | |
| – | – | 2.86e+11$^-$ | – | 1.61e+11$^-$ | 1.91e+11 | | $f_{20}$ |
| – | – | 2.87e+11 | – | 1.64e+11 | 1.95e+11 | | |
| – | – | 3.20e+10 | – | 2.04e+10 | 2.45e+10 | | |
| 4.58e+06$^-$ | – | 4.62e+06$^-$ | – | 4.66e+06$^-$ | 4.55e+06$^-$ | 4.56e+06 | $f_{19}$ |
| 4.62e+06 | – | 4.66e+06 | – | 4.71e+06 | 4.58e+06 | 4.58e+06 | |
| 3.06e+05 | – | 3.13e+05 | – | 3.17e+05 | 3.02e+05 | 3.03e+05 | |

**CEC'2013 suite**

| Type | Fun. | FDG | ID | RDG | DG2 | FII | GDG | DG |
|---|---|---|---|---|---|---|---|---|
| $T_1$ | $f_1$ | 1.86e+09 | 1.86e+09$^-$ | 1.93e+09$^-$ | – | 1.93e+09$^-$ | – | – |
| | | 1.87e+09 | 1.87e+09 | 1.94e+09 | – | 1.94e+09 | – | – |
| | | 3.62e+08 | 3.61e+08 | 3.68e+08 | – | 3.68e+08 | – | – |
| | $f_2$ | 1.70e+04 | 1.70e+04$^-$ | 1.71e+04$^-$ | – | 1.71e+04$^-$ | – | – |
| | | 1.73e+04 | 1.73e+04 | 1.74e+04 | – | 1.74e+04 | – | – |
| | | 9.62e+02 | 9.63e+02 | 9.63e+02 | – | 9.63e+02 | – | – |
| | $f_3$ | 2.13e+01 | 2.13e+01$^-$ | 2.13e+01$^-$ | – | 2.13e+01$^-$ | – | – |
| | | 2.13e+01 | 2.13e+01 | 2.13e+01 | – | 2.13e+01 | – | – |
| | | 9.75e-03 | 1.03e-02 | 1.04e-02 | – | 9.52e-03 | – | – |
| $T_{2\text{-}1}$ | $f_4$ | 9.91e+10 | 9.21e+10$^-$ | 1.01e+11$^-$ | – | 7.36e+10$^-$ | – | 8.90e+11 |
| | | 1.13e+11 | 1.06e+11 | 1.14e+11 | – | 7.47e+10 | – | 9.20e+11 |
| | | 4.45e+10 | 4.44e+10 | 4.45e+10 | – | 2.19e+10 | – | 3.74e+11 |
| | $f_5$ | 9.01e+06 | 9.01e+06$^-$ | 9.01e+06$^-$ | – | 9.40e+06$^-$ | – | – |
| | | 9.16e+06 | 9.15e+06 | 9.16e+06 | – | 9.37e+06 | – | – |
| | | 4.96e+05 | 4.96e+05 | 4.96e+05 | – | 5.11e+05 | – | – |
| | $f_6$ | 1.06e+06 | 1.07e+06 | 1.07e+06$^-$ | – | 1.06e+06$^-$ | – | – |
| | | 1.06e+06 | 1.07e+06 | 1.07e+06 | – | 1.06e+06 | – | – |
| | | 1.51e+03 | 1.36e+03 | 1.18e+03 | – | 1.43e+03 | – | – |
| | $f_7$ | 1.17e+09 | 1.17e+09$^-$ | 4.70e+09$^-$ | – | 1.41e+09$^-$ | – | 4.89e+09 |
| | | 1.20e+09 | 1.20e+09 | 5.87e+09 | – | 1.72e+09 | – | 5.16e+09 |
| | | 5.12e+08 | 5.12e+08 | 3.48e+09 | – | 8.69e+08 | – | 1.84e+09 |
| $T_{2\text{-}2}$ | $f_8$ | 7.96e+15 | 1.27e+16$^-$ | 7.42e+16$^-$ | – | 3.28e+15$^-$ | – | 5.62e+16 |
| | | 8.97e+15 | 1.29e+16 | 7.53e+16 | – | 4.06e+15 | – | 5.93e+16 |
| | | 4.33e+15 | 6.11e+15 | 2.98e+16 | – | 2.70e+15 | – | 2.70e+16 |
| | $f_9$ | 9.87e+08 | 9.87e+08$^-$ | 9.88e+08$^-$ | – | 9.88e+08$^-$ | – | 9.88e+08$^-$ |
| | | 9.79e+08 | 9.79e+08 | 9.79e+08 | – | 9.79e+08 | – | 9.79e+08 |
| | | 7.17e+07 | 7.17e+07 | 7.17e+07 | – | 7.17e+07 | – | 7.18e+07 |
| | $f_{10}$ | 9.52e+07 | 9.52e+07$^-$ | 9.54e+07$^-$ | – | 9.53e+07$^-$ | – | 9.54e+07$^-$ |
| | | 9.52e+07 | 9.52e+07 | 9.53e+07 | – | 9.52e+07 | – | 9.54e+07 |
| | | 2.72e+05 | 2.72e+05 | 3.90e+05 | – | 2.57e+05 | – | 3.11e+05 |
| | $f_{11}$ | 5.86e+11 | 5.86e+11$^-$ | 1.35e+11$^-$ | – | 9.82e+10$^-$ | – | 8.84e+12 |
| | | 7.28e+11 | 7.28e+11 | 1.42e+11 | – | 1.08e+11 | – | 1.46e+13 |
| | | 3.58e+11 | 3.58e+11 | 5.69e+10 | – | 4.31e+10 | – | 1.28e+13 |
| $T_{3\text{-}1}$ | $f_{12}$ | 2.04e+11 | 1.71e+11$^-$ | 3.01e+11$^-$ | – | – | – | – |
| | | 2.06e+11 | 1.73e+11 | 2.99e+11 | – | – | – | – |
| | | 2.28e+10 | 2.13e+10 | 2.51e+10 | – | – | – | – |
| | $f_{13}$ | 1.88e+10 | 1.80e+10$^-$ | 4.84e+10$^-$ | – | 1.87e+10$^-$ | – | 3.38e+11 |
| | | 1.90e+10 | 1.82e+10 | 5.24e+10 | – | 1.89e+10 | – | 4.23e+11 |
| | | 3.59e+09 | 3.50e+09 | 1.00e+10 | – | 3.57e+09 | – | 2.61e+11 |
| | $f_{14}$ | 2.19e+11 | 2.12e+11$^-$ | 2.26e+11$^-$ | – | 2.20e+11$^-$ | – | 6.16e+12 |
| | | 2.36e+11 | 2.28e+11 | 2.44e+11 | – | 2.38e+11 | – | 6.68e+12 |
| | | 5.75e+10 | 5.69e+10 | 5.83e+10 | – | 5.76e+10 | – | 7.66e+12 |
| $T_{3\text{-}2}$ | $f_{15}$ | 1.84e+08 | 1.84e+08$^-$ | 2.16e+08$^\approx$ | – | 2.07e+08$^-$ | – | 1.93e+08$^-$ |
| | | 2.01e+08 | 2.04e+08 | 2.44e+08 | – | 2.31e+08 | – | 2.13e+08 |
| | | 1.84e+08 | 1.84e+08 | 2.16e+08 | – | 2.07e+08 | – | 1.93e+08 |

TABLE S-V

THE MEDIAN, THE MEAN, AND THE STANDARD DEVIATION OF THE BEST SOLUTIONS OBTAINED BY MA-SW-CHAINS, MOS-CEC2013, DSPLSO, DLLSO, AND FCRACC$_{\text{FDG+SHADE}}$ ON EACH BENCHMARK FUNCTION IN THE CEC'2010 AND THE CEC'2013 SUITES. THE MEDIAN OF DSPLSO IS NOT REPORTED IN ITS ORIGINAL PAPER, AND THE CORRESPONDING ENTRIES ARE MARKED WITH "−". THE SUPERSCRIPTS "−," "≈," AND "+" INDICATE THAT THE CORRESPONDING RESULT IS WORSE THAN, SIMILAR TO, AND BETTER THAN THAT OF FCRACC$_{\text{FDG+SHADE}}$, RESPECTIVELY. THE RESULTS IN BOLDFACE INDICATE THEY ARE THE BEST.

| CEC'2010 MA-SW-Chains | MOS-CEC2013 | DSPLSO | DLLSO | FCRACC$_{\text{FDG+SHADE}}$ | Fun. | Type | Fun. | FCRA$_{\text{FDG+SHADE}}$ | DLLSO | DSPLSO | MOS-CEC2013 | CEC'2013 MA-SW-Chains |
|---|---|---|---|---|---|---|---|---|---|---|---|---|
| 2.67e−14 | 0.00e+00 | − | 2.93e−22 | 0.00e+00 | | | $f_1$ | 0.00e+00 | 3.86e−22 | − | 0.00e+00 | 6.15e−13 |
| **3.80e−14**$^{+}$ | **0.00e+00**$^{+}$ | **7.73e−20**$^{+}$ | **3.13e−22**$^{+}$ | **0.00e+00** | $f_1$ | | | **0.00e+00** | 3.99e−22$^{-}$ | **1.18e−19**$^{-}$ | **0.00e+00**$^{\approx}$ | **1.14e−12**$^{-}$ |
| 4.91e−14 | 0.00e+00 | 7.07e−21 | 8.03e−23 | 0.00e+00 | | | | 0.00e+00 | 1.32e−22 | 1.06e−20 | 0.00e+00 | 1.28e−12 |
| 8.47e+02 | 0.00e+00 | − | 9.85e+02 | 2.99e+00 | | | $f_2$ | 7.45e+00 | 1.14e+03 | − | 8.24e+02 | 1.13e+03 |
| 8.40e+02$^{+}$ | **0.00e+00**$^{+}$ | 4.45e+02$^{-}$ | 9.82e+02$^{-}$ | 3.08e+00 | $f_2$ | T1 | | **7.47e+00** | 1.14e+03$^{-}$ | 1.06e+03$^{-}$ | 8.23e+02$^{-}$ | 1.18e+03$^{-}$ |
| 4.88e+01 | 0.00e+00 | 1.65e+01 | 4.39e+01 | 2.77e−01 | | | | 4.24e−01 | 5.78e+01 | 4.45e+02 | 4.69e+01 | 1.84e+02 |
| 5.16e−13 | 1.67e−12 | − | 2.89e−14 | 1.48e−01 | | | $f_3$ | 2.01e+01 | 2.16e+01 | − | 1.69e−12 | 6.79e−13 |
| **5.76e−13**$^{+}$ | **1.65e−12**$^{+}$ | **2.52e−13**$^{+}$ | **2.76e−14**$^{+}$ | 1.49e−01 | $f_3$ | | | 2.01e+01 | 2.16e+01$^{-}$ | 2.16e+01$^{-}$ | **1.69e−12**$^{+}$ | **6.78e−13**$^{+}$ |
| 2.73e−13 | 6.44e−14 | 1.89e−14 | 2.38e−15 | 3.93e−01 | | | | 1.75e−03 | 4.07e−03 | 7.53e−03 | 9.16e−14 | 2.28e−13 |
| 3.10e+11 | 1.63e+10 | − | 4.37e+11 | 1.64e+07 | | | $f_4$ | 6.13e+05 | 6.34e+09 | − | 7.80e+07 | 2.70e+09 |
| 2.97e+11$^{-}$ | 1.70e+10$^{-}$ | 4.30e+11$^{-}$ | 4.40e+11$^{-}$ | **1.83e+07** | $f_4$ | | | **5.74e+05** | 6.68e+09$^{-}$ | 9.40e+09$^{-}$ | 8.73e+07$^{-}$ | 3.80e+09$^{-}$ |
| 6.19e+10 | 6.39e+09 | 8.31e+10 | 1.10e+11 | 1.26e+07 | | | | 1.68e+05 | 1.68e+09 | 1.89e+09 | 3.11e+07 | 2.70e+09 |
| 2.30e+08 | 1.08e+08 | − | 1.19e+07 | 4.58e+07 | | | $f_5$ | 3.64e+06 | 6.60e+05 | − | 6.95e+06 | 1.98e+06 |
| 2.18e+08$^{-}$ | 1.07e+08$^{-}$ | **6.30e+06**$^{+}$ | 1.22e+07$^{+}$ | 4.70e+07 | $f_5$ | T2-1 | | 3.62e+06 | 7.00e+05$^{-}$ | **6.30e+05**$^{+}$ | 6.89e+06$^{-}$ | 2.26e+06$^{+}$ |
| 5.75e+07 | 2.49e+07 | 1.76e+06 | 3.43e+06 | 7.77e+06 | | | | 2.59e+05 | 1.28e+05 | 1.02e+05 | 9.16e+05 | 1.36e+06 |
| 2.45e+00 | 9.28e−08 | − | 4.00e−09 | 1.36e+01 | | | $f_6$ | 1.05e+06 | 1.06e+06 | − | 1.39e+05 | 6.24e+02 |
| 1.42e+05$^{-}$ | 1.11e−07$^{+}$ | **9.45e−09**$^{+}$ | 5.20e−01$^{+}$ | 1.29e+01 | $f_6$ | | | 1.05e+06 | 1.06e+06$^{-}$ | 1.06e+06$^{-}$ | 1.43e+05$^{+}$ | **1.07e+04**$^{+}$ |
| 3.96e+05 | 5.88e−08 | 1.20e−09 | 7.46e−01 | 2.39e+00 | | | | 1.90e+03 | 8.28e+02 | 8.05e+02 | 6.86e+04 | 2.09e+04 |
| 7.94e−03 | 0.00e+00 | − | 1.22e+01 | 7.34e−22 | | | $f_7$ | 2.42e+06 | 1.33e+06 | − | 1.10e+03 | 3.99e+06 |
| 1.17e+02 | **0.00e+00**$^{+}$ | 4.76e+02$^{-}$ | 7.19e+02$^{-}$ | **8.14e−22** | $f_7$ | | | 5.00e+06 | 1.60e+06$^{-}$ | 5.50e+06$^{-}$ | **4.65e+03**$^{+}$ | 3.78e+06$^{-}$ |
| 2.37e+02 | 0.00e+00 | 1.31e+02 | 2.59e+03 | 6.03e−22 | | | | 1.20e+07 | 8.38e+05 | 2.26e+06 | 1.06e+04 | 8.46e+05 |
| 2.76e+06 | 1.38e−07 | − | 2.34e+07 | 1.21e−19 | | | | | | | | |
| 6.90e+06 | 1.40e+00$^{-}$ | 3.11e+07$^{-}$ | 2.34e+07$^{-}$ | **2.66e−19** | $f_8$ | T2-1 | | | | | | |
| 1.90e+07 | 7.01e+00 | 9.36e+04 | 2.46e+05 | 3.11e−19 | | | | | | | | |
| 1.48e+07 | 3.47e+06 | − | 4.33e+07 | 8.35e+05 | | | | | | | | |
| 1.49e+07 | 3.59e+06$^{-}$ | 4.59e+07$^{-}$ | 4.36e+07$^{-}$ | **8.85e+05** | $f_9$ | | | | | | | |
| 1.61e+06 | 4.89e+05 | 3.04e+06 | 4.28e+06 | 1.48e+05 | | | | | | | | |
| 2.02e+03 | 3.82e+03 | − | 8.89e+02 | 5.83e+03 | | | | | | | | |
| 2.01e+03 | 3.81e+03$^{+}$ | 7.99e+03$^{-}$ | **8.91e+02**$^{+}$ | 6.13e+03 | $f_{10}$ | T2-2 | | | | | | |
| 1.59e+02 | 1.62e+02 | 1.28e+02 | 3.66e+01 | 7.65e+02 | | | | | | | | |
| 3.77e+01 | 1.91e+02 | − | 2.75e+00 | 8.58e+00 | | | | | | | | |
| 3.86e+01 | 1.91e+02$^{-}$ | **3.04e−12**$^{+}$ | 5.80e+00$^{+}$ | 8.74e+00 | $f_{11}$ | | | | | | | |
| 8.06e+00 | 4.01e−01 | 2.89e−13 | 5.40e+00 | 6.99e−01 | | | | | | | | |
| 3.09e−06 | 0.00e+00 | − | 1.24e+04 | 3.34e−10 | | | | | | | | |
| 3.24e−06$^{-}$ | **0.00e+00**$^{+}$ | 9.52e+04$^{-}$ | | **3.59e−10** | $f_{12}$ | | | | | | | |
| 5.78e−07 | 0.00e+00 | 6.69e+03 | 1.46e+03 | 1.93e−10 | | | | | | | | |
| 8.61e+02 | 5.69e+02 | − | 7.28e+02 | 7.97e+00 | | | | | | | | |
| 9.83e+02 | 8.23e+02$^{-}$ | 5.48e+02$^{-}$ | 7.35e+02$^{-}$ | **7.18e+00** | $f_{13}$ | | | | | | | |
| 5.66e+02 | 6.77e+02 | 1.69e+02 | 1.93e+02 | 3.99e+00 | | | | | | | | |
| 3.23e+07 | 9.77e+06 | − | 1.25e+08 | 5.87e+06 | | | $f_8$ | 2.10e+08 | 1.16e+14 | − | 2.82e+12 | 4.65e+13 |
| 3.25e+07$^{-}$ | 9.69e+06$^{-}$ | 1.60e+08$^{-}$ | 1.24e+08$^{-}$ | **5.98e+06** | $f_{14}$ | T2-2 | | **2.35e+08** | 1.20e+14$^{-}$ | 1.55e+14$^{-}$ | 2.85e+12$^{-}$ | 4.63e+13$^{-}$ |
| 2.46e+06 | 6.71e+05 | 8.50e+06 | 7.38e+06 | 5.28e+05 | | | | 7.86e+07 | 3.35e+13 | 2.96e+13 | 1.44e+12 | 9.18e+12 |
| 2.67e+03 | 7.45e+03 | − | 8.40e+02 | 6.43e+03 | | | $f_9$ | 3.00e+08 | 1.26e+08 | − | 4.18e+08 | 1.16e+08 |
| 2.68e+03$^{+}$ | 7.44e+03$^{-}$ | 9.91e+03$^{-}$ | **8.33e+02**$^{+}$ | 6.56e+03 | $f_{15}$ | | | 2.99e+08 | 1.30e+08$^{-}$ | **8.07e+07**$^{+}$ | 3.99e+08$^{-}$ | 1.14e+08$^{+}$ |
| 9.95e+01 | 1.90e+02 | 6.70e+01 | 4.31e+01 | 7.54e+02 | | | | 2.07e+05 | 2.04e+07 | 2.24e+07 | 6.26e+07 | 2.05e+07 |
| 9.32e+01 | 3.87e+02 | − | 3.70e+00 | 1.42e−13 | | | $f_{10}$ | 9.21e+07 | 9.40e+07 | − | 1.17e+06 | 3.32e+02 |
| 9.95e+01$^{-}$ | 3.79e+02$^{-}$ | **4.68e−12**$^{+}$ | 4.25e+00$^{-}$ | **1.40e−13** | $f_{16}$ | | | 9.22e+07 | 9.40e+07$^{-}$ | 9.39e+07$^{-}$ | 9.38e+05$^{+}$ | **3.66e+04**$^{+}$ |
| 1.53e+01 | 1.83e+01 | 4.49e−13 | 2.41e+00 | 2.31e−15 | | | | 2.07e+05 | 2.11e+05 | 2.26e+05 | 4.79e+05 | 6.17e+04 |
| 1.28e+00 | 2.67e−07 | − | 9.02e+04 | 5.06e−03 | | | $f_{11}$ | 7.31e+08 | 9.29e+11 | − | 1.71e+07 | 2.10e+08 |
| 1.27e+00 | **2.73e−07**$^{+}$ | 6.84e+05$^{-}$ | 9.05e+04$^{-}$ | 5.57e−03 | $f_{17}$ | T3-1 | | 7.07e+08 | 9.30e+11$^{-}$ | 9.27e+11$^{-}$ | **1.73e+07**$^{+}$ | 2.10e+08$^{+}$ |
| 1.24e−01 | 7.67e−08 | 3.63e+04 | 3.53e+03 | 1.38e−03 | | | | 2.44e+08 | 9.50e+09 | 9.48e+09 | 5.04e+06 | 2.43e+07 |
| 1.41e+03 | 1.55e+03 | − | 2.29e+03 | 3.78e+02 | | | | | | | | |
| 1.57e+03 | 1.77e+03$^{-}$ | 1.35e+03$^{-}$ | 2.55e+03$^{-}$ | **3.85e+02** | $f_{18}$ | | | | | | | |
| 6.73e+02 | 9.57e+02 | 3.87e+02 | 8.32e+02 | 5.32e+01 | | | | | | | | |
| 1.04e+03 | 1.81e+02 | − | 1.85e+03 | 5.74e+04 | | | $f_{12}$ | 1.80e+05 | 1.79e+03 | − | 1.56e+01 | 1.24e+03 |
| 1.06e+03$^{+}$ | **2.93e+02**$^{+}$ | 1.06e+03$^{-}$ | 1.88e+03$^{-}$ | 2.82e+05 | $f_{20}$ | | | 2.83e+06 | 1.79e+03$^{+}$ | 1.05e+03$^{+}$ | **8.13e+01**$^{+}$ | 1.23e+03$^{+}$ |
| 9.38e+01 | 3.99e+02 | 1.79e+02 | 1.90e+02 | 5.16e+05 | | | | 1.12e+07 | 1.39e+02 | 5.37e+01 | 1.57e+02 | 8.32e+01 |
| | | | | | | | $f_{13}$ | 1.14e+07 | 2.70e+08 | − | 1.02e+06 | 1.91e+07 |
| | | | | | | T3-1 | | 2.34e+07 | 3.35e+08$^{-}$ | 1.20e+09$^{-}$ | **1.00e+06**$^{+}$ | 1.98e+07$^{+}$ |
| | | | | | | | | 5.34e+07 | 1.71e+08 | 4.99e+08 | 5.53e+05 | 2.30e+06 |
| | | | | | | | $f_{14}$ | 3.26e+07 | 1.03e+08 | − | 1.28e+07 | 1.47e+08 |
| | | | | | | | | 5.30e+07 | 1.72e+08$^{-}$ | 8.31e+09$^{-}$ | **1.24e+07**$^{+}$ | 1.45e+08$^{-}$ |
| | | | | | | | | 9.10e+07 | 1.38e+08 | 6.67e+09 | 2.86e+06 | 1.69e+07 |
| 2.95e+05 | 2.87e+04 | − | 1.79e+06 | 2.79e+05 | | | $f_{15}$ | 1.07e+06 | 4.45e+06 | − | 1.68e+06 | 5.76e+06 |
| 2.95e+05$^{-}$ | **2.92e+04**$^{+}$ | 8.20e+06$^{-}$ | 1.80e+06$^{-}$ | 2.85e+05 | $f_{19}$ | T3-2 | | **1.11e+06** | 4.48e+06$^{-}$ | 4.13e+07$^{-}$ | 1.71e+06$^{-}$ | 5.90e+06$^{-}$ |
| 1.34e+04 | 2.29e+03 | 4.69e+05 | 9.96e+04 | 2.74e+04 | | | | 1.37e+06 | 3.32e+05 | 3.11e+06 | 1.44e+05 | 1.36e+06 |

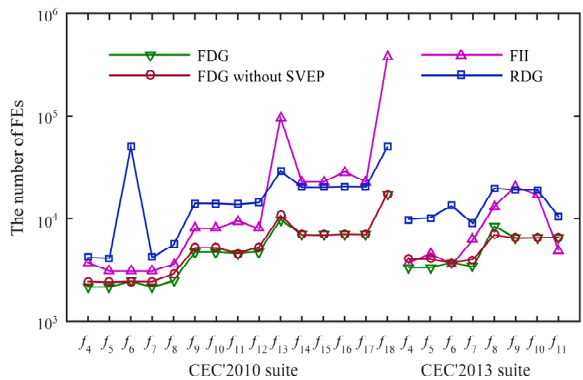

Fig. S-4. The variation of the number of FEs consumed by FDG on partially separable functions if SVEP is removed.

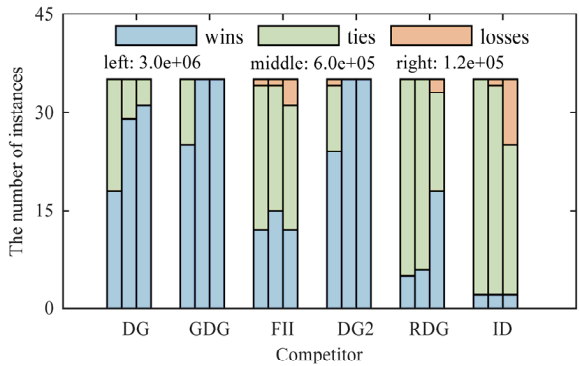

Fig. S-5. The numbers of wins, ties, and losses of FDG against DG, GDG, FII, DG2, RDG, and ID on a total of 35 benchmark functions under the cases of three different quantities of FEs.

From Tables S-II, S-III, S-IV, and Fig. S-5, it can be seen that, with the reduction of the number of available FEs, the advantage of FDG becomes more obvious. When $6.0 \times 10^5$ FEs are available, FDG performs no worse than all of its six competitors except for being surpassed by ID and FII on $f_{18}$ in the CEC'2010 suite and $f_8$ in the CEC'2013 suite, respectively; the numbers of the functions on which it outperforms DG, GDG, FII, DG2, and RDG increase to 29, 35, 15, 35, and 6, respectively. When the number of allowed FEs reduces to $1.2 \times 10^5$, FDG outperforms DG and RDG on more functions. Besides, it is interesting to observe that, in this case, FDG is not only surpassed by ID on more functions, but also is dominated by FII and RDG on some new functions such as $f_6$ in the CEC'2010 suite and $f_{11}$ in the CEC'2013 suite. Compared with ID, the failure of FDG on corresponding functions can be attributed to its FE consumption during the decomposition process. By investigating the functions on which FDG loses its dominance against FII or RDG, it can be revealed that it is the fine-grained decomposition result of FDG that causes its failure. The underlying reason mainly consists in that, for each of the corresponding functions, FDG generates more subfunctions than its competitors, and each subfunction is assigned with fewer FEs, which cannot support the optimizer to find a good enough solution.

*Remark:* Although the termination condition with a huge number of FEs is commonly suggested for the sake of algorithm research for LSBO, many real-world LSBO problems cannot afford so many FEs because their solution quality generally has to be measured by expensive simulations due to the absence of analytical objective functions [2], [3]. From this perspective, it is significant to reduce the FE requirement of a decomposition algorithm. On the other side, it is still a challenging task to judge which decomposition is better between two candidates, since the quality of a decomposition has some relation with the available FEs and a more accurate decomposition may lead to an inferior solution in practice. Therefore, although FDG makes a great progress in decomposition accuracy and efficiency, there is still a long way to develop more practical decomposition algorithms.

## References

[1] K. Tang, X. Li, P. N. Suganthan, Z. Yang, and T. Weise, "Benchmark functions for the CEC'2010 special session and competition on large scale global optimization," Nat. Inspired Comput. Lab., Univ. Sci. Technol. China, Hefei, China, Rep., 2010.

[2] P. Yang, K. Tang, and X. Yao, "Turning high-dimensional optimization into computationally expensive optimization," *IEEE Trans. Evol. Comput.*, vol. 22, no. 1: 143–156, Feb. 2018.

[3] Z. Ren, B. Pang, M. Wang, Z. Feng, Y. Liang, and Y. Zhang, "Surrogate model assisted cooperative coevolution for large scale optimization," *Appl. Intell.*, vol. 49, no. 2: 513–531, Feb. 2019.